# On Young Tableau Involutions and Patterns in Permutations

Erik Ouchterlony

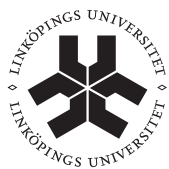

**Linköpings universitet**
**INSTITUTE OF TECHNOLOGY**

ii

# Contents

iv

# Abstract

This thesis deals with three different aspects of the combinatorics of permutations. In the first two papers, two flavours of pattern avoiding permutations are examined; and in the third paper Young tableaux, which are closely related to permutations via representation theory, are studied.

A permutation is alternating if it is alternatingly rising and falling and doubly alternating if both itself and its inverse are alternating. In the first paper we give solutions to several interesting problems regarding pattern avoiding doubly alternating permutations, such as finding a bijection between 1234-avoiding permutations and 1234-avoiding doubly alternating permutations of twice the size.

The matrix representation of a permutation is a square 0-1-matrix such that every row and column has exactly one 1. A partial permutation is a rectangular matrix which is a submatrix of a permutation matrix. In the second paper partial permutations which can be extended to pattern avoiding permutations are examined. A general algorithm is presented which is subsequently used to solve many different problems.

The third paper deals with involutions on Young tableaux. There is a surprisingly large collection of relations among these involutions and in this paper we make an effort to study them systematically in order to create a coherent theory. The most interesting result is that for Littlewood-Richardson tableaux, $a^3 = I$, where $a$ is the composition of three different involutions: the fundamental symmetry map, the reversal and rotation involutions.

# Acknowledgements

First of all, I thank my supervisor, Svante Linusson, who made it possible for me to start working in combinatorics. His support and his enthusiasm for the subject have been a great inspiration. I also thank Johan Wästlund, for standing in as supervisor when Svante was on parental leave, and Bruce Sagan, my opponent, for thoroughly examining the thesis and giving numerous suggestions for improvements.

For almost a year I was in France to work on my thesis, supported by the European network "Algebraic Combinatorics in Europe," three months in Lyon and six months in Bordeaux. I am truly grateful for getting this opportunity and I would like to thank all my colleagues there, Christian Krattenthaler, Johan Thapper, Olivier Guibert, Mireille Bousquet-Mélou, Mark Dukes, Theresia Eisenkölbl, Xavier Viennot, Philippe Duchon and many others. It was a great time, not only mathematically!

Finally, I would like to thank everyone at the Department of Mathematics, my family and friends. Thank you!

vi

# 1   Introduction

## Combinatorial objects and generating functions

Let $\mathcal{F}$ be a set of *combinatorial objects* (or *discrete objects*), on which a *size function*, $f : \mathcal{F} \to \{0, 1, \ldots\}$, is defined. The size function classifies the objects into finite sets of objects of the same size, $\mathcal{F}_n = \{a \in \mathcal{F} : f(a) = n\}$. Let $f_n = |\mathcal{F}_n|$ be the cardinality of these sets. The most common task in enumerative combinatorics is to calculate the integer sequence $f_n$.

A powerful tool for counting combinatorial objects is the use of (ordinary) *generating functions*, where the numbers $f_n$ are used as the coefficients in a *formal power series*,

$$F(x) = \sum_{n \geqslant 0} f_n x^n = \sum_{a \in \mathcal{F}} x^{f(a)}.$$

Much more information on generating functions and other subjects in combinatorics can be found in the books by Flajolet and Sedgewick [8], Stanley [27] and Wilf [32].

## Permutations

The type of combinatorial objects we concentrate on in this thesis is the *permutations*, which are linear orderings on a finite set of objects. Since it is not important what kind of objects we are permuting, we often use the set $\{1, 2, \ldots, n\}$ and with this set, the permutations are bijections, $\{1, 2, \ldots, n\} \to \{1, 2, \ldots, n\}$. The size of a permutation is the size of the set it acts on, and the number of permutations of size $n$ is $n!$. The set of permutations of size $n$, $\mathcal{S}_n$, is also known as the *symmetric group* of degree $n$, due to the properties we discuss in Section 3.

Even though permutations look like fairly simple objects, a surprisingly rich theory has been created around them. For the combinatorial aspects of permutations, see the book by Bóna [4] and for a more algebraic treatment, see Sagan's book [22].

A common way to represent a permutations is the *word form*, where we write, e.g., $\sigma = 1342$ for the permutation $\sigma \in \mathcal{S}_4$, which satisfies $\sigma(1) = 1$, $\sigma(2) = 3$, $\sigma(3) = 4$ and $\sigma(4) = 2$. Another representation is to use a 0-1-matrix, which has a 1 as the element in row $i$ and column $j$ iff $\sigma(i) = j$. Note that there is choice of how to orient the matrix, and a different convention has been used by many authors, but the one used here has the big advantage that matrix multiplication corresponds to composition of permutations. The third representation we are going to use is to regard the permutation as a perfect matching on the complete bipartite graph with nodes $\{1, 2, \ldots, n\}$ on the left side and $\{1', 2', \ldots, n'\}$ on the right side, where node $i$ and $j'$ are matched iff $\sigma(i) = j$. See Figure 1 for an example.

$$\sigma = 251634$$ 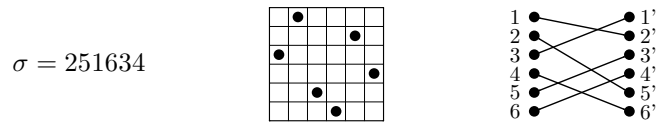

Figure 1: A permutation in word, matrix and perfect matching form.

## 2  Pattern avoidance

We can define pattern avoidance both using the word form and the matrix form of the permutation.

**Definition 2.1 (Pattern avoidance).** A permutation, $\sigma = \sigma_1 \sigma_2 \cdots \sigma_n$, is said to *contain* the *pattern* $\tau = a_1 a_2 \cdots a_k$ if there is a subsequence of $\sigma$ which is order equivalent to $\tau$. In other words, there is an integer sequence $1 \leqslant m_1 < m_2 < \cdots < m_k \leqslant n$, such that $\sigma_{m_i} < \sigma_{m_j}$ iff $a_i < a_j$ for all $i, j \in \{1, \ldots, k\}$.

Alternatively, $\sigma$ contains the pattern $\tau$ if $\tau$ is a submatrix of $\sigma$, i.e., if there is a way to remove $n - k$ rows and columns of $\sigma$ in order to get $\tau$. We say that $\sigma$ *avoids* $\tau$ if it does not contain $\tau$ and $\mathcal{S}_n(\tau_1, \tau_2, \ldots, \tau_t)$ denotes the set of permutations avoiding all the patterns, $\tau_i$, $i = 1, \ldots, t$. If $\tau = \{\tau_1, \tau_2, \ldots, \tau_t\}$ is a set of patterns, we use the shorthand notation $\mathcal{S}_n(\tau)$ for $\mathcal{S}_n(\tau_1, \tau_2, \ldots, \tau_t)$.

The problem of enumerating pattern avoiding permutations was first studied by Knuth [13] when examining stack sorting. It turns out that permutations which are sortable by a single stack is exactly those that avoid the pattern 231. A classical result in pattern avoidance is that the number of permutations avoiding a pattern of length three is the same for all such patterns.

**Theorem 2.2.** $|\mathcal{S}_n(\tau)| = C_n$ *for every pattern, $\tau$, of length three.*

The numbers, $C_n = \frac{1}{n+1}\binom{2n}{n} = 1, 1, 2, 5, 14, 42, 132, \ldots$, are called the *Catalan numbers* and are very important in combinatorics, since there is a huge number of objects counted by them, and there are, e.g., 66 different listed in [28]. As a curiosity, it can be mentioned that Catalan himself called the numbers *Segner numbers* since Segner had earlier found them when counting triangulations of convex polygons [25].

As seen from Theorem 2.2, there are many sets of patterns that are equally easy to avoid and one of the most fundamental problems in pattern avoidance is to identify them.

**Definition 2.3.** Two sets of patterns, $\tau$ and $\tau'$, are *Wilf-equivalent* if $|\mathcal{S}_n(\tau)| = |\mathcal{S}_n(\tau')|$ for all $n \geqslant 0$. Wilf-equivalence divides the pattern sets into equivalence classes, which are called *Wilf-classes*.

A large class of Wilf-equivalences can be found by using the symmetries of the permutations. If $\sigma = \sigma_1\sigma_2\cdots\sigma_n$, is a permutation, we define the *reverse* of $\sigma$ to be $\sigma^r = \sigma_n\sigma_{n-1}\cdots\sigma_1$ and the *complement* of $\sigma$ to be $\sigma^c = (n+1-\sigma_1)(n+1-\sigma_2)\cdots(n+1-\sigma_n)$. These operations correspond to reflecting the permutation matrix vertically and horizontally, respectively.

**Lemma 2.4.** *Let $\tau = \{\tau_1, \ldots \tau_t\}$ be a set of patterns. Then $\tau$, $\tau^r$, $\tau^c$ and $\tau^{-1}$ are all Wilf-equivalent.*

A more unexpected set of Wilf-equivalences is a generalisation of the fact the $12\cdots k$ and $k(k-1)\cdots 1$ are Wilf-equivalent. It turns out that if we append both of these patterns to the same permutation of $\{k+1, k+2, \ldots, n\}$, they still are Wilf-equivalent.

**Theorem 2.5 (Backelin, West and Xin [2]).** *The patterns $12\cdots k\tau$ and $k(k-1)\cdots 1\tau$ are Wilf-equivalent for all permutations, $\tau$, of $\{k+1, k+2, \ldots, n\}$.*

For the singleton patterns of length four, Lemma 2.4 and Theorem 2.5 imply that there can be at most four different Wilf-classes, represented by the patterns 1234, 1324, 1342 and 2413. But the last two are in fact also Wilf-equivalent.

**Theorem 2.6 (Stankova [26]).** *The patterns 1342 and 2413 are Wilf-equivalent.*

The remaining three Wilf-classes are indeed different, and we have

$$|\mathcal{S}_n(1234)| = 1, 1, 2, 6, 23, 103, 513, 2761, 15767, 94359, 586590, \ldots$$
$$|\mathcal{S}_n(1324)| = 1, 1, 2, 6, 23, 103, 513, 2762, 15793, 94776, 591950, \ldots$$
$$|\mathcal{S}_n(1342)| = 1, 1, 2, 6, 23, 103, 512, 2740, 15485, 91245, 555662, \ldots$$

The problem of enumerating the permutations in the first and third class have been solved, and we know explicit formulas, but the pattern 1324 has, however, turned out to be surprisingly hard and the best result so far is a rather complicated recursive formula by Marinov and Radoičić [20].

## Generalisations

One way to generalise pattern avoidance is to look at a larger class of restrictions. This was done by Babson and Steingrímsson [1], who studied restrictions where some of the letters in the pattern are required to be consecutive. This generalisation allow us, for example, to examine the number of permutations having no pattern 123 where the elements representing 1 and 2 are consecutive. Using the notation for *generalised patterns*, we say that such permutations are 12-3-avoiding. Here, unlike for normal pattern, the dash means that any distance is allowed and the absence of dash means these values need to be consecutive. For example, we have that 24153 contains 12-3, but 31524 avoids it. This case and many others where studied by Claesson [6], and he showed that they are counted by the *Bell numbers*, which also count the number of partitions of a set of size $n$.

Another way to make a generalisation is to relax the notation of permutation. In [17] matrices that are the upper-left corners of permutations where studied. This means that there no longer needs to be exactly one 1 in every row and column, since there may be empty rows and columns. We call these matrices *pre-permutations* and a pre-permutation is said to be *extendably avoiding* the pattern set $\tau = \{\tau_1, \ldots, \tau_t\}$ if it is an upper-left corner of a $\tau$-avoiding permutation. We denote by $\mathcal{S}_{d,c,r}(\tau)$, the set of all extendably $\tau$-avoiding pre-permutations with $c$ empty columns, $r$ empty rows and $d$ ones (or *dots*). Two pattern sets, $\tau$ and $\tau'$ are *extendably Wilf-equivalent* if $|\mathcal{S}_{d,c,r}(\tau)| = |\mathcal{S}_{d,c,r}(\tau')|$ for all nonnegative integers, $d$, $c$ and $r$.

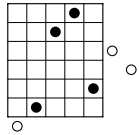

Figure 2: An example of an extendably 312-avoiding pre-permutation. The hollow circles indicate one of the two possible ways of extending it into a 312-avoiding permutation.

**Theorem 2.7 (Linusson).** *For $\tau = 132, 231, 312$ and $321$,*

$$|\mathcal{S}_{d,c,r}(\tau)| = \binom{n+d}{d} - \binom{n+d}{d-1}$$

*and for $\tau = 123$ and $213$,*

$$|\mathcal{S}_{d,c,r}(\tau)| = \sum_{i=0}^{d} \binom{n-c}{i}\binom{n-r}{i} - \binom{n+d}{d-1}.$$

The last case follows from a more general theorem, which is the extended version of Theorem 2.5.

**Theorem 2.8 (Linusson).** *The patterns $12\cdots k\tau$ and $k(k-1)\cdots 1\tau$ are extendably Wilf-equivalent for all permutations, $\tau$, of $\{k+1, k+2, \ldots, n\}$.*

In the second paper we examine extended pattern avoidance for many different pattern sets. We give generating functions for all different combinations of patterns of size three, a bijective proof of Theorem 2.7, and many identities for pattern-pairs, where the patterns are of size three and four, respectively.

In order to achieve this we have developed a method, called the *generating graph method*, which is a refined version of the generating tree method developed by Chung et al [5] and used in many research papers, see for example [31].

## Specialisations

Instead of generalising pattern avoidance, we can look at special types of permutations. One common specialisation is to study pattern avoiding involutions, i.e., permutations which are equal to its own inverse. Another type of permutations to study is the *alternating permutations*, for which the permutation is alternatingly rising and descending, starting with a rise. For example 27483615 is an alternating permutation. A collection of many of the results on the specialisations (and also generalisations) can be found in the survey paper by Kitaev and Mansour [12].

In the first paper we are studying *doubly alternating permutation*, for which both the permutation and its inverse are alternating. These have not received much attention and the only result so far being found in a paper by Guibert and Linusson [10], where doubly alternating *Baxter permutations* are studied. In terms of generalised patterns, a Baxter permutation is a permutation avoiding the pattern set $\{3\text{-}14\text{-}2, 2\text{-}41\text{-}3\}$.

**Theorem 2.9 (Guibert, Linusson).** *Doubly alternating Baxter permutations are counted by the Catalan numbers.*

We present an alternative proof of this theorem by using results on doubly alternating 2413-avoiding permutations. Furthermore, we examine the 1234-avoiding doubly alternating permutations and demonstrate a bijection, which uses the RSK correspondence defined in Section 3, in order to prove the following theorem:

**Theorem 2.10.** $\quad |\operatorname{DA}_{2n}(1234)| = |\mathcal{S}_n(1234)|.$

Here $\operatorname{DA}_n(\tau)$ denotes the set of doubly alternating permutations avoiding $\tau$. We have also been able to bijectively show a doubly alternating version of Theorem 2.5, for $k = 2$.

**Theorem 2.11.** *Let $\tau$ be any permutation of $\{3, 4, \ldots, n\}$, $m \geqslant 3$. Then*

$$|\operatorname{DA}_n(21\tau)| = |\operatorname{DA}_n(12\tau)|.$$

There is also a large collection of conjectures involving doubly alternating permutations presented in the paper, which indicates that work remains to be done on this subject.

# 3    Young tableaux

A *partition* of a positive integer $n$ is a way of writing $n$ as a sum of positive integers, where the terms are called *parts*. The order of the parts is not important so we write them in weakly decreasing order. Partitions are often displayed as *Young diagrams*, where each row corresponds to a part in the partition, see Figure 3.

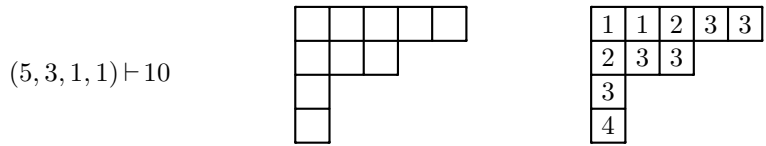

Figure 3: Examples of a partition, its Young diagram and a Young tableau of partition shape.

A *skew Young diagram*, $\lambda/\mu$, is the set theoretic difference between two Young diagrams, $\lambda$ and $\mu$, such that $\mu \subseteq \lambda$.

**Definition 3.1.** A skew Young diagram, $\lambda/\mu$, which is filled with positive integers is a *Young tableau* of shape $\lambda/\mu$ if every row is weakly increasing and every column is strictly increasing. If $\mu$ is the zero partition, we say that the tableau has a *partition shape*.

Another way of representing how a skew Young diagram is filled is to use a *recording matrix*, which is a matrix, such that the value in row $i$ and column $j$ is the number of $j$:s in row $i$ of the tableau. For example the Young tableau in Figure 3 has the recording matrix

$$\begin{pmatrix} 2 & 1 & 2 & 0 \\ 0 & 1 & 2 & 0 \\ 0 & 0 & 1 & 0 \\ 0 & 0 & 0 & 1 \end{pmatrix}.$$

We can see that for a Young tableau with partition shape, the recording matrix must be upper triangular. For a more extensive exposition on Young tableaux and their applications, we recommend the books by Fulton [9], Sagan [22] and Stanley [28].

## Littlewood-Richardson tableaux

Next we define a special type of tableau, which plays a crucial role in the theory of symmetric functions, as we will see in the next section. There are actually several different ways to define it, but here we use a definition using the reading word of a Young tableau.

**Definition 3.2.** A *word* is a sequence of positive integers and the *reading word* of a Young tableau is the word we get by starting at the top right corner and reading toward the left, then continuing with the row below and reading it from right to left and so on.

**Definition 3.3.** Let $w = w_1 w_2 \ldots w_n$ be a word. Then the *weight* of $w$ is $\mathsf{weight}(w) = (|\{s : w_s = i\}|)_{i=1}^{N}$, i.e., a vector which counts the number of times each value appears in $w$. We say that a vector is a partition if it is weakly decreasing. If $w$ is a reading word of a Young tableau $T$, then the weight of $T$ is $\mathsf{weight}(T) = \mathsf{weight}(w)$.

**Definition 3.4.** A word $w = w_1 w_2 \ldots w_n$ is *Yamanouchi* if the weight for every prefix, $w_1 w_2 \ldots w_k$, is a partition.

**Definition 3.5.** A Young tableau is a *Littlewood-Richardson tableau* if its reading word is Yamanouchi.

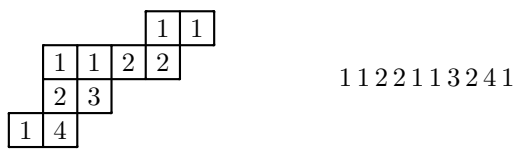

$$1\,1\,2\,2\,1\,1\,3\,2\,4\,1$$

Figure 4: Example of a Littlewood-Richardson tableau which has the shape $(6, 5, 3, 2)/(4, 1, 1)$ and weight $(5, 3, 1, 1)$. Its reading word is shown to the right.

We can see that the recording matrix of a Littlewood-Richardson tableau must be lower triangular. For example the in Figure 4 has the recording matrix

$$\begin{pmatrix} 2 & 0 & 0 & 0 \\ 2 & 2 & 0 & 0 \\ 0 & 1 & 1 & 0 \\ 1 & 0 & 0 & 1 \end{pmatrix}.$$

In particular, the first row of a Littlewood-Richardson tableau can only consist of ones, since otherwise the reading word would not be Yamanouchi, as it does not begin with a one.

## Symmetric functions and Schur functions

An introduction to symmetric functions and Schur functions can be found in [28] and for further study, see Macdonald's book [19]. Here we will just give a small sample of the theory in order to motivate why the Littlewood-Richardson tableau are so important.

Let $\mathcal{X} = x_1, x_2, \ldots$ be a finite or infinite *alphabet* of variables. A *symmetric function* is a formal power series (i.e. a possibly infinite polynomial) in $\mathcal{X}$, such that any permutation of the variables leaves it unchanged. For example, $x_1 x_2 + x_1 x_3 + x_2 x_3$ is a symmetric function in three variables. If $|\mathcal{X}| = n$ is finite we can let $\mathcal{S}_n$ act on the polynomials in $\mathcal{X}$ by permuting the variables, so a symmetric function is a function which is unchanged under the action of the symmetric group. The theory of symmetric functions is fundamental to algebraic combinatorics, and has wide spectrum of applications, such as group theory, Lie theory and algebraic geometry.

The symmetric functions are often expressed as a basis expansion, and the most fundamental basis is the set of monomial symmetric functions, $\{m_\lambda : \lambda \text{ is a partition}\}$, where

$$m_\lambda = m_\lambda(x) = \sum_a x_{i_1}^{a_1} x_{i_2}^{a_2} \cdots x_{i_N}^{a_N}$$

and the sum is over all distinct monomials having exponents $\lambda_1, \lambda_2, \ldots, \lambda_N$. For example we have $m_{(1)} = x_1 + x_2 + \cdots$ and $m_{(2,1)} = \sum_{i \neq j} x_i^2 x_j = x_1^2 x_2 + x_1 x_2^2 + x_1^2 x_3 + x_1 x_3^2 + x_2^2 x_3 + x_2 x_3^2 + \cdots$. Note that if the number of variables in the alphabet is infinite, the sum also becomes infinite, but we treat it as a formal power series, so that there is no problem with convergence.

Using this basis, we can create another basis, the *complete symmetric functions*, $h_\lambda$, defined by

$$h_\lambda = h_{(\lambda_1)} h_{(\lambda_2)} \cdots h_{(\lambda_N)},$$

where

$$h_{(n)} = \sum_{\nu \vdash n} m_\nu.$$

We have, for example, $h_{(1)} = m_{(1)}$, $h_{(1,1)} = m_{(1)} m_{(1)} = 2 m_{(1,1)} + m_{(2)}$ and $h_{(2)} = m_{(1,1)} + m_{(2)} = \sum_{i \geqslant j} x_i x_j$. These two bases are used to define a scalar product, $\langle\,,\,\rangle$, on the vector space of symmetric functions, such that $m_\lambda$ and $h_\lambda$ are dual, i.e.,

$$\langle m_\lambda, h_\mu \rangle = [\![ \lambda = \mu ]\!],$$

where we have used Iverson's bracket notation [11] for the characteristic function, $[\![ S ]\!] = 1$ if $S$ is true and 0 otherwise.

The set of *Schur functions*, is a basis for the symmetric functions which is orthonormal with respect to this scalar product. This makes them very attractive to work with and they have been used in a large number of applications. One way to define them is to use the Young tableaux discussed above:

$$s_\lambda = s_\lambda(x) = \sum_T x^{\mathsf{weight}(T)},$$

where the sum is over Young tableaux of shape $\lambda$. Here we use the shorthand notation $x^a = x_1^{a_1} x_2^{a_2} \cdots x_N^{a_N}$.

The Schur functions have many remarkable properties, like the following formula, known as the Cauchy identity, which can be proved using the RSK correspondence, defined later in this section.

**Theorem 3.6.**
$$\prod_{i,j} \frac{1}{(1 - x_i y_j)} = \sum_\lambda s_\lambda(x) s_\lambda(y).$$

From this it can be shown that $\langle s_\lambda, s_\mu \rangle = [\![\lambda = \mu]\!]$, which proves that the Schur functions constitute an orthonormal basis. Therefore any product of the Schur functions must be expressible in terms of the Schur functions themselves, hence

$$s_\mu s_\nu = \sum_\lambda c_{\mu\nu}^\lambda s_\lambda,$$

where the $c_{\mu\nu}^\lambda$ are called the *Littlewood-Richardson coefficients*. The reason for the name is that Littlewood and Richardson [18] conjectured that the coefficients could be found by counting the Littlewood-Richardson tableaux. Later Schützenberger [24] proved this conjecture and in more recent years many short and elegant proofs have been published [16, 21, 29, 33].

**Proposition 3.7.** *The Littlewood-Richardson coefficients $c_{\mu\nu}^\lambda$ are equal to the number of Littlewood-Richardson tableaux of shape $\lambda/\mu$ and with weight $\nu$.*

Since the product of symmetric functions is commutative, $c_{\mu\nu}^\lambda$ must be symmetric in $\mu$ and $\nu$. This, however, is not at all clear by looking at the Littlewood-Richardson tableaux. This is the kind of symmetry we are examining in the third paper. In fact there are many other identities, and moreover, these identities, which can often be express by involutions, are often connected and give rise to many commutation relations.

We next discuss two classical tools for finding such relations, jeu de taquin and the Robinson-Schensted-Knuth correspondence.

## Jeu de taquin

One of the most useful operations on Young tableaux is Schützenberger's *jeu de taquin*, which has got its name from the French for the 15-puzzle. The procedure is as follows: Suppose that $T$ is a Young tableau of shape $\lambda/\mu$, where $\mu$ is not zero so we don't have a partition shape, and pick a box belonging to a southeast corner of the Young diagram of $\mu$. This box will be used as the blank of the 15-puzzle, and we will do slides such that the blank moves to the right or downwards until it is on the southeast side of the tableau. By requiring that the Young tableau conditions (i.e., rows are weakly increasing and columns strictly increasing) are never broken, we are left with a unique slide for each move — the smallest available element is always moved — so the path for the blank is unique. Repeating the procedure with another box on the inside as the blank, we finally get a tableau of partition shape, which we denote jdt($T$). See Figure 5 for an example. It is a remarkable property of jeu de taquin that jdt($T$) is well defined.

**Lemma 3.8.** *The operation* jdt *does not depend on the order in which the boxes on the inside are chosen.*

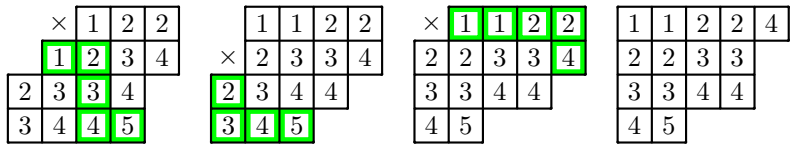

Figure 5: Examples of jeu de taquin slides, starting with $T$ and ending with jdt($T$).

Another operation which can be defined with the help of jeu de taquin is the *Schützenberger involution*, $\xi$, which is a function between two Young tableaux of the same partition shape but where the weight is reversed. This time we use a filled box as the blank, namely the box in the northwest corner. After the sliding is finished, we fill the blank with the value $N+1-v$, where $N$ is the maximal value of the Young tableau and $v$ is the value in the box we used as blank. This newly filled box, however, should not participate in the remainder of the procedure, so we draw it as being outside of the tableau. We then start over and pick a new box to use as the blank and continue until all the boxes have been used as the blank exactly once. See Figure 6 for an example.

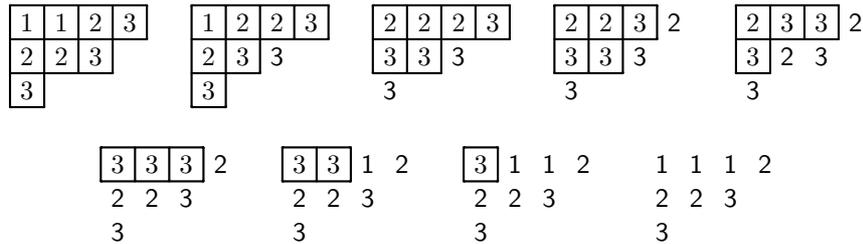

Figure 6: Example of the Schützenberger involution.

## The RSK correspondence

The Robinson-Schensted correspondence is a bijection between permutations and pairs of standard Young tableaux, i.e. Young tableau for which the numbers $1, \ldots, n$ occurs exactly once, of the same partition shape. It is defined through a bumping process, where each element of the permutation is inserted one at the time and left-to-right into a tableau, $P$. For each element, $m$, we start with the top row and if $m$ is larger than all the entries in the tableau, we insert $m$ at the end of the row. Otherwise there is a unique entry, $m'$, in the first row of the tableau, which can be replaced by $m$ without disrupting the increasing order of the row, namely the least element in the row greater than $m$. We say that $m$ bumped down $m'$. Next we repeat the procedure in the row below, but with $m'$ instead of $m$. The tableau, $P$, created this way is called the *insertion tableau*. The second tableau, $Q$, is the *recording tableau*, which is a tableau of the same shape as $P$, and where the entries tell in which order the boxes in $P$ were created. See Figure 7 for an example.

To see that this is a bijection, we can construct the inversion of the bumping process, where we take each box of the insertion tableau in the order determined by the recording tableau and bump it upwards. This process is uniquely determined by the requirement that the rows are increasing and one can check that this inverted bumping process recreates the original permutation.

Figure 7: Example of the bumping process for calculating the RSK correspondence of the permutation 53412.

The Robinson-Schensted correspondence was invented by Robinson [7], independently discovered by Schensted [23] and generalised by Knuth [14] into a bijection between square matrices with nonnegative integer entries and pairs of (not necessarily standard) Young tableaux of the same shape. This correspondence has been given the name RSK from the initials of the three inventors. For the definition of the RSK correspondence, see Section 4 of Paper 3 in this thesis.

The two fundamental operations we have looked at here, jeu de taquin and RSK, might seem like they are completely different in character. However, they are closely related by the following theorem.

**Theorem 3.9.** *Let $M$ be a recording matrix of a Young tableaux, $T$, and $\mathrm{RSK}(M^r) = (P, Q)$, where $M^r$ denotes the matrix $M$ turned upside down.*

*Then*

$$P = \mathrm{jdt}(T) \qquad and \qquad Q = \xi(\mathrm{jdt}(T')),$$

*where $T'$ is a tableau which has recording matrix $M^t$.*

This theorem, which in its first incarnation was proved by Schützenberger [24], is one of the reasons why Young tableaux are so interesting to study. There are many, often surprising, relationships and identities involving tableaux, and there are a large number of different ways of representing them which highlight different properties, e.g., by pictures [30], BZ-triangles [3] and honeycombs [15].

A classical theorem by Knuth [14] reveals how the RSK correspondence is affected by transposing and rotating the matrix.

**Theorem 3.10 (Knuth).** *Let $M$ be a matrix and $\mathrm{RSK}(M) = (P, Q)$. Then $\mathrm{RSK}(M^t) = (Q, P)$ and $\mathrm{RSK}(M^\bullet) = (\xi(P), \xi(Q))$, where $M^\bullet$ is the matrix $M$ rotated 180 degrees.*

In the third paper we study four involutions on Young tableaux, the *companion involution* and the *rotation involution*, which correspond to transposing and rotating the recording matrix, respectively; the *reversal*, which is a generalisation of the Schützenberger involution; and the *first fundamental symmetry map*, which is a manifestation of the symmetry of the Littlewood-Richardson coefficients. The final result is a commutative diagram, which summarises the different commutation relations that exist among these involutions.

**Theorem 3.11.** *The diagram in Figure 21 of Paper 3 is commutative.*

# Paper 1

**16**

# PATTERN AVOIDING DOUBLY ALTERNATING PERMUTATIONS

ERIK OUCHTERLONY

ABSTRACT. We study pattern avoiding doubly alternating (DA) permutations, i.e., alternating (or zigzag) permutations whose inverse is also alternating. We exhibit a bijection between the 1234-avoiding permutations and the 1234-avoiding DA permutations of twice the size using the Robinson-Schensted correspondence. Further, we present a bijection between the 1234- and 2134-avoiding DA permutations and we prove that the 2413-avoiding DA permutations are counted by the Catalan numbers.

## 1. Introduction

A permutation $\sigma \in \mathcal{S}_n$ is said to *contain* the pattern $\tau \in \mathcal{S}_m$ if there is a subsequence of (the word representation of) $\sigma$ which is order equivalent to (the word representation of) $\tau$. To distinguish between patterns and other permutations, we will use slightly different notation. For example, the permutation $(1, 3, 2, 4)$ will be written as 1324 if it is used as a pattern. We will often use the matrix representation of $\sigma$, which is the $n \times n$ 0-1-matrix having ones in the positions with matrix coordinates $(i, \sigma(i))$. It can also be written as $(\llbracket \sigma(i) = j \rrbracket)_{i,j=1}^n$, using Iverson's bracket notation [9] for the characteristic function, $\llbracket S \rrbracket = 1$ if $S$ is true and 0 otherwise. In the figures we will use dots instead of ones and leave the zeroes empty, as in Figure 1, to make the picture clearer. In this notation $\sigma$ contains the pattern $\tau$ if some submatrix of (the matrix representation of) $\sigma$ is equal to (the matrix representation of) $\tau$. The permutations not containing $\tau$ are called $\tau$-*avoiding*, and we write

$$\mathcal{S}_n(\tau_1, \tau_2, \ldots, \tau_t) = \{\sigma \in \mathcal{S}_n : \sigma \text{ is } \tau_i\text{-avoiding for all } i = 1, \ldots, t\}.$$

A word, e.g., a permutation, $w = (w_i)_{i=1}^n$, is *(up-down)-alternating* if $w_{2i-1} < w_{2i}$ and $w_{2i} \geqslant w_{2i+1}$ for all applicable $i$. This mean that the word alternates between *rises* and *descents*, beginning with a rise. If it instead starts with a descent, it is called *down-up-alternating*.

*Date*: December 17, 2005.

*Key words and phrases.* Permutation, Pattern avoidance, Doubly alternating, Bijection.

The author was supported by the European Commission's IHRP Programme, grant HPRN-CT-2001-00272, "Algebraic Combinatorics in Europe".

A permutation $\sigma$ is *doubly alternating* (DA) if both $\sigma$ and $\sigma^{-1}$ are alternating. The set of pattern avoiding doubly alternating permutation is denoted by

$$\mathrm{DA}_n(\tau_1, \ldots, \tau_t) = \{\sigma \in \mathcal{S}_n(\tau_1, \ldots, \tau_t) : \sigma \text{ is doubly alternating }\}.$$

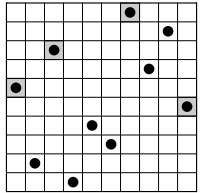

FIGURE 1. The permutation $(7, 9, 3, 8, 1, 10, 5, 6, 2, 4) \in \mathrm{DA}(1234)$ contains the pattern 3214, but avoids 1234.

Pattern avoiding permutations have been subject to much attention since the pioneering work by Knuth's [10], where he used them for studying stack sortable permutations. For a thorough summary of the current status of research, see Bóna's book [4]. Alternating permutations have a long history, they where studied already in the 19th century by André [1], and it is well know that they are counted by the tangent and secant numbers, also known as Euler numbers, $E_k$, and thus, their exponential generating function is $\tan(x)+\sec(x)$. Alternating permutation avoiding patterns have been studied by Mansour [11], but there are still many open questions remaining.

The doubly alternating permutations where first counted by Foulkes [6], up to $n = 10$, using Theorem 5.2. The only formula known is due to Stanley [13]:

$$\sum_{n \text{ odd}} \mathrm{DA}_n x^n = \sum_{k \text{ odd}} E_k^2 ((\log(1+x)/(1-x))/2)^k/k!,$$

$$\sum_{n \text{ even}} \mathrm{DA}_n x^n = (1-x^2)^{-1/2} \sum_{k \text{ even}} E_k^2 ((\log(1+x)/(1-x))/2)^k/k!,$$

from which we get that the first few numbers for $\mathrm{DA}_n, n \geqslant 1$, are

$$1, 1, 1, 2, 3, 8, 19, 64, 880, 3717, 18288, 92935, \ldots$$

The motivation for studying doubly alternating permutations came from work by Guibert and Linusson [8] who showed that doubly alternating Baxter permutations are counted by the Catalan numbers. It was a natural step to study other restrictions to see whether interesting results could be found.

Using computer simulations Guibert came up with several conjectures that indicated there are surprising connections between doubly alternating permutations and ordinary permutations. Some of these are proved in this paper, see proposition 4.3 and Theorems 5.5 and 6.6, whereas other still remain unproved and are listed in conjecture 7.1.

In this paper we study doubly alternating permutations avoiding patterns of lengths three and four. The patterns of length three are covered in Section 3. In Section 4, we show that doubly alternating permutations avoiding 2413 are counted by Catalan numbers, and are closely related to the doubly alternating Baxter permutations. Section 5 contains a bijection between $\mathrm{DA}_{2n}(1234)$ and $\mathcal{S}_n(1234)$ and in Section 6 we use a result by Babson and West [2] to construct a bijection between $\mathrm{DA}_n(12\tau)$ and $\mathrm{DA}_n(21\tau)$, where $\tau$ is any permutation of $\{3, 4, \ldots, m\}$, $m \geqslant 3$. In Section 7 other patterns giving the same sequence are investigated and in the final section some remarks on a few DA permutations avoiding two patterns of length four are given.

I like to thank Olivier Guibert for introducing me to the problem and for interesting discussions. Thanks also to Svante Linusson, Bruce Sagan and Mark Dukes for numerous comments and suggestions.

## 2. **Notation and basic facts**

First we define the *reverse*, the *complement* and the *rotation* of a permutation $\sigma$,

$$\sigma^r = (\sigma(n+1-i))_{i=1}^n$$
$$\sigma^c = (n+1-\sigma(i))_{i=1}^n$$
$$\sigma^\# = (\sigma^c)^r = (\sigma^r)^c = (n+1-\sigma(n+1-i))_{i=1}^n$$

In terms of matrices, the first two correspond to flipping the matrix vertically and horizontally, respectively, whereas the last operation rotates the matrix 180 degrees. However, these bijections do not in general preserve the doubly alternating property, which means that we lose some symmetry compare with ordinary permutations, so that more genuinely different patterns need to be examined. However, it is obvious from the definition that inverting and, if $n$ is even, rotating a permutation does preserve the property of being doubly alternating.

**Lemma 2.1.**
  (a) $\sigma \in \mathrm{DA}_n \iff \sigma^{-1} \in \mathrm{DA}_n$
  (b) $\sigma \in \mathrm{DA}_{2n} \iff \sigma^\# \in \mathrm{DA}_{2n}$

Another simple, but very useful, property that follows from the DA condition is that some areas on the border of the matrix can never have a dot, see Figure 2.

**Lemma 2.2.**
  (a) *Let* $\sigma \in \mathrm{DA}_{2n}$, *then*
    (i) $\sigma(1)$ *is odd,*
    (ii) $\sigma(2) \in \{3, 5, 7, \ldots, 2n-1, 2n\}$,
    (iii) $\sigma(2) = 2n$ *iff* $\sigma(1) = 2n-1$,
    (iv) $\sigma(2n)$ *is even,*
    (v) $\sigma(2n-1) \in \{1, 2, 4, 6, \ldots, 2n-2\}$,
    (vi) $\sigma(2n-1) = 1$ *iff* $\sigma(2n) = 2$.
  (b) *Let* $\sigma \in \mathrm{DA}_{2n+1}$, *then*
    (i) $\sigma(1)$ *is odd and less than* $2n+1$,
    (ii) $\sigma(2)$ *is odd and greater than* $1$,
    (iii) $\sigma(2n+1)$ *is even,*

    (iv) $\sigma(2n) \in \{4, 6, 8, \ldots, 2n, 2n+1\}$,
    (v) $\sigma(2n) = 2n+1$ *iff* $\sigma(2n+1) = 2n$.

*Proof.* First for the case a(i), if $\sigma(1) = k > 1$, then $\sigma^{-1}(k) = 1$, so $\sigma^{-1}(k-1) > 1$, which implies that $k$ is odd, since $\sigma \in$ DA. For a(ii), assume $\sigma(2) = m < 2n$, $m > \sigma(1) \geqslant 1$. Then $\sigma^{-1}(m-1) > 2$ or $\sigma^{-1}(m+1) > 2$, so $m$ is odd. The equivalence a(iii) is a direct consequence of the definition of doubly alternating, since $\sigma(2) > \sigma(1)$ and $\sigma^{-1}(2n-1) < \sigma^{-1}(2n)$. The other cases are similar. □

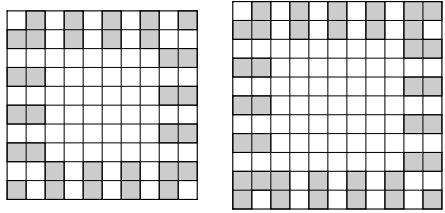

FIGURE 2. Illustration of Lemma 2.2. Shaded areas are forbidden.

    Note that this lemma could be applied to $\sigma^{-1}$ as well, because of Lemma 2.1. From the two lemmas it is clear that there is a difference between odd and even sizes, so they require separate treatment in many of the proofs. This disparity is also reflected in the fact that $DA_n(\tau_1, \tau_2, \ldots, \tau_t)$ is not increasing for all patterns. Some counterexamples are, $DA_4(321) = 2 > 1 = DA_5(321)$, $DA_6(2431) = 6 > 5 = DA_7(2431)$ and, if Conjecture 8.2 is true, $DA_{27}(1234, 2134) = 2681223 > 2674440 = DA_{28}(1234, 2134)$.

## 3. Patterns of length three

    For normal permutations, patterns of length three are the first non-trivial cases; they are all counted by the Catalan numbers. However, for the doubly alternating permutations, it turns out that all the patterns of length three are (more or less) trivial.

**Proposition 3.1.**
    (i) $|DA_n(123)| = |DA_n(213)| = |DA_n(231)| = |DA_n(312)| = 1$
    (ii) $|DA_n(132)| = [\![n \text{ even or } n = 1]\!]$
    (iii) $|DA_n(321)| = 1 + [\![n \text{ even and } n \geqslant 4]\!]$

*Proof.* First it is easy to see that the proposition holds for $n \leqslant 3$, since $(1, 2)$ and $(1, 3, 2)$ are the single DA permutations of size two and three, respectively, so we can assume that $n \geqslant 4$ for the rest of the proof.

    **Pattern** $123$**,** $n$ **even:** Assume $\sigma$ avoids $123$. We want to show that $\sigma(1) = n-1$, which, by Lemma 2.2a (iii), implies $\sigma(2) = n$, and, by induction,

that $\sigma = (n-1, n, n-3, n-2, \ldots, 3, 4, 1, 2)$. But $\sigma(1) < n-1$ is not allowed, since it gives us, $\sigma(1) < \sigma(2) < n$, so that $\sigma$ has the pattern 123.

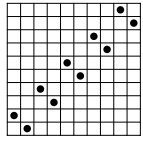

**Pattern 123, $n$ odd:** Firstly, $\sigma(2) = n$, since if $\sigma(2) < n$, we have $\sigma(1) < \sigma(2) < n$. Also $\sigma(1) = n-2$, since otherwise $1 < \sigma^{-1}(n-2) < \sigma^{-1}(n-1)$. By symmetry $\sigma(n-2) = 1$ and $\sigma(n) = 2$.

We want to show that $\sigma(2i-1) = n-2i$ for all $2i-1 = 1, 3, \ldots, n-2$. To avoid 123, it is clear that $1 = \sigma(n-2) < \sigma(n-4) < \cdots < \sigma(3) < \sigma(1) = n-2$, so it suffices to show that $\sigma(2i-1)$ is odd for all $2i-1 = 3, 5, \ldots, n-4$. But $\sigma(2j-1) = 2k$ implies that $\sigma^{-1}(2k-1) < \sigma^{-1}(2k) < 2j$, where $2k-1 < 2k < \sigma(2j)$, so $\sigma$ has pattern 123. Since the even rows must also be descending, we get only one option, namely $\sigma = (n-2, n, n-4, n-1, n-6, n-3, \ldots, 5, 8, 3, 6, 1, 4, 2)$

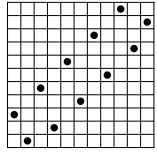

**Pattern 132 or 213, $n$ even:** Similar arguments as to those above also hold in these cases, since if $\sigma \in \mathrm{DA}_n(132)$ and $\sigma(1) = a < n-1$, then $\sigma(1) = a < a+1 < \sigma(2)$ and $1 < 2 < \sigma^{-1}(a+1)$ gives a contradicting 132-pattern. By rotational symmetry, we have the same result for $\mathrm{DA}_n(213)$.

**Pattern 213, $n$ odd:** First, $\sigma(2) = n$, as otherwise $\sigma(3) < \sigma(2) < n$, where $2 < 3 < \sigma^{-1}(n)$, which is a 213-pattern. By symmetry $\sigma(n) = 2$. Furthermore, $\sigma(4) = n-1$, since otherwise $\sigma(5) < \sigma(4) < n-1$ and $4 < 5 < \sigma^{-1}(n-1)$, since $\sigma^{-1}(n-1)$ is even. Next, $\sigma^{-1}(n-2) = 3$, since it must be less than $\sigma^{-1}(n-1) = 4$ and if $\sigma^{-1}(n-2) = 1$ we get the pattern 213 in rows 1, 3 and 4. Continuing, by using induction, we get $\sigma(2i) = n-2i+3$ and $\sigma(2i-1) = n-2i+2$ for $2i = 4, 6, \ldots, n-1$. We are now forced to put $\sigma(1) = 1$, and it is clear that by doing this we get $\sigma \in \mathrm{DA}_n(213)$.

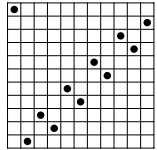

**Pattern 231:** If $\sigma(1) > 1$, then $\sigma(2) > \sigma(1) > 1$, so we get the pattern we were supposed to avoid. Hence $\sigma(1) = 1$. If $n$ is even, by the same argument, $\sigma(n) = n$. Now let $\tilde{\sigma}$ be the permutation where the first row and column, and, if $n$ is even, also the last row and column, are removed.

Now $\sigma \in \mathrm{DA}_n(231)$ iff $(\tilde{\sigma}^c)^{-1}$, i.e. $\tilde{\sigma}$ rotated $90°$ counterclockwise, is in $\mathrm{DA}_{2k}(213)$, hence $\sigma = (1, 3, 2, 5, 4, \ldots, n-1, n-2, n)$, if $n$ is even and $\sigma = (1, 3, 2, 5, 4, \ldots, n, n-1)$, if $n$ is odd.

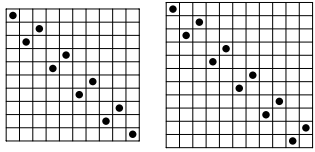

**Pattern** 312: By symmetry, this is the same as for 231.

**Pattern** 321: If $\sigma(1) = 1$ and $\sigma(n) = n$ we can rehash the proof for the pattern 231. But here the case $\sigma(1) \neq 1$ is not excluded. We can assume $\sigma(1) \neq 1$, since if $\sigma(1) = 1$ and $\sigma(n) \neq n$, we can study the rotated permutation instead.

If $\sigma(1) \neq 1$, the only other option is $\sigma(1) = 3$, since otherwise $\sigma^{-1}(2) > \sigma^{-1}(3) > 1$ and $2 < 3 < \sigma(1)$ gives the forbidden pattern. Similarly $\sigma(2) > 5$ implies $\sigma^{-1}(4) > \sigma^{-1}(5) > 2$. So we get $\sigma(2) = 5$, and, by symmetry $\sigma(3) = 1$ and $\sigma(5) = 2$. Continuing the same argument and using induction we get $\sigma(2i) = 2i + 3$ and $\sigma(2i + 3) = 2i$, for $1 \leqslant i \leqslant \lfloor (n-3)/2 \rfloor$.

If $n$ is odd, only row and column $n-1$ are still empty, so $\sigma(n-1) = n-1$. This permutation, however, contains the pattern 321, since $n = \sigma(n-4) > \sigma(n-1) > \sigma(n) = n-4$.

For $n$ even, there are two possibilities to fill the remaining two rows and columns, namely by putting $\sigma(n-2) = n$ and $\sigma(n) = n-2$ or $\sigma(n-2) = n-2$ and $\sigma(n) = n$. Only the former, however, avoids 321.

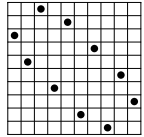

$\square$

## 4. **2413-avoiding doubly alternating permutations**

The doubly alternating 2413 permutations where conjectured by Guibert to be counted by the Catalan numbers. We prove this by showing them to possess a fairly simple block structure. First we need a technical lemma.

**Lemma 4.1.** *Let $\sigma \in \mathrm{DA}_n(2413)$, then*

(i) *$n$ odd $\implies$ $\sigma(1) = 1$*
(ii) *$n$ even $\implies$ $\sigma(1) = 1$ and $\sigma(n) = n$ or there is a $k$, $2 < k \leqslant n-1$, such that $\sigma(i) > \sigma(j)$ for all $i < k \leqslant j$.*

*Proof.* We can assume that $n \geqslant 3$, the smaller cases are trivial. Let $a = \sigma(1)$, and assume $\sigma(1) \neq 1$. By Lemma 2.2, $a$ must be odd. Now let $b = \sigma(\beta)$, where $\beta$ is the smallest number such that $\sigma(\beta) < a$. Note that $\beta \geqslant 3$, since $\sigma(2) > a$. Also, $\beta$ is odd since $\sigma(\beta - 1) > \sigma(\beta)$.

Let $c$ be the largest number such that $\gamma = \sigma^{-1}(c) < \beta$, see Figure 3. Thus $c > a$ and $c$ must be odd or $c = n$, since $\sigma^{-1}(c+1) > \beta > \sigma^{-1}(c) = \gamma$ if $c < n$.

If $\kappa > \beta$, then $\sigma(\kappa) < a$ or $\sigma(\kappa) > c$, otherwise we get the 2413 pattern. Therefore the rectangle, with NW corner $(2, a+1)$ and SE corner $(\beta - 1, c)$ contains exactly one dot in each row and column, so it is square, and hence $c = a + \beta - 2$ is even and thus $c = n$ is the only possibility. This proves the first assertion.

If $n$ is even, $\sigma(1) = 1$ implies $\sigma(n)$, by rotational symmetry, so the second assertion follows from $i < \beta \leqslant j \Rightarrow \sigma(j) < a \leqslant \sigma(i)$. □

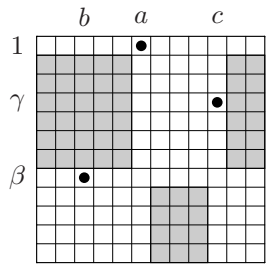

Figure 3. Illustration of Lemma 4.1, the shaded areas are empty.

As a direct consequence, we get the following corollary, which gives a very explicit description of what the DA 2413-avoiding permutations look like.

**Corollary 4.2.**

(i) $\sigma \in \mathrm{DA}_{2n+1}(2413)$ iff $\sigma = (1, \tilde{\sigma})$, where $(\tilde{\sigma}^r)^{-1} \in \mathrm{DA}_{2n}(2413)$.

(ii) $\sigma \in \mathrm{DA}_{2n}(2413)$ iff the permutation matrix of $\sigma$ is a block matrix, where all but the anti-diagonal blocks are empty. Any non-empty block $\nu$ has even size, $2k$, and can be written $\nu = (1, \tilde{\nu}, 2k)$, where $(\tilde{\nu}^r)^{-1} \in \mathrm{DA}_{2k-2}(2413)$.

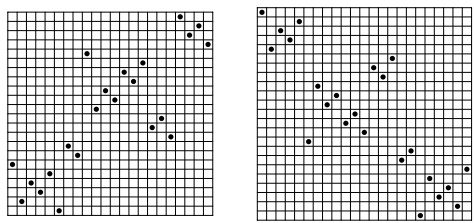

Figure 4. Example of the block structure of 2413-avoiding DA permutations.

The block structure condition in the corollary is in fact invariant under taking inverses, even though the pattern 2413 is not, therefore we get the following, slightly surprising result: $\mathrm{DA}_n(2413)$ and $\mathrm{DA}_n(3142)$ are not only the same size, but are actually the same sets. Therefore also $\mathrm{DA}_n(2413, 3142)$ is the same set. They are all counted by the Catalan numbers.

**Proposition 4.3.** $|\operatorname{DA}_n(2413)| = C_{\lfloor n/2 \rfloor}$.

*Proof.* First if $n$ is odd, we have as a direct consequence of Corollary 4.2(i) that $|\operatorname{DA}_n(2413)| = |\operatorname{DA}_{n-1}(2413)|$. If $n$ is even, then Corollary 4.2(ii) tells us that $\sigma \in \operatorname{DA}_n(2413)$ can be factored into blocks $\sigma_1, \sigma_2, \ldots, \sigma_m$, where $\sigma_i = (1, \tilde{\sigma}_i, |\sigma_i|)$ and $\tilde{\sigma}_i^r \in \operatorname{DA}(2413)$.

Let $D(x)$ be the generating function $D(x) = \sum_k |\operatorname{DA}_{2k}(2413)|x^k$. Then

$$D(x) = \sum_{i=0}^{\infty} (xD(x))^i = \frac{1}{1 - xD(x)}$$

which implies $xD(x)^2 - D(x) + 1 = 0$, i.e., the well know equation for the generating function of the Catalan numbers. Since $D(0) = 1$, we get $|\operatorname{DA}_{2k}(2413)| = C_k$. $\quad\square$

Another way to prove this is to construct a bijection with Dyck paths. We define $\Theta : \operatorname{DA}_{2n}(2413) \leftrightarrow \{$Dyck paths of length $2n\}$ recursively, by using Corollary 4.2.

(i) $\Theta(\emptyset) = \emptyset$

(ii) If $\sigma$ consists of a single block, so that $\sigma = (1, \tilde{\sigma}, 2n)$, then $\Theta(\sigma)$ is the Dyck path starting with a rise, ending with a descent and having the Dyck path $\Theta(\tilde{\sigma}^r)$ as the middle part.

(iii) If $\sigma$ can be factored into $k$ blocks, $\sigma_1, \ldots, \sigma_k$ (starting with the leftmost block), then $\Theta(\sigma)$ is the concatenation of the Dyck paths $\Theta(\sigma_1), \ldots, \Theta(\sigma_k)$.

The inverse is similarly defined, using recursion.

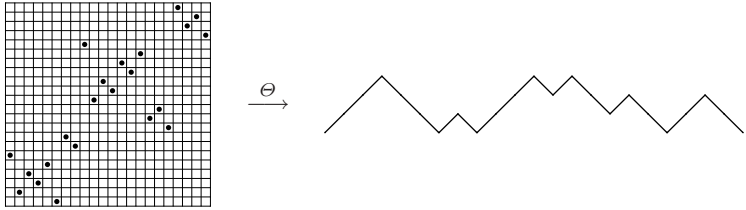

FIGURE 5. Example of the bijection between $\operatorname{DA}_{2n}(2413)$ and Dyck paths.

## 4.1. Doubly alternating Baxter permutations

A *Baxter permutation* is defined to be a permutation, $\sigma = (\sigma_i)_{i=1}^n$, such that for all $1 \leqslant i < j < k < l \leqslant n$,

$$\sigma_i + 1 = \sigma_l \text{ and } \sigma_j > \sigma_l \implies \sigma_k > \sigma_l \quad \text{and}$$
$$\sigma_l + 1 = \sigma_i \text{ and } \sigma_k > \sigma_i \implies \sigma_j > \sigma_i.$$

It is clear from this definition that if $\sigma$ avoids both 2413 and 3142 then it is a Baxter permutation, so we have

$$\operatorname{DA}_n(2413, 3142) \subset \{\sigma \in \operatorname{DA}_n : \sigma \text{ is Baxter}\}.$$

However, in [8], Guibert and Linusson showed that the doubly alternating Baxter permutations are counted by the Catalan numbers, so the sets must in fact be the same:

**Corollary 4.4.**

$$\{\sigma \in \mathrm{DA}_n : \sigma \text{ is Baxter}\} = \mathrm{DA}_n(2413, 3142) = \mathrm{DA}_n(2413) = \mathrm{DA}_n(3142).$$

It is also possible to prove this directly, without referring to the result by Guibert and Linusson.

**Lemma 4.5.** $\{\sigma \in \mathrm{DA}_n : \sigma \text{ is Baxter}\} \subset \mathrm{DA}_n(2413)$.

*Proof.* Assume $\sigma$ is Baxter, but not 2413-avoiding, and $d_1, d_2, d_3, d_4$, with $d_k = (i_k, j_k)$, constitute a 2413 pattern, such that $j_4 - j_1$ is as small as possible and given $j_1$ and $j_4$, $i_3 - i_2$ is as small as possible , as in Figure 6. The four areas shaded in the figure are empty, otherwise we would use one of those dots for the 2413-pattern. Now, let $\nu$ be the permutation having as permutation matrix the submatrix of $\sigma$, consisting of the rows $i_2 + 1, i_2 + 2, \ldots, i_3 - 1$ and columns $j_1 + 1, j_1 + 2, \ldots, j_4 - 1$. Since $\sigma$ is Baxter, $\nu$ cannot be empty. The DA condition implies that $\nu$ is up-down-alternating, $\nu^{-1}$ is down-up-alternating and $|\nu|$ is even. Also $\nu$ is 2413-avoiding, since otherwise this occurrence of 2413 would have been used instead of $d_1, \ldots, d_4$, and $\nu(1) \neq 1$, since $\nu^{-1}$ is down-up-alternating.

Rehashing the argument for Lemma 4.1, we can see that in fact no such permutation $\nu$ can exist. Defining $a$, $b$, $c$, $\beta$ and $\gamma$ in the same way as in the proof of Lemma 4.1 we get once again Figure 3. The difference is that now $a$ and $c$ must be even, whereas $\beta$ is still odd. But then $c = a + \beta - 2$ is odd, a contradiction. $\qquad\square$

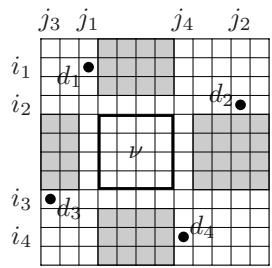

FIGURE 6. Illustration of Lemma 4.5. The shaded areas do not contain any dots.

## 5. 1234-avoiding doubly alternating permutations

In this section we construct a bijection between the doubly alternating 1234-avoiding permutations of size $2n$ and the ordinary 1234-avoiding permutations of size $n$ by using the Robinson-Schensted correspondence. Let $\lambda$ be a *Young diagram*, and denote by $\mathrm{SYT}(\lambda)$ the set of standard Young tableaux of shape $\lambda$. For a standard Young tableaux, $T$, let $\mathsf{row}_k(T)$ ($\mathsf{col}_k(T)$) denote the number of the row (column) for the entry $k$, counting from the top row (leftmost column), which is given the number one. The vector $\mathsf{row}(T)$ ($\mathsf{col}(T)$) is called the row (column) reading of $T$.

We define the set of *alternating standard tableaux* as

$$\mathrm{Alt}(\lambda) = \{T \in \mathrm{SYT}(\lambda) : \mathsf{col}(T) \text{ is up-down-alternating}\}$$
$$= \{T \in \mathrm{SYT}(\lambda) : \mathsf{row}(T) \text{ is down-up-alternating}\},$$

where the second equality is a consequence of the relative positions of two consecutive entries in a standard tableau, as shown in Figure 7.

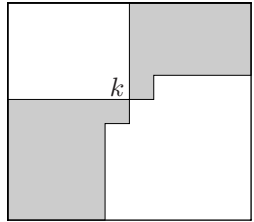

FIGURE 7. The shaded area denotes the possible positions for the entry $k + 1$, relative to the entry $k$, in a standard Young tableau.

The Robinson-Schensted correspondence, RSK, is a well known bijection between a permutation and a pair of standard Young tableaux of the same shape, see for example the book by Fulton [7]. An interesting fact is that doubly alternating permutations can be recognised by their RSK tableaux: They are both alternating if and only if the permutation is DA. In fact, Foulkes [5] proved a more general theorem, in which he counts the number of permutations with any given sequences of ups and downs for the permutation and its inverse. The following lemma is the key for proving Foulkes theorem.

**Lemma 5.1.** *Let* $\sigma \in \mathcal{S}_n$, $\mathrm{RSK}(\sigma) = (P, Q)$ *and* $1 \leqslant k < n$, *then*

$$k \text{ comes before } k + 1 \text{ in } \sigma \quad \Longleftrightarrow \quad \mathsf{row}_k(P) \geqslant \mathsf{row}_{k+1}(P).$$

*Proof.* First assume $k$ is inserted before $k + 1$ in the RSK bumping process. This means that $k + 1$ can never end up below $k$, since whenever they are in the same row only $k$ can be bumped down.

For the converse, assume $k + 1$ is inserted before $k$. Then $k$ will always be strictly above $k + 1$, since if $k$ is in the row exactly above $k + 1$ and is being bumped down (or $k + 1$ is in the first row and $k$ is about to be inserted), $k$ has to bump down $k + 1$ in the next step and thus stay above. $\qquad\square$

Let the *signature* of a word, $w = w_1 w_2 \ldots w_n$, be a sequence of $+$'s and $-$'s which has a $+$ in position $i$ iff $w_i < w_{i+1}$. For example, $\mathrm{signature}(4, 1, 5, 5, 6, 2, 2) = (-, +, -, +, -, -)$. We now get Foulkes theorem as a consequence of Lemma 5.1 and the fact that $\mathrm{RSK}(\sigma) = (P, Q)$ iff $\mathrm{RSK}(\sigma^{-1}) = (Q, P)$:

**Theorem 5.2** (Foulkes)**.** *Let $\sigma \in \mathcal{S}_n$ and* $\text{RSK}(\sigma) = (P, Q)$. *Then*

$$\text{signature}(\sigma^{-1}) = \text{signature}(\mathsf{col}(P)) = -\text{signature}(\mathsf{row}(P)),$$
$$\text{signature}(\sigma) = \text{signature}(\mathsf{col}(Q)) = -\text{signature}(\mathsf{row}(Q)).$$

The doubly alternating permutations are a special case:

**Corollary 5.3.** *Let $\sigma \in \mathcal{S}_n$ and* $\text{RSK}(\sigma) = (P, Q)$. *Then*

$$\sigma \in \text{DA}_n \quad \Longleftrightarrow \quad P, Q \in \text{Alt}(\lambda).$$

Let $T$ be an alternating standard tableau with $2n$ entries and at most three columns. We define the *pair column reading* $\mathsf{colpair}(T) = (w_i)_{i=1}^n$, where $w_i = \mathsf{col}_{2i-1}(T) + \mathsf{col}_{2i}(T) - 2$, i.e, $(1, 2) \mapsto 1, (1, 3) \mapsto 2$ and $(2, 3) \mapsto 3$, since Lemma 5.1 tells us that the only possibilities for the pairs are $(1, 2), (1, 3)$ and $(2, 3)$.

Let $w = (w_i)_{i=1}^l$ be a word, and $\mathsf{weight}(w) \stackrel{\text{def}}{=} (|\{i : w_i = k\}|)_{k \geqslant 1}$ be the *weight* vector of $w$. We call $w$ *Yamanouchi* (or a *ballot sequence*) if the weight of each prefix of $w$ is a partition, i.e., it is weakly decreasing.

**Lemma 5.4.** $\mathsf{colpair}$ *is a bijection between alternating standard tableaux with $2n$ elements and at most three columns and Yamanouchi words of length $n$ on three letters.*

*Proof.* Let $T$ be an alternating standard tableaux, with at most three columns. Then $\mathsf{colpair}(T) = (w_i)_{i=1}^n$ is, as noted above, a word on the letters 1,2 and 3, so we need to show that $\mathsf{colpair}(T)$ is Yamanouchi iff $\mathsf{col}(T)$ is an alternating Yamanouchi word.

First assume $\mathsf{colpair}(T)$ is Yamanouchi and let $v = (w_i)_{i=1}^k$ be an arbitrary prefix of $\mathsf{colpair}(T)$. Then $\mathsf{weight}(v) = (a, b, c)$, is a partition, i.e., $a \geqslant b \geqslant c$. Hence $\mathsf{weight}((\mathsf{col}_j(T))_{j=1}^{2k}) = (a + b, a + c, b + c)$ is also a partition. Since $\mathsf{col}(T)$ is alternating and $\mathsf{weight}((\mathsf{col}_i(T))_{i=1}^{2k+1})$ is a partition if $\mathsf{weight}((\mathsf{col}_i(T))_{i=1}^{2k+2})$ is, it follows that $\mathsf{col}(T)$ is an alternating Yamanouchi word.

For the converse, assume $\mathsf{col}(T)$ is an alternating Yamanouchi word and let $u = (\mathsf{col}_i(T))_{i=1}^{2k}$ be a prefix of $\mathsf{col}(T)$. Then $\mathsf{weight}(u) = (d, e, f)$ is a partition, so $\mathsf{weight}((\mathsf{colpair}_i(T))_{i=1}^k) = \frac{1}{2}(d + e - f, d + f - e, e + f - d)$ is a partition, which proves that $\mathsf{colpair}(T)$ is Yamanouchi. $\qquad\square$

Now we are ready to combine the bijections to get the bijection $\Phi : \mathcal{S}_n(1234) \to \text{DA}_{2n}(1234)$, defined by

$$\Phi(\sigma) = \text{RSK}^{-1}(\mathsf{colpair}^{-1}(\mathsf{col}(P)), \mathsf{colpair}^{-1}(\mathsf{col}(Q))),$$

where $\text{RSK}(\sigma) = (P, Q)$. See Figure 8 for an illustrative example.

**Theorem 5.5.** $\Phi$ *is a bijection, hence*

$$|\text{DA}_{2n}(1234)| = |\mathcal{S}_n(1234)|.$$

*Proof.* Let $\sigma \in \mathcal{S}_n(1234)$ and $\text{RSK}(\sigma) = (P, Q)$. RSK is a bijection between permutations and pairs of standard tableaux of the same shape such that if the permutation is 1234-avoiding iff the shape does not have more than three columns. From the definitions we know that $\mathsf{col}(P)$ and $\mathsf{col}(Q)$ are Yamanouchi words, so, by Lemma 5.4, $\mathsf{colpair}^{-1}(\mathsf{col}(P))$ and $\mathsf{colpair}^{-1}(\mathsf{col}(Q))$ are alternating standard tableaux with at most three columns. Their shapes are the same since the weights of $\mathsf{col}(P)$ and $\mathsf{col}(Q)$ are the same, which is a consequence of $P$ and $Q$ having the same shape. Applying the inverse of RSK and using Corollary 5.3, we get that $\Phi(\sigma) \in \mathrm{DA}_{2n}(1234)$.

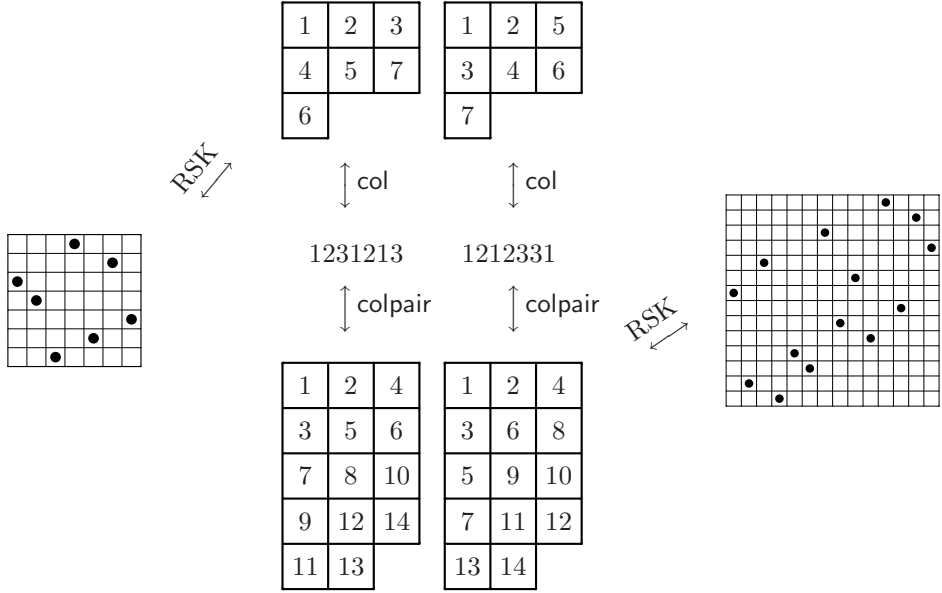

Figure 8. Example of the bijection $\Phi : \mathcal{S}_n(1234) \to \mathrm{DA}_{2n}(1234)$.

The converse is similar. $\qquad\square$

## 6. $21\tau$-avoiding doubly alternating permutations

The goal of this section is to find a bijection between $\mathrm{DA}_n(12\tau)$ and $\mathrm{DA}_n(21\tau)$, where $\tau$ is any permutation of $\{3, 4, \ldots, m\}$, $m \geqslant 3$, by using a well known bijection due to Babson and West [2]. The problem is that it is far from obvious that it will preserve the property of being doubly alternating. To show this we need a few definitions. During this section we assume $\tau$ to be fixed.

A dot, $d$, is called *active* if $d$ is the 1 or 2 in any $12\tau$ or $21\tau$ pattern in $\sigma$ and other dots are called *inactive.* Also the pair of dots, $(d_1, d_2)$, is called an *active pair* if $d_1 d_2$ is the 12 in a $12\tau$-pattern or the 21 in a $21\tau$-pattern.

**Lemma 6.1.** *Assume $\sigma \in \mathrm{DA}_n(12\tau) \cup \mathrm{DA}_n(21\tau)$ and $d = (i, j)$ is any active dot. Then $i$ and $j$ are odd.*

*Proof.* First assume $\sigma \in \mathrm{DA}_n(12\tau)$ and that $\sigma$ has a $21\tau$-pattern, otherwise there are no active dots. By inversion symmetry, we can assume that $d$ is the 1 in a $21\tau$ pattern. If $j = 1$, i.e. $\sigma^{-1}(1) = j$, then $j$ is odd by Lemma 2.2, and if $j > 1$, then $\sigma^{-1}(j-1) > \sigma^{-1}(j)$, since a dots to the north-west of $d$ would give a $12\tau$ pattern. Hence $j$ is odd. Also, to avoid the $12\tau$, $\sigma(i-1) > \sigma(i)$, so $i$ is odd as well.

Now assume instead $\sigma \in \mathrm{DA}_n(21\tau)$. Let $d_1, d_2, \ldots, d_m$, be the dots in a $12\tau$ pattern, with $d_k = (i_k, j_k)$. If $i_1$ is even then there is a descent from $i_1$ to $i_1 + 1$ and so the corresponding points along with tau will make the forbidden pattern. So $i_1$ is odd and the same argument applies to $i_2$, $j_1$, and $j_2$. $\qquad\square$

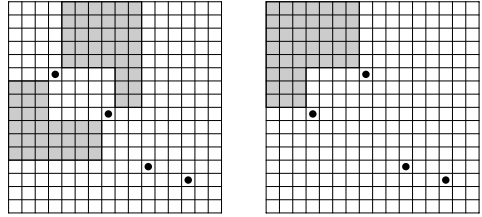

FIGURE 9. Illustration of the proof for Lemma 6.1, with $\tau = (3, 4)$. Shaded areas are forbidden.

We now define a Young diagram, $\lambda_\sigma$, consisting of the part of the board which contains the active dots. For a pair of dots, $d_1, d_2$, let $R_{d_1,d_2}$ to be the smallest rectangle with top left coordinates $(1, 1)$, such that $d_1, d_2 \in R_{d_1,d_2}$. Define

$$\lambda_\sigma \overset{\mathrm{def}}{=} \bigcup R_{d_1,d_2},$$

where the union is over all active pairs $(d_1, d_2)$. It is clear from the definition that $\lambda_\sigma$ is indeed a Young diagram (see Figure 10).

A *rook placement* (also known as *traversal* or *transversal*) of a Young diagram, $\lambda$, is a placement of dots, such that all rows and columns contain exactly one dot. If some of the rows or columns are empty we call it a *partial rook placement*. Furthermore, we say that a rook placement on $\lambda$ avoids the pattern $\tau$ if no rectangle, $R \subset \lambda$, contain $\tau$.

The definition of $\lambda_\sigma$ implies the following useful fact:

**Lemma 6.2.** *Let $\sigma \in \mathrm{DA}_n$ and $\mathrm{rp}(\lambda_\sigma)$ be the partial rook placement on $\lambda_\sigma$ induced by $\sigma$. Then*

$$\sigma \in \mathrm{DA}_n(12\tau) \quad \Longleftrightarrow \quad \mathrm{rp}(\lambda_\sigma) \text{ is } 12\text{-avoiding,}$$
$$\sigma \in \mathrm{DA}_n(21\tau) \quad \Longleftrightarrow \quad \mathrm{rp}(\lambda_\sigma) \text{ is } 21\text{-avoiding.}$$

The bijection we will use is due to Babson and West [2], which built on work by Simion and Schmidt [12] and West [14]. But we give here the more general result by Backelin, West, and Xin [3].

**Theorem 6.3** (Backelin, West, Xin). *Let $\tau$ be any permutation of $\{t+1, t+2, \ldots, m\}$. Then for every Young diagram $\lambda$, the number of $(t, t-1, \ldots, 1, \tau)$-avoiding rook placements on $\lambda$ equals the number of $(1, 2, \ldots, t, \tau)$-avoiding rook placements on $\lambda$.*

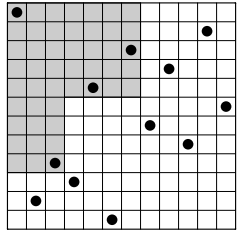

FIGURE 10. Example of $\lambda_\sigma$ for $\sigma = (1, 11, 7, 9, 5, 12, 8, 10, 3, 4, 2, 6)$
and $\tau = (3, 4)$.

We call two permutations of the same size *a-equivalent* if all the inactive dots are the same, and write $\sigma_1 \sim_a \sigma_2$. We shall see in Lemma 6.5 that this implies $\lambda_{\sigma_1} = \lambda_{\sigma_2}$.

**Lemma 6.4.** *If $\sigma \in \mathrm{DA}_n(12\tau) \cup \mathrm{DA}_n(21\tau)$ and $\nu \sim_a \sigma$, then $\nu$ is doubly alternating.*

*Proof.* Let $d = (i, j) \in \nu$ be a dot in an odd row, so that, by Lemma 6.1, both the dots in row $i - 1$ and row $i + 1$ are inactive (if they exists). If $d$ is inactive then all three dots also belong to $\sigma$ so $\nu(i - 1) > \nu(i) < \nu(i + 1)$. If $d$ is active there is a $\tau$-pattern to the SE of $d$, so if either of the dots in row $i - 1$ or row $i + 1$ are to the left of $d$, then this dot is active, since it creates either a $12\tau$ or $21\tau$ pattern together with $d$ and the $\tau$ pattern, giving a contradiction, so again $\nu(i-1) > \nu(i) < \nu(i+1)$. The same applies, by symmetry, to dots in odd columns. □

**Lemma 6.5.** *If $\sigma, \nu \in \mathrm{DA}_n$, then*

$$\sigma \sim_a \nu \quad \Longrightarrow \quad \lambda_\sigma = \lambda_\nu.$$

*Proof.* Let $\sigma \sim_a \nu$ be two arbitrary $a$-equivalent DA permutations and assume $s = (i, j)$ is a SE corner of $\lambda_\sigma$. We need to show that $s \in \lambda_\nu$, so that $\lambda_\sigma \subseteq \lambda_\nu$ and thus, since $\sim_a$ is reflexive, $\lambda_\nu = \lambda_\sigma$.

Let $R_{d_1, d_2} \subset \lambda_\sigma$ be a rectangle, such that $s \in R_{d_1, d_2}$. Such a rectangle must exist, otherwise $s$ could not belong to $\lambda_\sigma$. Hence there is a $\tau$-pattern to the SE of $s$ and one of the $d_k$ is in row $i$ and one (possibly the same one) is in column $j$. But, as $\nu \sim_a \sigma$, they have the same inactive dots, so there must also exist a dot $d'_1 \in \lambda_\nu$ in row $i$ and a dot $d'_2 \in \lambda_\nu$ in column $j$. If $d'_1$ is east of $s$ or if $d'_2$ is south of $s$ then $s \in \lambda_\nu$. Hence we can assume $d'_1$ and $d'_2$ to be weakly NW of s. If $d'_1 \neq d'_2$, then $s \in R_{d'_1, d'_2} \subset \lambda_\nu$, since the $\tau$-pattern is still SE of $s$, and if $d'_1 = d'_2 = s$ then clearly $s \in \lambda_\nu$. □

Now we are ready to construct a bijection $\Psi : \mathrm{DA}_n(12\tau) \to \mathrm{DA}_n(21\tau)$. Let $\sigma \in \mathrm{DA}_n(12\tau)$, so that the restriction of $\sigma$ to $\lambda_\sigma$ is a partial 12-avoiding rook placement. By Theorem 6.3 (ignoring the empty rows and columns) and Lemma 6.5, there exists a unique 21-avoiding (partial) rook placement on $\lambda_\sigma$, with the same rows and columns empty, which we combine with the inactive dots of $\sigma$ to get $\Psi(\sigma)$. By Lemma 6.4, $\Psi(\sigma)$ is DA, and Lemma 6.2 says that it avoids $21\tau$. It is also clear from Theorem 6.3 that it is indeed a bijection. We have thus bijectively shown:

**Theorem 6.6.** *Let $\tau$ be any permutation of $\{3, 4, \ldots, m\}$, $m \geqslant 3$. Then*

$$|\operatorname{DA}_n(21\tau)| = |\operatorname{DA}_n(12\tau)|.$$

As a special case we have

**Corollary 6.7.** $\quad |\operatorname{DA}_n(2134)| = |\operatorname{DA}_n(1234)|.$

## 7. Other patterns with the same number sequence as $\mathcal{S}_n(1234)$

By examining all the patterns of length four with computer, Guibert found 15 different cases that all seem to give rise to the same sequence, $|\mathcal{S}_n(1234)|$. Using Theorems 5.5 and 6.6, inversion, rotation and Proposition 7.2 below, we get altogether ten bijections, see Figure 11. However, to prove that all of them are indeed the same we would need five more bijections. In fact, we conjecture that the number of permutations are the same in all the cases given below.

**Conjecture 7.1** (Guibert)**.**

$$|\operatorname{DA}_{2n}(1234)| = |\operatorname{DA}_{2n+1}(1243)|$$
$$= |\operatorname{DA}_{2n}(1432)|$$
$$= |\operatorname{DA}_{2n+1}(1432)|$$
$$= |\operatorname{DA}_{2n}(2341)|$$
$$= |\operatorname{DA}_{2n}(3421)|$$

One can note that many of the patterns in the conjecture are of the same type as treated in Theorem 6.3, but the proof does not work here, except for 2134, since the bijections destroy the DA property.

**Proposition 7.2.** $\quad |\operatorname{DA}_{2n}(2143)| = |\operatorname{DA}_{2n+1}(3412)| = |\operatorname{DA}_{2n+2}(3412)|.$

*Proof.* Let $\sigma \in \operatorname{DA}_n(3412)$, with $n \geqslant 4$. If $\sigma(1) > 1$, we get the forbidden pattern on the rows 1, 2, $\sigma^{-1}(1)$, $\sigma^{-1}(2)$, so $\sigma(1) = 1$. Let $\tilde{\sigma}$ be the permutation with the first row and column of $\sigma$ removed. It is clear that if $n$ is odd then $\tilde{\sigma}^c \in \operatorname{DA}_{n-1}(2143)$ iff $\sigma \in \operatorname{DA}_n(3412)$, and if $n$ is even then $\tilde{\sigma}^\# \in \operatorname{DA}_{n-1}(3412)$ iff $\sigma \in \operatorname{DA}_n(3412)$. $\quad\square$

## 8. Avoiding pairs of patterns of length four

When we have two patterns of length four, there are a huge number of cases. We have not yet studied many of these, but would like to give a flavour of what can happen by presenting one result and two conjectures. Combining the results in Sections 4 and 5 we get

**Proposition 8.1.**

$$|\operatorname{DA}_n(1234, 2413)| = \begin{cases} F_{n/2}, & \textit{if } n \textit{ is even}, \\ 2, & \textit{if } n = 5, \\ 1, & \textit{otherwise}, \end{cases}$$

| $\mathcal{S}_n(1234)$ | $\xleftrightarrow{\text{Th. 5.5}}$ | $\mathrm{DA}_{2n}(1234)$ $\xleftrightarrow{\text{Th. 6.6}}$ $\mathrm{DA}_{2n}(2134)$ $\xleftarrow{\#}$ $\mathrm{DA}_{2n}(1243)$ |
| | $\xleftrightarrow{\text{Th. 6.6}}$ | $\mathrm{DA}_{2n}(2143)$ $\xleftrightarrow{\text{Pr. 7.2}}$ $\mathrm{DA}_{2n+1}(3412)$ $\xleftrightarrow{\text{Pr. 7.2}}$ $\mathrm{DA}_{2n+2}(3412)$ |
| $\mathrm{DA}_{2n+1}(1243)$ | $\xleftrightarrow{\text{Th. 6.6}}$ | $\mathrm{DA}_{2n+1}(2143)$ |
| $\mathrm{DA}_{2n}(1432)$ | $\xleftarrow{\#}$ | $\mathrm{DA}_{2n}(3214)$ |
| $\mathrm{DA}_{2n}(2341)$ | $\xleftrightarrow{-1}$ | $\mathrm{DA}_{2n}(4123)$ |
| $\mathrm{DA}_{2n}(3421)$ | $\xleftarrow{\#}$ | $\mathrm{DA}_{2n}(4312)$ |
| $\mathrm{DA}_{2n+1}(1432)$ | | |

FIGURE 11. Known bijections between the sequences conjectured to be $|\mathcal{S}_n(1234)|$.

*where the $F_n$ are the Fibonacci numbers.*

*Proof.* Let $\sigma \in \mathrm{DA}_{2n}(1234, 2413)$. By Corollary 4.2, $\sigma$ can be factored into blocks, $\sigma_1, \sigma_2, \ldots \sigma_k$. As $\sigma$ avoids 1234 must each block be either 12 or 1324, since each of them have a dot in the NW corner and the SE corner. Hence

$$|\mathrm{DA}_{2n}(1234, 2413)| = |\mathrm{DA}_{2n-2}(1234, 2413)| + |\mathrm{DA}_{2n-4}(1234, 2413)|,$$

and since $|\mathrm{DA}_0(1234, 2413)| = |\mathrm{DA}_2(1234, 2413)| = 1$, we get the Fibonacci numbers.

If $\sigma \in \mathrm{DA}_{2n+1}(1234, 2413)$, then $\sigma(1) = 1$. Let $\tilde{\sigma} = (2n+1-\sigma(i))_{i=2}^{2n+1}$ be the permutation constructed from $\sigma$ by removing the first row and column and then flipping horizontally. Then $\tilde{\sigma} \in \mathrm{DA}_{2n}(321, 2413)$, which by Proposition 3.1(iii) gives two possibilities if $2n \geqslant 4$, namely $(1, 3, 2, 5, 4, \ldots, 2n-1, 2n-2, 2n)$ and $(3, 5, 1, 7, 2, 9, 4, \ldots, n, n-4, n-2)$. However, only the former avoids 2413 if $n \geqslant 6$, so we get the desired result. $\square$

The following two conjectures have been verified by computer calculations up to $n = 23$.

**Conjecture 8.2.**

$$|\mathrm{DA}_n(1234, 3214)| = \begin{cases} F_{n-1}, & \text{if } n \text{ is even,} \\ 1, & \text{if } n = 1 \text{ or } n = 3, \\ F_{n-1} - F_{n-7}, & \text{otherwise.} \end{cases}$$

**Conjecture 8.3.**

$$|\mathrm{DA}_n(1234, 2134)| = \begin{cases} C_{n/2}, & \text{if } n \text{ is even,} \\ 1, & \text{if } n = 1 \text{ or } n = 3, \\ C_{(n-5)/2}^{(4)}, & \text{otherwise.} \end{cases}$$

Here $C_n^{(4)}$ is the fourth difference of the Catalan numbers, defined recursively by $C_n^{(0)} = C_n$ and $C_n^{(i+1)} = C_{n+1}^{(i)} - C_n^{(i)}$. By collecting the terms and simplifying

we get

$$C_n^{(4)} = C_{n+4} - 4C_{n+3} + 6C_{n+2} - 4C_{n+1} + C_n$$
$$= 9C_n \frac{9n^4 + 54n^3 + 135n^2 + 122n + 40}{(n+2)(n+3)(n+4)(n+5)}.$$

Linköpings Universitet, 581 83 Linköping, Sweden
*E-mail address*: `erouc@mai.liu.se`

**34**

# Paper 2

**36**

# EXTENDED PATTERN AVOIDANCE AND GENERATING GRAPHS

ERIK OUCHTERLONY

ABSTRACT. The matrix representation of a permutation is a square 0-1-matrix such that every row and column has exactly one 1. A partial permutation is a rectangular matrix which is a sub-matrix of a permutation matrix. In this paper we are studying partial permutations which can be extended to pattern avoiding permutations, i.e., partial permutations which are the upper left corner of a pattern avoiding permutation. Many cases are solved, including all combinations of patterns of length three and some patterns pairs where the lengths are three and four respectively. The method used is generating graphs, a refined version of the generating tree method, which has help us to find all the solutions very efficiently.

## 1. Preliminaries

### 1.1. Introduction

The subject of pattern avoiding permutations has been gathering a lot of interest lately and generalisations in many different directions have been suggested. Much of the previous interest have been emphasised on generalising the patterns, as was done by Babson and Steingrímsson [2]. In the current paper we instead generalise the permutations.

The idea of extended pattern avoidance originated in the study of Fulton's essential set for 321-avoiding permutations in a paper by Eriksson and Linusson [4]. This was later developed by Linusson [5] to cover all patterns of size three. In this paper we build upon these results and have developed a method which allow us to find generating functions for many different pattern sets in a very efficient manner. The method, called the generating graph method, is based on generating trees and combines an algorithm for finding the generating tree with an automaton representation which gives an important visual aid in more complicated cases. We believe the method could be applied to many other problems as well, both in pattern avoidance and other fields.

In Section 2 we introduce the generating graphs and present algorithms to find them with computer assistance and in Section 3 this is generalised to extended pattern avoidance. The final three sections are concerned with special cases of extended pattern avoidance, in Section 4 we study singleton patterns of length three, in Section 5 (and in Appendix A) all combinations of patterns of length three are treated and in Section 6 pairs of patterns of lengths three and four, respectively, are examined.

I thank Svante Linusson for introducing me to the problem and for numerous fruitful discussions. Thanks also to Bruce Sagan for providing many corrections, suggestions and comments.

## 1.2. **Notation**

Let $\sigma \in \mathcal{S}_n$ be a permutation. Then the *permutation matrix* of $\sigma$ is a 0-1-matrix, $(a_{ij})_{i,j=1}^n$, such that $a_{ij} = 1$ iff $\sigma(i) = j$. We use the word *dots* for the ones in the permutation matrix in accordance with the usual way of illustrating permutation matrices, where the zeroes are left empty and the ones are changed to dots, see Figure 1.

$$\sigma = 41523 \qquad \qquad$$

FIGURE 1. A permutation in word and matrix form.

As a slight abuse of notation, we often use the same notation for both the permutation and its matrix and, for example, speak about the NE corner of a permutation. A position $k < n$ of the permutation $\sigma$ is an *ascent* (or *rise*) if $\sigma(k) < \sigma(k+1)$ and a *descent* otherwise. A permutation without ascents is *decreasing* and if it has no descent it is *increasing*.

The permutation $\sigma$ *contains* the pattern $\tau \in \mathcal{S}_m$ if the permutation matrix of $\tau$ is a submatrix of $\sigma$ and it *avoids* $\tau$ otherwise. This also generalises to multiple patterns and if $\tau = \{\tau_1, \tau_2, \ldots, \tau_t\}$ is a *pattern set* we write

$$\mathcal{S}_n(\tau) = \mathcal{S}_n(\tau_1, \tau_2, \ldots, \tau_t) = \{\sigma \in \mathcal{S}_n : \sigma \text{ avoids } \tau_i \text{ for all } i = 1, 2, \ldots, t\}$$

for the set of permutations avoiding all of the patterns $\tau_1, \tau_2, \ldots, \tau_t$. The ordinary generating function for the $\tau$-avoiding permutations is

$$\mathcal{S}_\tau(x) = \sum_n |\mathcal{S}_n(\tau)| x^n.$$

Two sets of patterns, $\tau$ and $\tau'$, are called *Wilf-equivalent* if $|\mathcal{S}_n(\tau)| = |\mathcal{S}_n(\tau')|$ for all $n$.

The *reverse* and the *complement* of a permutation $\sigma$ is defined by

$$\sigma^r = (\sigma(n+1-i))_{i=1}^n$$
$$\sigma^c = (n+1-\sigma(i))_{i=1}^n,$$

which corresponds to reflecting the matrix vertically and horizontally, respectively. We use these operations, as well as the inverse, on sets of permutations by applying them component wise. It is clear that from the symmetry of the permutation matrices that applying these operations does not affect the number of permutations avoiding these patterns, hence $\tau$, $\tau^{-1}$, $\tau^r$ and $\tau^c$ are all Wilf-equivalent for any set of patterns $\tau$.

## 2. Generating trees and generating graphs

### 2.1. Definitions

A common method for counting pattern avoiding permutations is to use *generating trees*, first introduced by Chung et al [3] for studying Baxter permutations. The idea is to use the fact that if a permutation, $\sigma$, of size $n > 0$ avoids the set of patterns, $\tau$, and if $\sigma'$ is the permutation created from $\sigma$ by removing (e.g.) the rightmost dot and its row and column, then $\sigma'$ is also a $\tau$-avoiding permutation. Thus we can create a tree, with the $\tau$-avoiding permutations as the nodes, by putting the zero permutation as the root and for each $\tau$-avoiding permutation, associate a parent node, by removing the row and column with the rightmost dot.

To define this formally, we need a *parent rule*, which is a function, $\mathsf{p} : \mathcal{S}_{n+1}(\tau) \to \mathcal{S}_n(\tau)$, where $n \geqslant 0$, which determines how to build the tree. The standard parent rule is the one described above, but it is also possible to create other parent rules, for example by removing dots in any of the other three directions. This choice can be very important to make the generating tree a simple as possible, as we will see in Example 2.4.

**Definition 2.1.** Let $\tau$ be a set of patterns and $\mathsf{p} : \mathcal{S}_{n+1}(\tau) \to \mathcal{S}_n(\tau)$ a parent rule. Then the *unlabelled generating tree*, $T_{\mathsf{p}}(\tau)$, is the directed tree with $\mathcal{S}(\tau)$ as the node set and where $(\sigma, \sigma')$ is an edge iff $\mathsf{p}(\sigma') = \sigma$.

With this definition, we have

$$\mathcal{S}_\tau(x) = \sum_w x^{|w|},$$

where the sum is over all walks in $T_{\mathsf{p}}(\tau)$ starting at the zero permutation. Hence we could use the generating tree to calculate the generating function. To do this, however, we need to examine the structure of the unlabelled generating tree and add labels to the nodes so that each label uniquely determines the labels of the children. This method has proved to be very successful, and has helped solving many problems, see for example [7]. However, for more complicated problems it could be very difficult to guess a good way to label the nodes, so we will demonstrate algorithms to do this by computer in Section 2.2.

To define what a proper labelling is, we need an equivalence relation. For any node, $\sigma \in V(T_{\mathsf{p}}(\tau))$, let $T_{\mathsf{p}}^\sigma(\tau)$ be the subgraph of $T_{\mathsf{p}}(\tau)$ containing all nodes that can be reached from $\sigma$.

**Definition 2.2.** Let $\sigma_1, \sigma_2 \in V(T_{\mathsf{p}}(\tau))$, then $\sigma_1$ and $\sigma_2$ are $\mathsf{p}$-*equivalent* iff $T_{\mathsf{p}}^{\sigma_1}(\tau)$ and $T_{\mathsf{p}}^{\sigma_2}(\tau)$ are isomorphic as graphs.

Using this equivalence relation, we can now get the generating tree by labelling each node with its equivalence class. In previous work, the generating tree has been represented using *rewrite rules* or *succession rules*, where the nodes are listed together with their outgoing edges. We will instead use a graph, since this gives a more visual picture of the generating tree, which is especially important when we study extended pattern avoidance, where the tree is weighted.

**Definition 2.3** (Generating graphs)**.** Let $T_{\mathsf{p}}(\tau)$ be the unlabelled generating tree for the pattern set $\tau$. We define the *generating graph*, $G_{\mathsf{p}}(\tau)$, to be the weighted directed graph having the set of equivalence classes under $\mathsf{p}$-equivalence as vertex set and where $e = (v, w)$ is an edge with weight $k$ iff for every $\sigma \in v$ there are exactly $k$ different edges $(\sigma, \mu_i) \in E(T_{\mathsf{p}}(\tau))$, such that $\mu_1, \ldots, \mu_k \in w$. The edges with weight $k > 1$ can also be represented as a $k$-tuple edge.

Now we can calculate the generating graph by counting the walks in the generating graph instead,

$$\mathcal{S}_\tau(x) = \sum_w \mathsf{weight}(w)\, x^{|w|},$$

where the sum is over all walks in $G_{\mathsf{p}}(\tau)$ starting at the class containing the zero permutation and the weight is the product of the weights of the all the edges in the walk.

We illustrate the method by looking at some well known results. Here we use the term *site* for a position to the right of the permutation where insertion of a dot could be made, i.e., between the rows. We refer to site $k$ as the site between row $k-1$ and row $k$, with site 1 and site $n+1$ being the top and bottom site, respectively. A site is called *active* if inserting a dot there gives a permutation which still avoids the pattern set in question.

**Example 2.4** (Even Fibonacci)**.** The permutations avoiding $\{213, 4123\}$ were enumerated in [7] using generating trees, and we will see that the generating graph in this case is

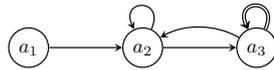

where the class $a_k$ contains the permutations with $k$ active sites. Let $i$ be the row with the rightmost dot. The sites $i$ and $i + 1$ are always active, since both the rightmost dot and a newly inserted dot above or below it could not both be part of either of the patterns. The sites below $i + 1$ are not active due to the pattern 213 and of the sites above site $i$, only $i - 1$ could be active, since the sites higher up are prohibited by 4123, since the two rows above the rightmost dot must be increasing. Therefore, a permutation has three active sites iff the site $i - 1$ is active, which occurs when the dot a row above the rightmost dot is not the 1 in a 123 pattern. By analysing the permutations with three active sites, one can see that all of them have one child with two active sites and two children with three active sites. Similarly, all the permutations with two active sites have one child with two active sites and one with three active sites. Hence the class $a_k$ contains the permutations with $k$ active sites. See Figure 2 for an example.

When calculating the generating function we use the name of the node as the generating function for walks starting at $a_1$ and ending at the given node. The graph immediately gives us a linear system of equations in these variables,

$$a_1 = 1$$
$$a_2 = x(a_1 + a_2 + a_3)$$
$$a_3 = x(a_2 + 2a_3),$$

which we can solve to get

$$F(x) = a_1 + a_2 + a_3 = 1 + \frac{x(1-2x)}{1-3x+x^2} + \frac{x^2}{1-3x+x^2} = \frac{1-2x}{1-3x+x^2},$$

which in turn implies that $|\mathcal{S}_n(213, 4123)| = F_{2n}$, where $F_n$ is the $n$:th Fibonacci number.

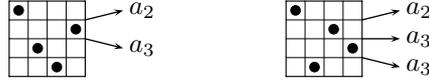

Figure 2. Illustration of Example 2.4. Two permutations in the classes $a_2$ and $a_3$, respectively, where an arrow indicates an active sites and to which class the permutation would belong if a dot is inserted at the site.

Another way to find the generating function when the graph is finite is to regard the graph as a *finite automata* as defined in [1] and use the *transfer matrix* method (see [6]) for calculating the generating function.

We label the nodes with $1, 2, \ldots, N$ and let $A$ be the $N \times N$-matrix, such that the entry $a_{ij}$ is the sum of the weights of the edges between node $a_j$ and node $a_i$. For the current example, we get

$$A = \left( \begin{array}{ccc} 0 & 0 & 0 \\ 1 & 1 & 1 \\ 0 & 1 & 2 \end{array} \right)$$

as the transfer matrix, which tells us how many walks of lengths one we have in the graph. To get the generating function, walks of arbitrary length must be calculated, and hence

$$T = \sum_{n \geqslant 0} x^n A^n = (I - xA)^{-1} = \frac{1}{1-3x+x^2} \left( \begin{array}{ccc} 1 & 0 & 0 \\ x(1-2x) & 1-2x & x \\ x^2 & x & 1-x \end{array} \right).$$

Since we begin on node $a_1$ and may end at any node, we again get the generating function

$$F(x) = (\, 1 \ 1 \ 1 \,) \, T \, (\, 1 \ 0 \ 0 \,)^t = \frac{1-2x}{1-3x+x^2}.$$

The transfer matrix method can be applied on any finite graph, which makes them in a sense trivial to solve. However, not all generating graphs are finite. If, for example, we had changed the parent rule, so that the permutation was growing downwards instead of towards the right, which corresponds to avoiding the inverse pattern set, $\{213, 2341\}$, with the normal parent rule, the generating graph becomes more complicated, namely

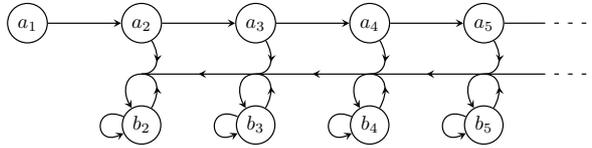

Note that we allow edges which follow the same route to be joined in order to make the figure clearer. For example, the node $a_4$ has edges to $a_5$, $b_2$, $b_3$ and $b_4$.

**Example 2.5** (The Catalan graph). It is well known that all patterns of length three have the same succession rule,

$$(k) \to (2)(3) \cdots (k+1),$$

where each class is labelled by the number of children in the tree. This is known as the Catalan tree, and its representation as a generating graph, which we call the *Catalan graph* is:

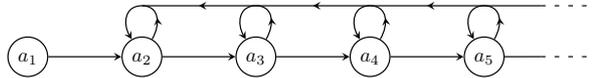

The generating function for the Catalan numbers is $c(x) = \frac{1-\sqrt{1-4x}}{2x}$, which can be seen from the graph by observing that any nonzero walk on the graph can be split into three parts: the part of the walk until just before the last time node $a_2$ is reached, the step to $a_2$ and the walk from $a_2$, never reaching $a_2$ again. This gives the equation $c(x) = 1 + c(x)xc(x)$, which, under the condition $c(0) = 1$ has the only solution given above.

One important feature of the generating graphs is that if we happen to get the same (or isomorphic) graphs for two different pattern sets, we automatically get a bijection between the permutations avoiding these sets, since both are in bijection with the walks on the graph. If the graph has multiple edges (i.e. weight greater than one), the bijection will not be unique, however, since we need for each of these edges to associate a rule which defines which active site will be used in this case.

2.2. **The generating graph algorithm**

The reason the generating graphs are useful is the fact that it is straight forward to implement an algorithm to calculate candidates for them with the help of computers. Even though it is possible to find the generating graph by hand in many cases, one quickly realises the need for computers when the graph is more complicated.

In order to calculate the generating graph, we first use computer program to find a good candidate for the generating graph by considering a finite part of the unlabelled generating tree. This can be done by a very simple algorithm: Assume we want to find candidates for the classes for all nodes at depth $N$ in the tree and make sure all nodes in each potential class have the same subtrees at least to depth $T$.

**Algorithm 2.6.**

    (1) Calculate the unlabelled generating tree to a fixed depth, $N + T$.

(2) Label each node with its subtree down to depth $T$.

(3) Classify the nodes by the labels.

To find which value of $T$ is large enough, we can increase it until no new potential classes are created. Also the variable $N$ should be selected appropriately, so that all classes will be included or, if the number of classes are infinite, a large enough number of classes are found, so that the structure of the graph can be guessed. Unfortunately, this crude algorithm is rather slow, since the width of the tree typically increases exponentially, which makes it unsuitable for more complicated generating graphs. On the other hand, many of the interesting cases do give rise to fairly simple generating graphs for which a candidate can be found without a problem even with this algorithm.

A slightly more sophisticated version of the algorithm takes into account that whenever we have guessed what class a given node belongs to, this will determine the classes for all of its children, so there is no need to continue calculating the unlabelled generating tree below this node.

**Algorithm 2.7.** The algorithm is recursive, starting with the root as $\nu$. We also have an, initially empty, associative array, $A$, which is a set of pairs, where each pair consists of a key and a value. Here the key is the subtree and the value is the equivalence class.

(1) Calculate the subtree of $\nu$ down to depth $T$.

(2) If the subtree is isomorphic to a subtree in $A$, put $\nu$ into that class. Otherwise, create a new class, store it in $A$ and apply the algorithm recursively for each child of $\nu$, unless we have reached the depth $N$.

In these two algorithms it is clear that two nodes which are in the same class of the generating graph must end up in the same class in the guessed graph. The opposite is, however, not true, so this needs to be check by hand. This is, however, much easier to do once we have a good candidate to work with, so that the only problem remaining is to characterise the different classes and confirm that they behave as the graph stipulates.

2.3. **Three lemmas on graphs**

The following three lemmas are used in many of the propositions in the following sections, since structures of these types are very common in the generating graphs of pattern avoiding permutations. In the graphs below, the $f_k$ nodes are thought of as part of a larger graph, and $f_k$ denotes, as before, the generating function for walks beginning at node $a_1$ and ending in $f_k$.

**Lemma 2.8.** *The generating function for the walks ending at* $b_k$ *in the graph*

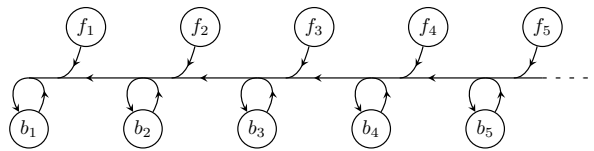

*is* $\frac{x}{1-x} \sum_{j=0}^{\infty} \frac{f_{j+k}}{(1-x)^j}$.

*Proof.* Let $b_k^{(j)}$ be the generating function for the walks going through $f_j$ and ending at $b_k$. By the word "through", we understand that $f_j$ is the last of the $f$-nodes that the walk visits. Then $b_k^{(j)} = x/(1-x)^{j-k+1} f_j$ if $j \geqslant k$, since the walk could go around any of $j - k + 1$ loops on its way from $f_j$ to $b_k$. Hence,

$$b_k = \sum_{j=k}^{\infty} b_k^{(j)} = \frac{x}{1-x} \sum_{j=0}^{\infty} \frac{f_{j+k}}{(1-x)^j}. \qquad \square$$

**Lemma 2.9.** *The generating function for the walks ending at $b_k$ in the graph*

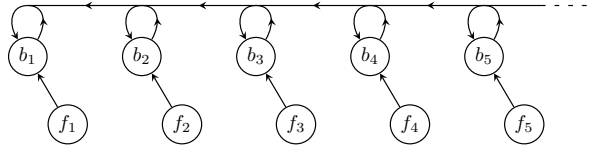

*is* $\frac{x}{1-x} f_k + (\frac{x}{1-x})^2 \sum_{i=0}^{\infty} \frac{f_{i+k+1}}{(1-x)^i}$.

*Proof.* Let $b_k^{(j)}$ be the generating function for walks through $f_j$ ending at $b_k$. Then, if $k < j$,

$$b_k^{(j)} = x \sum_{i=k}^{j} b_i^{(j)} f_j = \frac{x}{1-x} \sum_{i=k+1}^{j} b_i^{(j)} f_j = \frac{1}{1-x} b_{k+1}^{(j)}$$

$$= \frac{1}{(1-x)^{j-k-1}} b_{j-1}^{(j)} = \frac{x^2}{(1-x)^{j-k+1}} f_j$$

By summing the terms we get the given formula. $\qquad \square$

**Lemma 2.10.** *The generating function for the walks ending at $b_k$ in the graph*

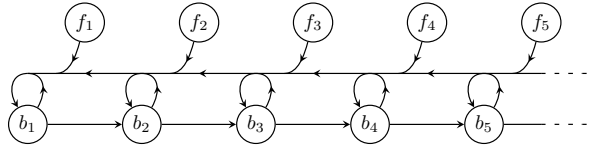

*is* $\sum_{m=0}^{k-1} (xc(x))^{m+1} \sum_{i=0}^{\infty} f_{k+i-m} c(x)^i$, *where $c(x)$ is the generating function for the Catalan numbers.*

*Proof.* Let $b_k^{(j,m)}, m \leqslant j, m \leqslant k$, be the generating function for walks through $f_j$ ending at $b_k$ which passes $b_m$ but not $b_{m-1}$. In order to split these walk in manageable pieces, we introduced virtual nodes, $v_i$, which are located right above $b_i$ in the figure above, i.e., at the switching points. The walk is split in the following way:

- Walk to $v_j$, ending at the moment when the walk for the first time is about go to the left from $v_j$: $c(x)$
- Walk to $v_{j-1}$, ending at the moment when the walk for the first time is about go to the left from $v_{j-1}$: $c(x)$
- Etc, until:
- Walk to $v_m$: 1

- Walk to $b_m$: $x$
- Walk until the last time it visits $b_m$: $c(x)$
- Walk to $b_{m+1}$: $x$
- Walk until the last time it visits $b_{m+1}$: $c(x)$
- Etc, until:
- Walk to $b_k$: $x$
- Walk until the last time it visits $b_k$: $c(x)$

so we get

$$b_k^{(j,m)} = (xc(x))^{k+1-m} c(x)^{j-m} f_j.$$

Now, let $b_k^{(m)}, m \leqslant k$, be the generating function for walks ending at $b_k$, which passes $b_m$ but not $b_{m-1}$. Then

$$b_k^{(m)} = \sum_{j=m}^{\infty} b_k^{(j,m)} = (xc(x))^{k-m+1} \sum_{j=0}^{\infty} f_{j+m} c(x)^j$$

and hence

$$b_k = \sum_{m=1}^{k} b_k^{(m)} = \sum_{m=0}^{k-1} (xc(x))^{m+1} \sum_{j=0}^{\infty} f_{k+j-m} c(x)^j. \qquad \square$$

## 3. Extended pattern avoidance

If we want to a generalise permutations, it is helpful to regard permutations as perfect matchings of complete bipartite graphs. It then comes very natural to simply abandon the "perfection" and allow any matching of the graph. In the matrix representation this means that we no longer require exactly one entry in every row and column since empty rows and columns are allowed.

**Definition 3.1.** A *partial permutation* is a rectangular 0-1-matrix, such that each row and column has at most one dot. The set of all partial permutations with $d$ dots, $r$ empty (i.e. dot-free) rows and $c$ empty columns, where $d,r,c \geqslant 0$, is denoted $\mathcal{S}_{d,c,r}$.

Note that we allow matrices with zero row and/or columns and, e.g., the partial permutations in $\mathcal{S}_{0,1,0}$ and $\mathcal{S}_{0,0,1}$ are not regarded as being equal even though they both are matrices without elements. As with permutations, partial permutations can be expressed in word form, but here empty rows are marked with an underscore, see Figure 3. This allow us to talk about ascents and descents as well as increasing and decreasing partial permutations in the same way as for permutations, where the only change is that if either $k$ or $k+1$ is an empty row, then $k$ is neither an ascent or a descent. Note that the word form does not uniquely determine a partial permutation, since empty columns on the right does not affect it, so we need to indicate the number of columns to make it unique.

We say that a permutation $\sigma \in \mathcal{S}_{d+c+r}$ is an *extension* of the partial permutation $\rho \in \mathcal{S}_{d,c,r}$ if $\rho$ is the NW corner of $\sigma$. This can be understood as extending $\rho$ by filling the empty rows and columns to the right and below.

When studying pattern avoidance on partial permutations, it would not be a true generalisation if we just used the same definition as for permutations since in that case the empty rows and columns do not come into play at all. Instead we utilise extensions in the following way:

**Definition 3.2** (Extended pattern avoidance). Let $\tau = \{\tau_1, \tau_2, \ldots, \tau_t\}$ be a set of patterns. A partial permutation $\rho \in \mathcal{S}_{d,c,r}$ is said to *extendably avoid* the pattern set $\tau$, if $\rho$ has an extension in $\mathcal{S}_{d+c+r}(\tau)$. The set of extendably $\tau$-avoiding partial permutations is $\mathcal{S}_{d,c,r}(\tau) = \mathcal{S}_{d,c,r}(\tau_1, \tau_2, \ldots, \tau_t)$.

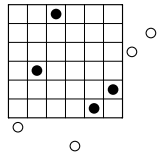

FIGURE 3. The partial permutation $(3, \_, \_, 2, 6, 5) \in \mathcal{S}_{4,2,2}$ is extendably 123-avoiding. The hollow circles indicate one possible way of extending it to a 123-avoiding permutation.

The multivariate generating function corresponding to $\mathcal{S}_{d,c,r}(\tau)$ is

$$\mathrm{Ext}_\tau(x, y, z) = \sum_d \sum_c \sum_r |\mathcal{S}_{d,c,r}(\tau)| x^d y^c z^r.$$

One important difference between normal pattern avoidance and the extended version is the amount of symmetry. Whereas the former has the eightfold symmetry of a square, in the latter only the symmetry of transposition remains.

**Lemma 3.3.** *Let $\tau$ be a set of patterns, then*

$$|\mathcal{S}_{d,c,r}(\tau)| = |\mathcal{S}_{d,r,c}(\tau^{-1})|.$$

*Proof.* It suffice to show that $|\mathcal{S}_{d,c,r}(\tau)| \geqslant |\mathcal{S}_{d,r,c}(\tau^{-1})|$, since the converse then follows from symmetry. Let $\rho \in \mathcal{S}_{d,c,r}(\tau)$, hence there exists a $\sigma \in \mathcal{S}_{d+c+r}(\tau)$, which extends $\rho$. But then is $\sigma^{-1}$ an extension of $\rho^t$ (the transpose of the matrix) and since $\sigma^{-1} \in \mathcal{S}_{d+c+r}(\tau^{-1})$, it follows that $\rho^t \in \mathcal{S}_{d,r,c}(\tau^{-1})$. $\square$

We call $\tau$ *extendably symmetric* if $|\mathcal{S}_{d,c,r}(\tau)| = |\mathcal{S}_{d,r,c}(\tau)|$. Examples of nontrivial (i.e., where $\tau \neq \tau^{-1}$) extendably symmetric pattern sets are 231 (Proposition 4.3), $\{123, 231\}$ (Proposition 5.1), $\{123, 132, 231\}$ (Proposition 5.2), $\{231, 4123\}$ (Proposition 6.7) and $\{231, 1234\}$ (Proposition 6.10).

### 3.1. Generating trees and graphs from extended pattern avoidance

To count the number of partial permutations extendably avoiding a set of patterns, $\tau$, we use algorithm 2.7, adapted to the situation. The first step is to determine how to create the unlabelled generating tree, which means we need a way to assign a parent to each nonzero partial permutation. There are again essentially four different ways of doing this, each corresponding to a direction in the plane. Due to the special asymmetry of the problem, it is preferable to build towards the right or bottom, and these in turn are interchangeable in the sense that using the former is equivalent to using the latter on the inverse pattern, with rows and columns changing roles. We therefore in this paper focus solely on building towards the right.

**Definition 3.4** (Standard parent rule). Let $\rho$ be a nonzero partial permutation, then the *parent* of $\rho$ is constructed by

- removing the lowest empty row if there are any empty rows.
- removing the rightmost column if there are no empty rows and the rightmost column is empty.
- removing the row and column of the rightmost dot, otherwise.

The unlabelled generating tree obtained by this rule is a weighted (or coloured tree, since it is important to keep track of which kind of addition was made when going down the tree. We use a solid edge for an added dot, a dashed for a new empty column and a dotted for an empty row. Alternatively, we can put the weight $x$, $y$ or $z$ on the edges to indicate which variable should be used for the generating function of the corresponding automata.

Note that it is of course possible to use the highest empty row instead of the lowest, and this give another parent rule. However, in many cases this would give the same graph, and in the other cases it is only a relatively minor change in the dotted part of the graph, so one rarely needs to consider this variation. In fact, we use it only once, namely in the proof of Proposition 6.1.

In the same way as for regular pattern avoidance, a *site*, i.e., a position where a dot could be inserted, is *active* if the resulting partial permutation is extendably avoiding the pattern. The site is *r-active* if inserting an empty row at the site gives a permutation which extendably avoids the pattern. We also say that the partial permutation is *c-active* if adding an empty column makes it extendably avoid the pattern. A quick observation gives us the following lemma.

**Lemma 3.5.** *Let* $\rho \in \mathcal{S}_{d,c,r}(\tau)$ *for a set of patterns* $\tau$. *Then*

(i) *If a site is active in* $\rho$ *then it is also r-active in* $\rho$.
(ii) *If the bottom site is active* $\rho$ *then* $\rho$ *is c-active.*

*Proof.* If adding a dot at a site is allowed, then an empty row added at the same site can be extended using the same permutation. Similarly, if we can add a dot at the bottom site, then an empty column added to the right can be extended using the same permutation. □

From the weighted tree, we can create a weighted generating graph in exactly the same way as before, by labelling each node with its subtree and putting isomorphic subtrees into equivalence classes. Since we will use the same kind of graphs many times, it is important be consistent. These are the rules we follow when drawing the graphs:

(1) Solid edges correspond to an added dot, dashed edges to an added empty column and dotted edges to an added empty row.
(2) The class containing the zero partial permutation is always named $a_1$. This is where all walks in the graph start.

(3) All partial permutations containing an empty row are in classes named $r_k$. Note that from the parent rule, any edges from an $r_k$ class will go to an $r_k$ class, so that the dotted part of the graph has no outgoing edges.

(4) All class names have an index which is the number of r-active sites, which is also the same as the number of active sites for partial permutations without empty rows.

(5) If there is a node with no outgoing edges, we will not draw it, instead, for edges to this node, we draw a short arrow, blocked by a bar. See for example the graph for the pattern set $\{123, 132, 312\}$ in Proposition 5.2.

As an example, we examine the extended version of Example 2.4.

**Example 3.6.** We want to count the partial permutations extendably avoiding $\{213, 4123\}$. The generating graph in this case is

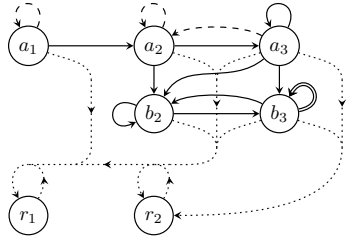

Now, let us classify the classes. Here, and in the remainder of the paper, we will not explicitly say that only the classes $r_k$ contain partial permutations with empty rows.

- $a_1$: partial permutations without dots.
- $a_2$: Increasing partial permutations for which the last column is empty or, if there is only one row, the dot is at the far right.
- $a_3$: Other increasing partial permutations.
- $b_2$: Partial permutations that are not increasing and not in class $b_3$
- $b_3$: Pre-permutations that are not increasing and the dot in the row above the rightmost dot exists and is not the 1 in a 123 pattern, where the dot representing 3 is allowed to be an empty column.
- $r_i$: Partial permutations which have empty rows and $i$ r-active sites.

In Figure 4 these classes are illustrated. From this we can see that the graph is correct. Since the graph is finite it is easy to calculate the generating function, using either the transfer matrix method or by solving the corresponding system of linear equations.

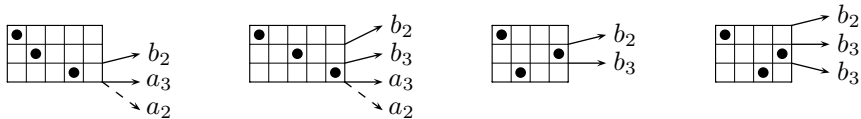

FIGURE 4. Illustration of Example 3.6. The permutations are examples of partial permutations in the classes $a_2$, $a_3$, $b_2$ and $b_3$.

The dotted part of the graph often has a very special structure, and the reason for this is that in many cases the empty rows can be added independently from one another.

**Lemma 3.7.** *If the patterns in $\tau$ are such that no two empty rows added to any partial permutation, when extended increasingly, can both be used to get an occurrence of a pattern, then the dotted part of the graph has the structure below. The same hold if "increasingly" is changed to "decreasingly".*

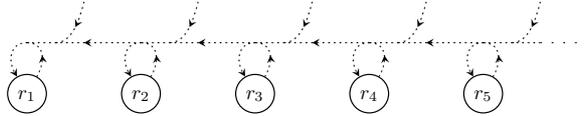

*Proof.* We will show that whenever a site is r-active, then itself and all r-active sites below will remain r-active when inserting an empty row there and hence the graph has the desired structure. Let $\rho \in \mathcal{S}_{d,r,c}(\tau)$, $i$ and $j$ be r-active sites in $\rho$ such that $i \leqslant j$ and $\tilde{\rho}$ be the partial permutation we get when adding empty rows to $\rho$ at both of these sites. By the condition, $\rho$ can be extended to a permutation for which all the dots to the right of the partial permutation are increasing. Hence, by using the same extension for $\tilde{\rho}$, except that there are more increasing dots to the right, we see that this extended permutation is $\tau$-avoiding, since otherwise the condition in the lemma would fail. Therefore, by the standard parent rule, the site $j$ remains r-active when inserting an empty row at $i$. The second statement is analogous. $\square$

The notions of Wilf-equivalence and graph-equivalence have natural generalisations to extended pattern avoidance.

**Definition 3.8.** Two sets of patterns, $\tau$ and $\tau'$, are *extendably Wilf-equivalent* if $|\mathcal{S}_{d,c,r}(\tau)| = |\mathcal{S}_{d,c,r}(\tau')|$ for all nonnegative integers, $d$, $c$ and $r$.

**Definition 3.9.** Two sets of patterns, $\tau$ and $\tau'$, are *extendably graph-equivalent* if there exists parent rules, $\mathsf{p}$ and $\mathsf{p}'$, such that the corresponding generating graphs are equal.

Exactly as before, if two patterns are extendably graph-equivalent, then there is a natural bijection since both of them are in bijection with the walks on the graph. We can also note that the generating graph algorithms also work for extended avoidance. It is, however, a bit more complicated to implement in this case due to the additional problem of checking whether a partial permutation can be extended to a $\tau$-avoiding permutation. This can by done by examining all $\tau$-avoiding permutations from which we get the $\tau$-avoiding partial permutations as all the possible upper left corners. Thereafter, we must make sure to remove duplicates, since the partial permutations may have several different extensions.

## 4. Singleton patterns of length three

The partial permutations avoiding singleton patterns of length three was enumerated by Eriksson and Linusson [4], for the pattern 321, and Linusson [5] for the remaining, by essentially treating all the patterns separately and using clever tricks to find the solutions. With the help of the automata we have been able to do a unified treatment of all cases. As a side benefit, we also get three bijections, due to graph equivalences.

## 4.1. 123-avoidance and 213-avoidance

**Proposition 4.1.** *The patterns* 123 *and* 213 *are extendably graph-equivalent and have the generating function*

$$\text{Ext}_{123} = \text{Ext}_{213} = \frac{1}{(1-y)(1-z)-x}\left(1 + \frac{x^2 c(x)}{(1-y-xc(x))(1-z-xc(x))}\right),$$

*where* $c(x) = \frac{1-\sqrt{1-4x}}{2x}$ *is the generating function for the Catalan numbers.*

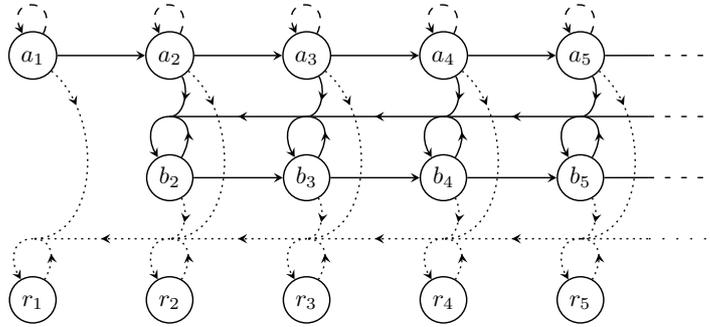

FIGURE 5. Generating graph for 123.

*Proof.* First we consider extendably 123-avoiding partial permutations, where we have the following classes:

- $a_k$ is the class of decreasing partial permutations with $k-1$ dots.
- $b_k$ is the class of partial permutations which have its first ascent from row $k-1$ to row $k$.
- $r_k$ is the class of partial permutations with empty rows and $k$ r-active sites.

By studying the properties of these classes we show that the generating graph is as in Figure 5.

A partial permutation, $\sigma \in a_k$, has all the $k$ sites active and is $c$-active. From the descriptions of the classes it is clear that adding a dot to the top site give a partial permutation in class $a_{k+1}$ and adding a dot in site $m$ from the top, with $m \geqslant 2$, gives a partial permutation in class $b_m$. Adding an empty column will not change the class of $\sigma$.

The class $b_k$ has the top $k$ sites active, and is not c-active. Adding a dot at the top site creates a partial permutation in $b_{k+1}$, and adding the dots at the site $m \geqslant 2$ gives a partial permutation of class $b_m$.

Finally, we need to consider adding empty rows. Since the empty rows can be extended to the right in decreasing order, we can conclude by Lemma 3.7, that the empty rows can be added independently of each other, as two empty rows can never both be used to get the 123 pattern.

Since none of the classes have the same number of outgoing edges of all three different types, we can conclude that none of the classes are equivalent and therefore the graph in Figure 5 is the generating graph for extendably 123-avoiding partial permutations.

For the extendably 213-avoiding partial permutations, the classes have a very similar description:

- $a_k$ is the class of increasing partial permutations with $k - 1$ dots.
- $b_k$ is the class of partial permutations which have its first descent from row $k$ to row $k - 1$.
- $r_k$ is the class of partial permutations with empty rows and $k$ r-active sites.

From this it is easy to see, by rehashing the arguments for the pattern 123, that 123 and 213 are indeed extendably graph-equivalent. This implies that there is a bijection induced by the graph-equivalence and in Figure 6 there is an example of this bijection. Finally, the generating function can now be calculated, using Lemmas 2.9 and 2.10. $\qquad\square$

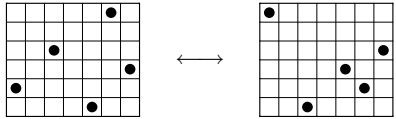

FIGURE 6. Example of the bijection between $\mathcal{S}_{d,c,r}(123)$ and $\mathcal{S}_{d,c,r}(213)$, induced by the graph-equivalence.

## 4.2. 132-avoidance and 312-avoidance

**Proposition 4.2.** *The patterns* 132 *and* 312 *are extendably graph-equivalent, with generating function*

$$\mathrm{Ext}_{132} = \mathrm{Ext}_{312} = \frac{1 - xc(x)}{(1 - y - xc(x))(1 - z - xc(x))} = \frac{c(x)}{(1 - yc(x))(1 - zc(x))}.$$

*Proof.* For the case of extendably 132-avoiding partial permutations, we get the following classes:

- The class $a_1$ contains the zero permutation and partial permutations which have the last column empty.
- The class $a_k, k \geqslant 2$ consists of partial permutation whose matrices are anti-diagonal block matrices, where each block is a 132-avoiding permutation. The blocks are separated by one or more empty columns and there may be empty columns to left as well. The last block has $k$ active sites. See Figure 8 for an example.

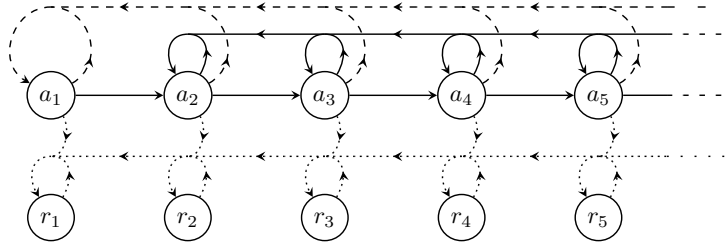

FIGURE 7. Generating graph for 132 and 312.

- $r_k$ is the class of partial permutations with empty rows and $k$ r-active sites.

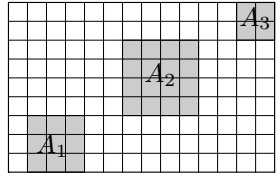

FIGURE 8. Example of a partial permutation in class $a_k$, $k \geqslant 2$. The blocks $A_1$, $A_2$ and $A_3$ are 132-avoiding permutations.

It is clear from the class description that partial permutations in $a_1$ only have the top most site active and the partial permutations in the $a_k$-classes are c-active. For the classes $a_k$, $k \geqslant 2$, the active sites are all within the top most block of the block matrix, so the solid part of the graph is the Catalan graph. Again the empty rows can be added independently, so we get the desired generating graph.

The extendably 312-avoiding case is very similar, except that the block matrix in the class $a_k$ is diagonal, instead of anti-diagonal, with 312-avoiding permutations as blocks. The generating graph is therefore the same in this case.

The generating function for this graph is easy to calculate, since the solid part of the graph is the Catalan graph in Example 2.5 and for the dotted part, Lemma 2.8 can be used. □

### 4.3. 231-avoidance and 321-avoidance

**Proposition 4.3.** *The patterns* 231 *and* 321 *are extendably graph equivalent and have the generating function*

$$\mathrm{Ext}_{231} = \mathrm{Ext}_{321} = \frac{1 - xc(x)}{(1 - y - xc(x))(1 - z - xc(x))} = \frac{c(x)}{(1 - yc(x))(1 - zc(x))}.$$

*Proof.* We begin by showing that the pattern 321 has the generating graph show in Figure 9. The classes are as follows:

- The class $a_k$ contains 321-avoiding permutations with $k$ active sites.

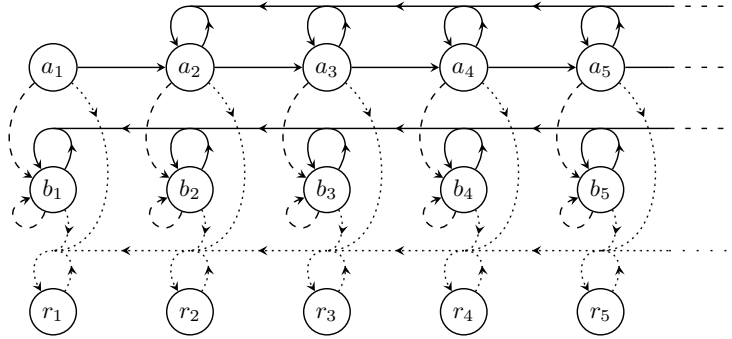

FIGURE 9. Generating graph for 231 and 321.

- The class $b_k$ consists of partial permutations with at least one empty column and $k$ active sites.
- $r_k$ is the class of partial permutations with empty rows and $k$ r-active sites.

Since the classes $a_k$ contain the regular 321-avoiding permutations, we immediately get that this part of the graph is the Catalan graph. When adding an empty column, the number of active sites remains the same, but the part to the right of the first empty column must be increasing, since any 21-pattern would give an 321-pattern when extending the empty column. This means that the active sites for the class $b_k$ are precisely the $k$ bottom sites and that no new sites can be created. The empty rows can clearly be added independently, so we get the desired graph.

For the extendably 231-avoiding partial permutations, we have the same descriptions of the classes, except $a_k$ contains 231-avoiding permutations. This time, however, the active sites are not at the bottom, but this does not affect the graph structure. For the classes $b_k$, we can say that whenever we have an empty column, all dots to its right must be decreasing to avoid the 231 pattern, which means only the sites above a newly inserted dot will remain active. Finally, the empty rows can again be added independently, so we get the graph in Figure 9.

The generating function for 231 must be the same as for the pattern 312, since the patterns are inverses. □

Combining Propositions 4.1, 4.2 and 4.3, we get

**Corollary 4.4** (Linusson). *The patterns 132, 312, 231 and 321 are extendably Wilf-equivalent and the patterns 123 and 213 are extendably Wilf-equivalent.*

## 5. Combinations of patterns of length three

In general, taking more patterns of length three makes the task easier and most of the cases have a finite generating graph. Using the method described before, the problem of finding the generating functions can be solved more or less automatically, so we will only give explicit proofs in a the few cases where there are non-trivial

Wilf-equivalences to be found. The generating graphs and generating functions for all combinations can be found in the appendix.

In many instances from the remaining two sections, we have classes which only differ by some empty columns at the left or right, so we will introduce a notation for these. Suppose $\rho \in a$, then $\rho' \in a^*$ iff $\rho'$ has the same word form as $\rho$ and $\rho'$ has more empty columns than $\rho$. Further, $\rho' \in {}^*a$ iff $(\rho')^c$ and $\rho^c$ have the same word form and $\rho'$ has more empty columns than $\rho$.

**Proposition 5.1.** *The pattern sets* $\{123, 231\}$ *and* $\{123, 312\}$ *are extendably Wilf-equivalent, hence they are both extendably symmetric.*

*Proof.* We choose to study the pattern set $\{123, 312\}$, since even though $\{123, 231\}$ has a finite generating graph, it is more complicated to describe the classes.

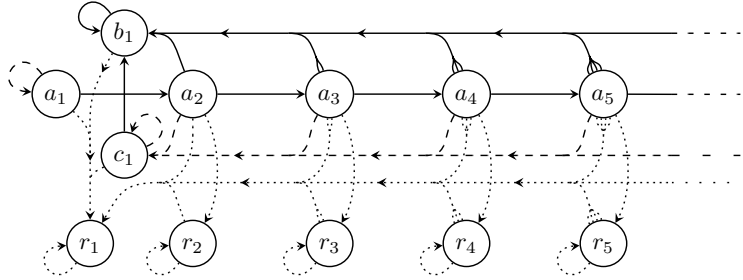

The classes are:

- $a_i$: Decreasing partial permutations with $i-1$ dots, such that any empty columns are to the left.
- $b_1$: Non-decreasing partial permutations.
- $c_1$: Decreasing partial permutations not in class $a_i$ for all $i$.
- $r_i$: Partial permutations which have empty rows and $i$ r-active sites.

Here we see that for partial permutations in $a_i$, all sites are active and adding a dot at the top site gives us a partial permutation in $a_{i+1}$. The other sites make the partial permutation non-decreasing so it would ends up in class $b_1$. A new column would bring us to class $c_1$. Both $b_1$ and $c_1$ consists of partial permutations with exactly one active site, namely the site right above the rightmost dot for $b_1$ and the bottom site for $c_1$. Thus $c_1$ is c-active and adding an empty column clearly does not change the class.

For the dotted part of the graph, we only need to analyse the $a_i$-classes, since both $b_1$ and $c_1$ has only one r-active site, where we can add any number of empty rows. For $a_i$, it is clear that if we add empty rows at two different sites, one of them must be the top site. This immediately gives the desired graph.

The generating function of the automaton,

$$\frac{1}{1-y}\frac{1}{1-z}\left(\frac{1}{1-x}+\frac{x}{(1-x)^2}\left(\frac{z}{1-z}+\frac{y}{1-y}\right)+\frac{x^2}{(1-x)^3}\right),$$

is straight forward to calculate and it is clearly symmetric in $y$ and $z$. $\qquad\square$

**Proposition 5.2.** *The pattern sets* $\{123, 132, 231\}$, $\{123, 132, 312\}$ *and* $\{213, 231,$ $312\}$ *are extendably Wilf-equivalent.*

*Proof.* The pattern set $\{123, 132, 312\}$ has a very simple graph:

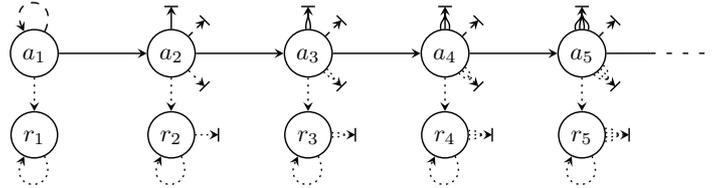

Here the class $a_i$ consists of only the decreasing partial permutations with $i - 1$ rows, such that any empty column is at the left edge. All sites are active, but every site, except for the top site, gives a permutation which has no active site. The same goes for adding an empty column or adding an empty row below the top site. Adding an empty row at the top bring us to a partial permutation in $r_i$ which has the same properties as the permutation in $a_i$, except only empty rows can be added.

The generating function for this graph is

$$\frac{1}{(1-x)^2} \frac{1}{1-y} \frac{1}{1-z} - \frac{x}{1-x},$$

which is symmetric, so $\{123, 132, 231\}$ and $\{123, 132, 312\}$ are Wilf-equivalent. Now it remains to examine the pattern set $\{213, 231, 312\}$, which we show has the graph

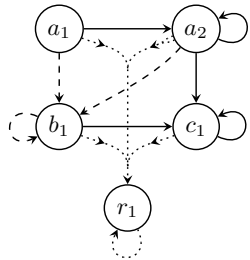

where the classes are:

- $a_1$: The size zero permutation.
- $a_2$: Nonzero increasing permutations.
- $b_1 = a_1^* \cup a_2^*$.
- $c_1$: Partial permutations which are non-increasing, or which are increasing and have an empty column to the left of a dot.
- $r_i$: Partial permutations which have empty rows and $i$ r-active sites.

It is easy to verify that the number and type of the active sites are in agreement with the graph. For example, the class $a_2$ has only the two bottom sites active, and inserting a dot would give a permutation in class $c_1$ and $a_2$ respectively.

The generating function for this graph is the same as for the two previous patterns, hence they are all Wilf-equivalent. $\qquad\square$

**Proposition 5.3.** *The pattern sets* $\{123, 231, 312\}$ *and* $\{132, 213, 321\}$ *are extendably graph-equivalent.*

*Proof.* We show that both patterns have the generating graph:

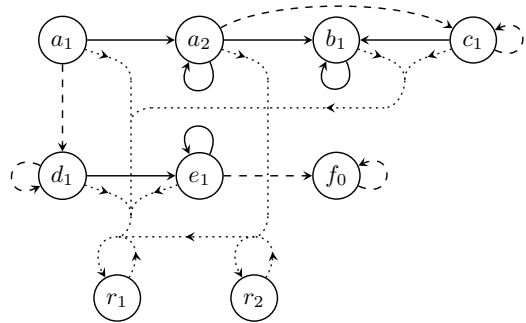

First, for the pattern set $\{123, 231, 312\}$, we have the following classes:

- $a_1$: The zero permutation.
- $a_2$: Nonzero decreasing permutations.
- $b_1$: Non-decreasing partial permutations.
- $c_1 = a_2^*$.
- $d_1 = a_1^*$.
- $e_1 = {}^*a_2$.
- $f_0 = {}^*a_2^*$.
- $r_i$: Partial permutations which have empty rows and $i$ r-active sites.

Again, it is easy to check that these descriptions give the graph above. One can notice that only nonzero decreasing permutations have more than one active site.

As for the second pattern set, $\{132, 213, 321\}$, we can see that it is the reverse of $\{123, 231, 312\}$, and it turns out that all the class descriptions happen to be same as for the first pattern set, but turned up-side down, so that increasing becomes decreasing and vice versa. $\qquad\square$

**Proposition 5.4.** *The pattern sets* $\{123, 132, 231, 312\}$ *and* $\{132, 213, 231, 312\}$ *are extendably Wilf-equivalent.*

*Proof.* First we show that the generating graph for $\{123, 132, 231, 312\}$ is

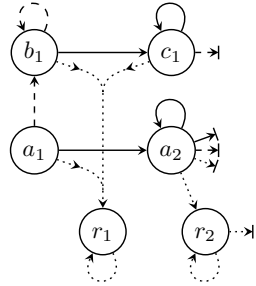

with equivalence classes

- $a_1$: The zero permutation.
- $a_2$: Nonzero decreasing permutations.
- $b_1$: Nonzero partial permutations without dots.
- $c_1 = {}^*a_2$.
- $r_i$: Partial permutations which have empty rows and $i$ r-active sites.

All classes have very simple structure so it is not hard to verify that they give rise the graph above.

Next, it is equally easy to check that the pattern set $\{132, 213, 231, 312\}$ has the graph

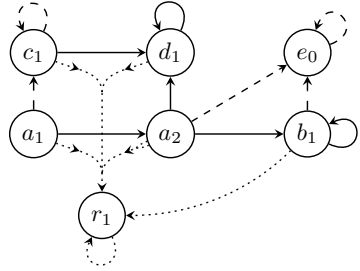

with equivalence classes

- $a_1$: The zero permutation.
- $a_2$: The permutation of size one.
- $b_1$: Increasing permutations of size two or greater.
- $c_1$: Nonzero partial permutations without dots.
- $d_1$: Decreasing permutations having size two or more or with at least one empty column added to the left.
- $e_0 = a_2^* \cup b_1^*$.
- $r_i$: Partial permutations which have empty rows and $i$ r-active sites.

The generating function for both the graph is

$$\frac{1}{1-x}\frac{1}{1-y}\frac{1}{1-z} + \frac{x}{1-x}\left(\frac{1}{1-y} + \frac{1}{1-z} - 1\right) - x,$$

so we conclude that they are Wilf-equivalent. $\qquad\square$

| | 123 | 132 | 213 | 231 | 312 | 321 |
|---|---|---|---|---|---|---|
| 1234 | - | * | $E_1$ | $G_1$ | $G_1^{-1}$ | $H_1$ |
| 1243 | - | - | * | * | * | * |
| 1324 | - | - | - | * | * | * |
| 1342 | - | - | $A_1$ | - | * | * |
| 1423 | - | - | $A_1^{-1}$ | * | - | * |
| 1432 | * | - | * | $B_2^{-1}$ | $B_2$ | - |
| 2134 | - | * | - | * | * | * |
| 2143 | * | - | - | * | * | * |
| 2314 | - | * | - | - | $C_1^{-1}$ | * |
| 2341 | - | * | $A_2$ | - | $D_1^{-1}$ | $F_1$ |
| 2413 | $A_3^{-1}$ | - | - | - | - | $B_1$ |
| 2431 | * | - | * | - | $B_3$ | - |
| 3124 | - | * | - | $C_1$ | - | * |
| 3142 | $A_3$ | - | - | - | - | $B_1^{-1}$ |
| 3214 | $E_2$ | * | - | $C_2$ | $C_2^{-1}$ | - |
| 3241 | * | * | - | - | $B_4^{-1}$ | - |
| 3412 | * | * | * | - | - | * |
| 3421 | * | * | * | - | * | - |
| 4123 | - | * | $A_2^{-1}$ | $D_1$ | - | $F_1^{-1}$ |
| 4132 | * | - | * | $B_3^{-1}$ | - | - |
| 4213 | * | * | - | $B_4$ | - | - |
| 4231 | * | * | * | - | - | - |
| 4312 | * | * | * | * | - | - |
| 4321 | $H_2$ | * | * | $F_2$ | $F_2^{-1}$ | - |

FIGURE 10. Relations between different pattern-pairs. See below for an explanation of the entries.

## 6. Partial permutations extendably avoiding one pattern of length three and one pattern of length four

The number of different possible combinations of one pattern of length three and one pattern of length four is very large, so we concentrate on cases where non-trivial symmetries are found. Figure 10 is a table of all these combinations and the relations among them. The dash indicates that this pattern pair isn't really a pair, since the smaller pattern is contained in the larger. Entries with the same letters are cases which are non-trivially Wilf-equivalent to each other, so these are the ones we will study more closely in this section. Two pattern sets which are inverses to each other are marked with the same letter and index, but with an inverse sign on one of them. Finally, the rest, which lack interesting Wilf-equivalences, are marked with asterisks.

In the following seven subsections we prove all the non-trivial identities, where each subsection contains the identities indicated by the corresponding letter in Figure 10.

### 6.1. Case A: The pattern sets $\{213, 1423\}$, $\{213, 4123\}$ and $\{123, 2413\}$

**Proposition 6.1.** *The pattern sets* $\{213, 1423\}$, $\{213, 4123\}$ *and* $\{123, 2413\}$ *are extendably Wilf-equivalent, and their generating function is*

$$\frac{1}{(1-y)(1-z)} \left( 1 + \frac{x(1-x)}{(1-x-y)(1-3x+x^2)} \left( \frac{1}{1-z} - x \right) \right).$$

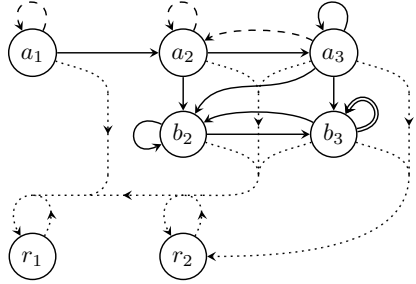

FIGURE 11. Generating graph for $\{213, 4123\}$ and $\{213, 1423\}$.

*Proof.* The first equality we prove bijectively by showing that they are graph-equivalent. This time, however, we must use the alternative parent rule, discussed after Definition 3.4, for the pattern set $\{213, 1423\}$.

For the pattern set $\{213, 1423\}$, with the alternative parent rule, we get the following equivalence classes:

- $a_1$: Partial permutations without dots.
- $a_2$: Partial permutations that are increasing and where the last column is empty or there is only one dot.
- $a_3$: Other increasing partial permutations.
- $b_2$: Non-increasing partial permutations that have a dot in the NE corner.
- $b_3$: Non-increasing partial permutations which are not in class $b_2$.
- $r_i$: Partial permutations which have empty rows and $i$ r-active sites.

The pattern set $\{213, 4123\}$ was studied in Example 3.6.

The last case, the pattern set $\{123, 2413\}$, has these classes:

- $a_k$: Decreasing partial permutations such that $k = m + 2$ if the number of rows is greater than $m$, the number of consecutive nonempty columns at the right edge, and $k = m + 1$ otherwise, i.e., if all empty columns are at the left edge.

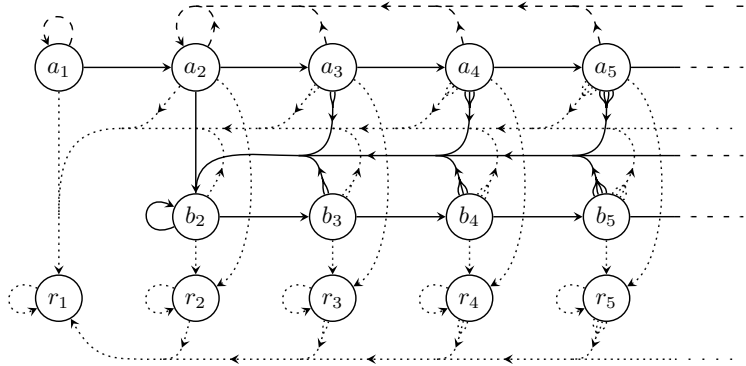

Figure 12. Generating graph for $\{123, 2413\}$.

- $b_k$: Partial permutations which have a 12-pattern and where $k$ is the largest value such that the NE corner has a $(k-2)(k-1)\cdots 1$-pattern.
- $r_i$: Partial permutations with empty rows and $i$ r-active sites.

Again, it is not difficult to see that the given classes in the two cases give rise to the graphs in Figure 11 and 12, respectively.

The generated function for the automaton in Figure 12 can be calculated using

$$a_1 = 1/(1-y), a_i = x^{i-2}a_2 \text{ and } b_i = x^{i-2}b_2, \text{ for } i \geqslant 2,$$

$$a_2 = xa_1 + y\sum_{i\geqslant 2} a_i = \frac{x}{1-y} + \frac{ya_2}{1-x} = \frac{x(1-x)}{(1-y)(1-x-y)},$$

$$b_2 = x\sum_{i\geqslant 2}(i-1)(a_i + b_i) = \frac{x^2(1-x)}{(1-y)(1-x-y)(1-3x+x^2)}.$$

$\square$

## 6.2. Case B: $\{321, 2413\}$, $\{312, 1432\}$, $\{312, 2431\}$ and $\{231, 4213\}$

We show that these four sets of patterns are extendably Wilf-equivalent by three lemmas. First we show that the first and last pattern sets have the same graph. The lemma thereafter shows that the inverse of the first and third pattern sets are extendably graph-equivalent. Finally, the second pattern set is analysed and shown to be Wilf-equivalent to the other three.

**Lemma 6.2.** *The pattern sets $\{321, 2413\}$ and $\{231, 4213\}$ are extendably graph-equivalent.*

*Proof.* First, for $\{321, 2413\}$, we have that the bottom $k$ sites are active (or r-active) for all classes. The descriptions are

- $a_1$: The zero permutation.
- $a_k, k \geqslant 2$: Permutations with the pattern $12\cdots(k-1)$ in the SE corner.
- $b_k$: Permutations, such that the dot in the last column is in row $n-k+1$.

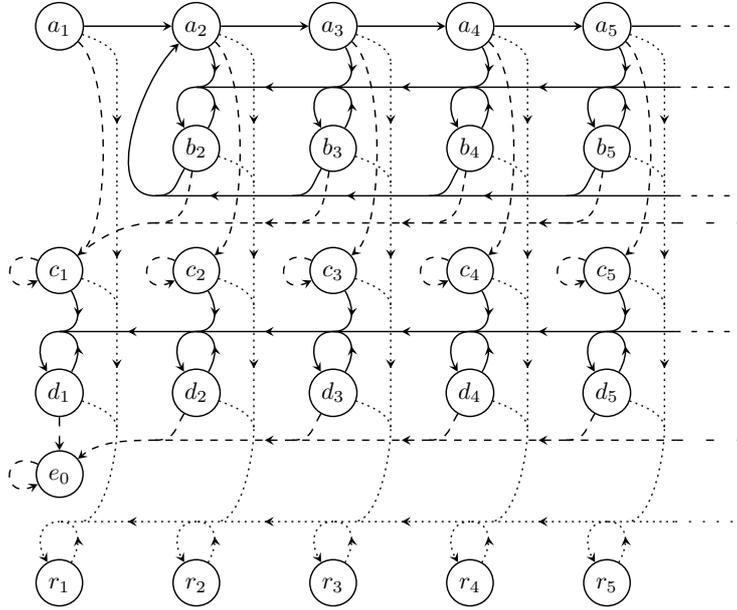

FIGURE 13. Generating graph for $\{321, 2413\}$ and $\{231, 4213\}$.

- $c_1 = a_1^* \cup b_2^* \cup b_3^* \cup \cdots$.
- $c_k = a_k^*, k \geqslant 2$.
- $d_k$: Partial permutations with at least one empty column. The last column is not empty and it has a dot in row $n - k + 1$.
- $e_0 = \bigcup_{k \geqslant 1} d_k^*$.
- $r_i$: Partial permutations which have empty rows and $i$ r-active sites.

For the pattern set $\{231, 4213\}$, the description is similar, but the active sites are not necessarily at the bottom. The classes with different descriptions are:

- $b_k$: Permutations not of class $a_i$ for any $i$. Let $m$ be the row of the dot in the last column. If the top $m$ rows are increasing, then $k = m + 1$, otherwise, $k$ is the largest value, such that the rows $(m - k + 1)(m - k + 2) \cdots m$ are increasing.
- $d_k$: Partial permutations with $c \geqslant 1$ and the last column nonempty, which have $k$ active sites, namely the $k$ sites directly above the dot in the last column.

By analysing these classes, one can see that we get that both pattern sets have the graph in Figure 13 as generating graphs. $\qquad \square$

**Lemma 6.3.** *The pattern sets $\{321, 3142\}$ and $\{231, 4132\}$ are extendably graph-equivalent.*

*Proof.* Here both cases have the same class descriptions, but the placements of the active sites differ.

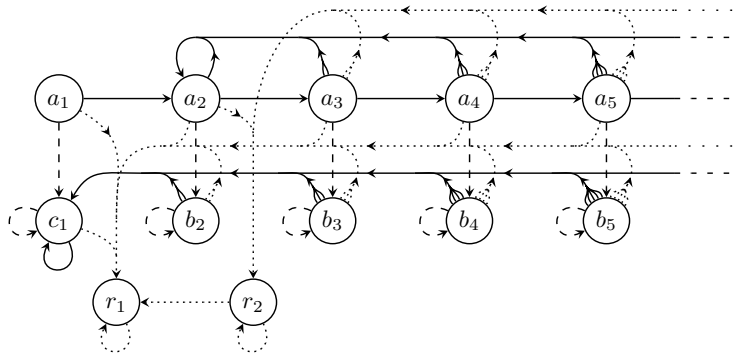

FIGURE 14. Generating graph for $\{231, 4132\}$ and $\{321, 3142\}$.

- $a_1$: The zero permutation.
- $a_k, k \geqslant 2$: Nonempty permutations where $m$ is the largest number such that there is a $12 \cdots m$-pattern in the SE corner, where $k = m + 1$ if the permutation is increasing and $k = m + 2$ otherwise.
- $b_k = a_k^*$.
- $c_1$: Nonzero partial permutations with one active site.
- $r_i$: Has empty rows and $i$ r-active sites.

From these descriptions, we can infer that both have the same generating graph, namely the one shown in Figure 14. □

**Lemma 6.4.** *The pattern set $\{312, 1432\}$ has the generating graph shown in Figure 15.*

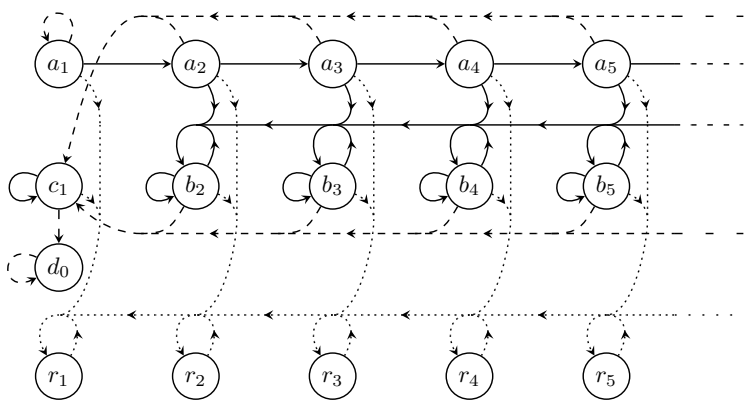

FIGURE 15. Generating graph for $\{312, 1432\}$.

*Proof.* The classes are as follows.

- $a_k$: Decreasing partial permutations such that any empty column is at the left edge.
- $b_k$: Partial permutations ending with $m\,i(i-1)\cdots(i-k+3)$, where $0 \leqslant k-2 < i+1 < m \leqslant n$ and all empty columns are at the left edge.
- $c_1$: Partial permutations $a_k$ or $b_k$, $k \geqslant 2$, with one empty column and a $12\cdots i$-pattern to the right of the empty column added to the SE corner, where $i \geqslant 0$.
- $d_0 = c_1^*$.
- $r_i$: Has empty rows and $i$ r-active sites.

From this we can conclude that the graph is the one in Figure 15. $\qquad\square$

**Corollary 6.5.** *The pattern sets* $\{321, 2413\}$, $\{312, 1432\}$, $\{312, 2431\}$ *and* $\{231, 4213\}$ *are all extendably Wilf-equivalent, with generating function*

$$\frac{1}{1-z} + \frac{(1-x)^2}{(1-x-z)(1-3x+x^2)}\left(\frac{x}{1-z} + \frac{y}{1-y}\right) + \frac{y^2}{(1-y)^2}\frac{x}{1-3x+x^2}.$$

*Proof.* What remains is to show that both $\{321, 2413\}$ and $\{312, 1432\}$ have the desired generating function. For the pattern set $\{321, 2413\}$, this is a fairly straight forward application of Lemmas 2.8 and 2.9, even though the graph looks a bit complicated. For the pattern set $\{312, 1432\}$, Lemma 2.8 is not directly applicable, since there are loops on the $b_k$ nodes. However, by following the method of the proof of that lemma, we get that $b_k = \frac{1}{1-y}\frac{x^k}{1-3x+x^2}$, from which the generated function is easily verified. $\qquad\square$

### 6.3. Case C: The pattern sets $\{312, 2314\}$ and $\{312, 3214\}$

**Proposition 6.6.** *The pattern sets* $\{312, 2314\}$ *and* $\{312, 3214\}$ *are Wilf-equivalent, and their generating function is*

$$\frac{1}{(1-z)(1-3x+x^2)(1-2x-y+xy)}\left(1-5x+6x^2+\frac{x(1-2x)}{1-z}\right).$$

*Proof.* First, for the pattern set $\{312, 2314\}$, we have the following classes:

- $a_1$: The zero permutation or the (normally, not necessarily extendably) 231-avoiding pre-permutations such that the last column is empty.
- $a_k$: Partial permutations which are 231-avoiding and end with $n(n-1)\cdots(n-k+2)$.
- $b_k$: Partial permutations that have the pattern $2(k+1)k\cdots31$ in consecutive rows, where the $(k+1)$ is in the last column.
- $r_i$: Partial permutations with at least one empty row and $i$ sites r-active.

From this we can see that we get the graph shown in Figure 16. The pattern set $\{312, 3214\}$ has a finite generating graph which is remarkably similar to the one for the pattern set $\{213, 4123\}$, as seen in the Figures 11 and 17.

- $a_1$: The zero permutation or 321-avoiding partial permutations where the last column is empty.

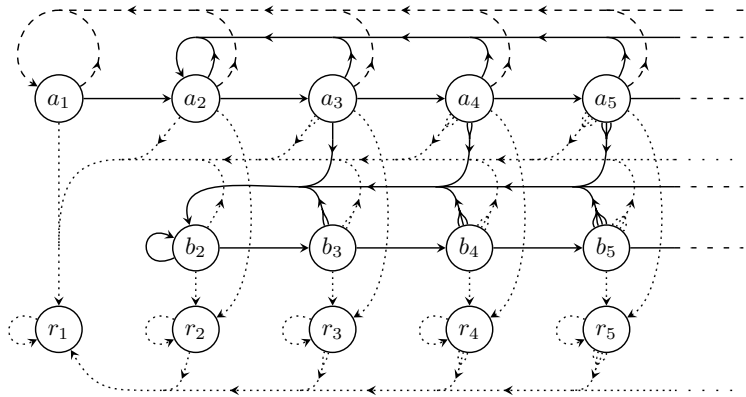

FIGURE 16. Generating graph for $\{312, 2314\}$.

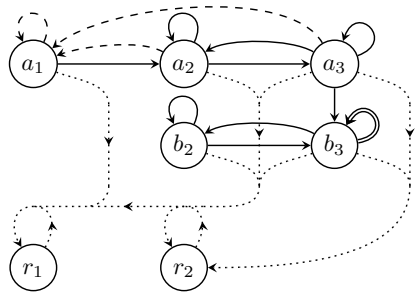

FIGURE 17. Generating graph for $\{312, 3214\}$.

- $a_2$: Partial permutations with one row and a dot in the last column or 321-avoiding partial permutations for which the dots in the last two rows are increasing.
- $a_3$: 321-avoiding partial permutations which do not belong to the previous two classes.
- $b_i$: Partial permutations which contain a 321-pattern and have $i$ active sites. It has two active sites iff the dot in the row after the rightmost dot is the 1 in a 321-pattern.
- $r_i$: The class $r_i$ consists as usual of partial permutations with at least one empty row and $i$ sites r-active.

The descriptions give the necessary information to conclude that the graph are as asserted. We can also easily calculate the generating functions from these graphs to find that they indeed agree and are as stated. □

### 6.4. **Case D: The pattern set** $\{231, 4123\}$

**Proposition 6.7.** *The pattern set* $\{231, 4123\}$ *is extendably symmetric with generating function*

$$\frac{1}{1 - 4x + 5x^2 - 3x^3} \left( x^3 + \frac{x^2(y^2(1-3x) + 3xy - x)}{(1-y)^3} + \frac{x^2(z^2(1-3x) + 3xz - x)}{(1-z)^3} \right.$$
$$\left. + \frac{1}{(1-y)(1-z)} \left( 1 - 5x + 5x^2 + x(1-x) \left( \frac{1}{1-y} + \frac{1}{1-z} \right) \right) \right)$$

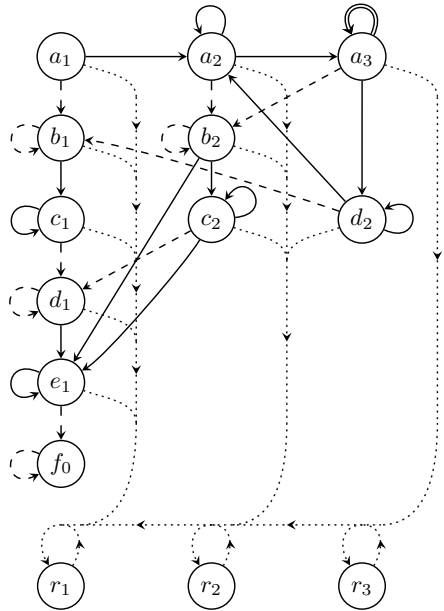

FIGURE 18. Generating graph for $\{231, 4123\}$.

*Proof.* The classes can be described as follows:

- $a_1$: The zero permutation.
- $a_2$: Permutations, not in class $a_3$, ending with $n(n-1)\cdots(l+1)$, where $1 \leqslant l < n$.
- $a_3$: Permutations ending with $l(l-1)\cdots(k+1)\, 1n(n-1)\cdots(l+1)$, where $1 \leqslant k < l < n$.
- $b_1 = a_1^* \cup d_2^*$.
- $b_2 = a_2^* \cup a_3^*$.
- $c_1$: Partial permutations of the form $(b_1)\, n(n-1)\cdots(k+1)$, where $1 \leqslant k < n$.
- $c_2$: Partial permutations of the form $(b_2)\, n(n-1)\cdots(k+1)$, where $1 < k < n$.

- $d_1 = c_1^* \cup c_2^*$.
- $d_2$: Permutations ending with $n(n-1)\cdots(l+1)\,k(k-1)\cdots(j+1)\,l(l-1)\cdots(k+1)$, where $1 \leqslant j < k < l < n$.
- $e_1$: Partial permutations in either $b_2$, $c_2$ or $d_1$, with at least one dot added and all dots added at a single site of the initial partial permutation.
- $f_0 = e_1^*$.
- $r_i$: Partial permutations with at least one empty row and $i$ sites r-active.

From these descriptions we can conclude that the graph is as in Figure 18. Since the generating graph is finite, we can use the transfer matrix method to find the generating function, which is clearly symmetric in $y$ and $z$. □

### 6.5. **Case E: The pattern sets** $\{213, 1234\}$ **and** $\{123, 3214\}$

**Proposition 6.8.** *The pattern sets* $\{213, 1234\}$ *and* $\{123, 3214\}$ *are extendably graph-equivalent and their generating function is*

$$\frac{1}{(1-y)(1-z)} + \frac{x}{(1-y)^2(1-z)^2(1-3x+x^2)}\left(1 - 2x + \frac{x(1-xyz)}{(1-y)(1-z)}\right).$$

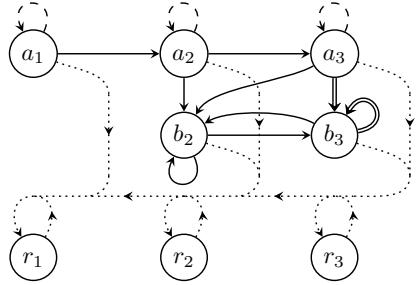

FIGURE 19. Generating graph for $\{123, 3214\}$ and $\{213, 1234\}$.

*Proof.* In both cases only the top three sites can be active. The classes for the pattern set $\{213, 1234\}$ are:

- $a_1$: Partial permutations without dots.
- $a_2$: Partial permutations with exactly one dot.
- $a_3$: Partial permutations with two dots which are in increasing order.
- $b_2$: Partial permutations with two or more dots which have a dot in the NE corner.
- $b_3$: Partial permutations with three or more dots but no dot in the NE corner.
- $r_i$: Partial permutations with at least one empty row and $i$ sites r-active.

And for $\{123, 3214\}$, we get:

- $a_1$: Partial permutations without dots.
- $a_2$: Partial permutations with exactly one dot.
- $a_3$: Partial permutations with two dots in decreasing order.

- $b_2$: Partial permutations with two or more dots for which there is an ascent in the top two rows.
- $b_3$: Partial permutations with three or more dots that do not belong to $b_2$.
- $r_i$: Partial permutations with at least one empty row and $i$ sites r-active.

Using this, we can see that both patterns give the graph in Figure 19. □

### 6.6. Case F: The pattern sets $\{312, 4321\}$ and $\{321, 4123\}$

**Proposition 6.9.** *The pattern sets $\{312, 4321\}$ and $\{321, 4123\}$ are extendably Wilf-equivalent and their generating function is*

$$\frac{1}{1-z} + \frac{1}{1-3x+x^2} \left( y^2 \left( \frac{1}{1-y} + \frac{1}{1-z} - 1 \right) + \frac{y(1-2x)}{1-z} \right.$$
$$\left. + \frac{x(1-2x+y)}{(1-z)^2} + \frac{x^2}{(1-z)^3} \right).$$

*Proof.* These two cases are not graph-equivalent, but their graphs are very similar. In both cases only the bottom three sites are active.

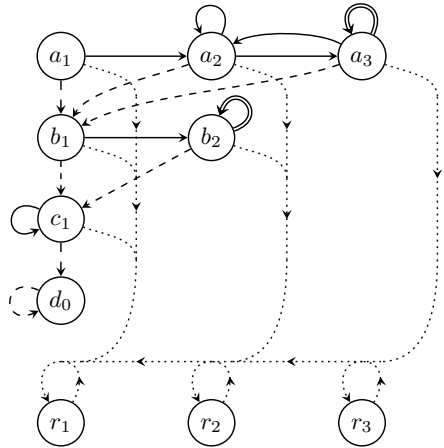

FIGURE 20. Generating graph for $\{312, 4321\}$.

First for $\{312, 4321\}$, we have

- $a_1$: The zero permutation.
- $a_2$: Permutations with a dot in the SE corner.
- $a_3$: Permutations where the last two rows have a 21-pattern.
- $b_1$: Partial permutations with one empty column, which is at the right edge.
- $b_2$: Partial permutations with one empty column not at the right edge.
- $c_1$: Partial permutations with two empty columns.
- $d_0$: Partial permutations with at least three empty columns.

• $r_i$: Partial permutations with at least one empty row and $i$ sites r-active.

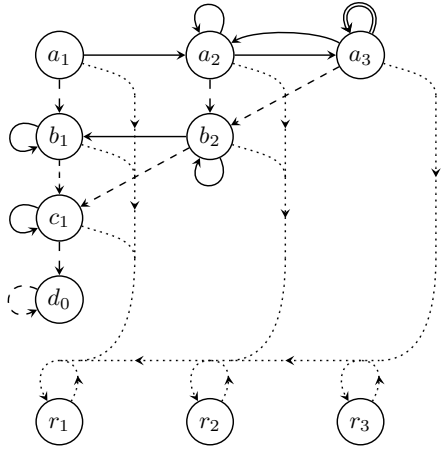

FIGURE 21. Generating graph for $\{321, 4123\}$.

The pattern set $\{321, 4123\}$ has the following classes:

• $a_1$: The zero permutation.
• $a_2$: The permutation of size one or permutations where the two last rows have a 21-pattern.
• $a_3$: Permutations where the last two rows have a 12-pattern.
• $b_1$: Partial permutations with one empty column and a dot in the SE corner if it has any dots.
• $b_2$: Partial permutations with one empty column but not in class $b_1$.
• $c_1$: Partial permutations with two empty columns.
• $d_0$: Partial permutations with at least three empty columns.
• $r_i$: Partial permutations with at least one empty row and $i$ sites r-active.

From the description, it is straight forward to check that the graph are correct and since both graphs have the same generating function, we conclude that they are Wilf-equivalent. □

### 6.7. Case G: The pattern set $\{231, 1234\}$

**Proposition 6.10.** *The pattern set $\{231, 1234\}$ is extendably symmetric with generating function*

$$\frac{1}{(1-x)^3(1-y)(1-z)} \left( x^2 \left( \frac{1}{(1-y)^2} + \frac{1}{(1-y)(1-z)} + \frac{1}{(1-z)^2} \right) \right.$$
$$\left. + \frac{1 - 5x + 9x^2 - 8x^3 + 5x^4}{(1-x)^2} + \frac{x(1 - 4x + 5x^2)}{1-x} \left( \frac{1}{1-y} + \frac{1}{1-z} - 1 \right) \right).$$

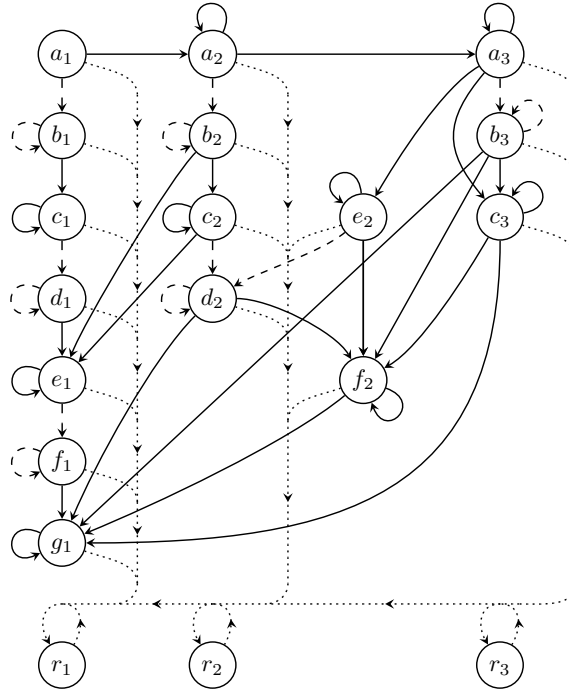

FIGURE 22. Generating graph for $\{231, 1234\}$.

*Proof.* Even though the generating graph in this case is finite, there are 19 different classes which makes the proof rather long. We show that the graph in Figure 22 is a generating graph for $\{231, 1234\}$. The classes are:

- $a_1$: The size zero permutation.
- $a_2$: Decreasing permutations.
- $a_3$: Permutations of the form $k(k-1)\cdots 1\, n(n-1)\cdots(k+1)$, where $1 < k < n$.
- $b_i = a_i^*$.
- $c_1 = {}^*a_2$.
- $c_2$: $a_3$, with at least one empty column added between the two decreasing sequences.
- $c_3$: Partial permutations of the form $k(k-1)\cdots 1\, l(l-1)\cdots(k+1)\, n(n-1)\cdots(m+1)$, where $1 < k < l \leqslant m < n$.
- $d_1 = c_1^*$.
- $d_2 = c_2^* \cup e_2^*$.
- $e_2$: Permutations of the form $n(n-1)\cdots(l+1)\, k(k-1)\cdots 1\, l(l-1)\cdots(k+1)$, where $1 < k < l < n$.
- $e_1$: Partial permutations of one of the forms $n(n-1)\cdots(m+1)\, l(l-1)\cdots(k+1)$, where $0 \leqslant k < l < m < n$, or $n(n-1)\cdots(m+1)\, k(k-1)\cdots 1\, m(m-1)\cdots(l+1)$, where $1 < k < l < m < n$.

- $f_2$: Partial permutations of one of the forms $j(j-1)\cdots 1\,n(n-1)\cdots(m+1)\,l(l-1)\cdots(k+1)$, where $1 < j \leqslant k < l < m < n$, $l(l-1)\cdots(k+1)\,j(j-1)\cdots 1\,k(k-1)\cdots(j+1)\,n(n-1)\cdots(l+1)$, where $1 < j < k < l < n$ or $j(j-1)\cdots 1\,n(n-1)\cdots(m+1)\,k(k-1)\cdots(j+1)\,m(m-1)\cdots(l+1)$, where $1 < j < k \leqslant l < m < n$.
- $f_1 = e_1^*$.
- $g_1$: Partial permutations which are not c-active and where only the top site is active.
- $r_i$: Partial permutations with at least one empty row and $i$ sites r-active.

It is an easy (but rather tedious) exercise to check that these classes do in fact give the graph above. For example the class $c_3$ has three active sites: the top site, which gives a partial permutation in class $g_1$, and the two sites between the descending sequences, which give partial permutations in classes $f_2$ and $c_3$ respectively.

Since the graph is finite, the generating function can be calculated with the transfer matrix method. $\qquad\square$

## Appendix A. Generating graphs and generating functions for combinations of patterns of length three

In the appendix we will give the generating graph and generating function for all possible combinations of patterns of length three. The pattern sets which contain both the patterns 123 and 321 are excluded since there are only finitely many pre-permutations avoiding these. Further, if there is a trivial symmetri between two sets of patterns, we will only give one of them, so for example {123, 231} is excluded, since it is the inverse of {123, 312}.

### A.1. Two patterns of size three

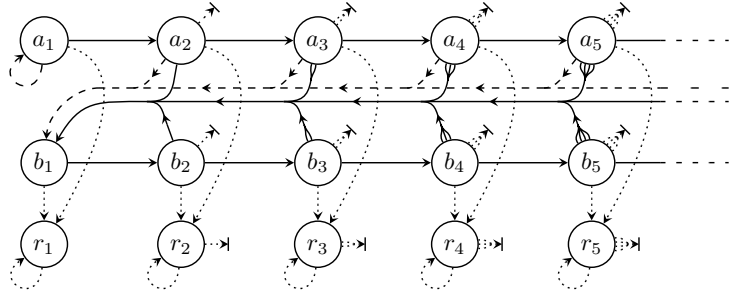

$$\text{Ext}_{123,132} = \frac{1}{(1-x)(1-2x)} \frac{xy+1-x}{1-y} \frac{xz+1-x}{1-z}$$

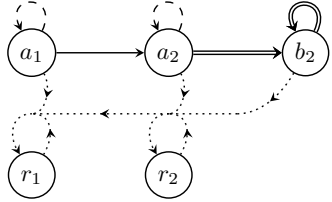

$$\text{Ext}_{123,213} = \frac{1}{1-y} \frac{1}{1-z} + \frac{x}{1-2x} \frac{1}{(1-y)^2} \frac{1}{(1-z)^2}$$

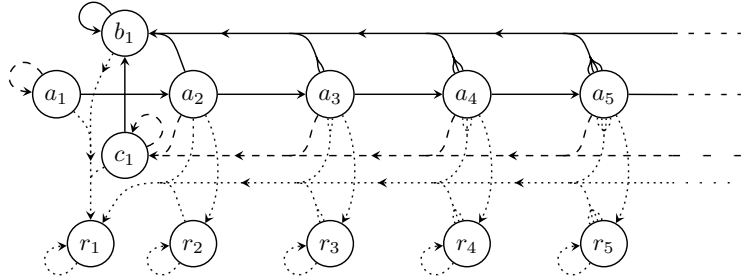

$$\text{Ext}_{123,312} = \frac{1}{1-y} \frac{1}{1-z} \left( \frac{1}{1-x} + \frac{x}{(1-x)^2} \left( \frac{z}{1-z} + \frac{y}{1-y} \right) + \frac{x^2}{(1-x)^3} \right)$$

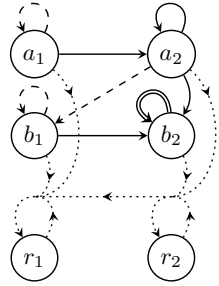

$$\text{Ext}_{132,213} = \frac{1}{1-x}\frac{1}{1-y}\frac{1}{1-z}\left(1-2x+\frac{x}{1-y}+\frac{x}{1-z}+\frac{x^2}{1-2x}\frac{1}{1-y}\frac{1}{1-z}\right)$$

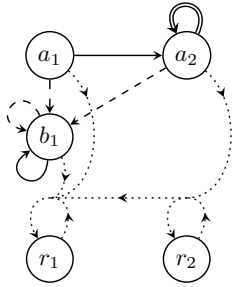

$$\text{Ext}_{132,231} = \frac{1}{1-z}\left(1+\frac{(1-x)y}{(1-2x)(1-x-y)}+\frac{x}{1-2x}\frac{1}{1-z}\right)$$

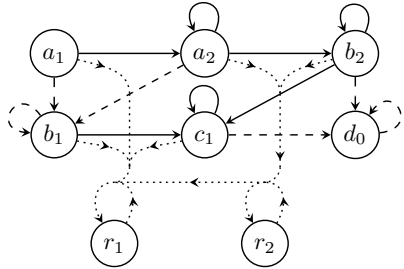

$$\text{Ext}_{132,321} = \frac{1}{1-x}\frac{1}{1-y}\frac{1}{1-z}+\frac{x^2}{(1-x)^3}\left(\frac{1}{1-y}+\frac{1}{1-z}-1\right)$$
$$+\frac{x}{(1-x)^2}\left(\frac{y}{1-y}\frac{z}{1-z}+\frac{y}{(1-y)^2}+\frac{z}{(1-z)^2}\right)$$

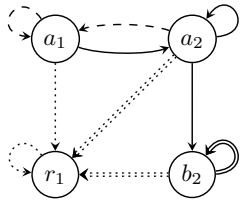

$$\text{Ext}_{213,312} = \frac{1-x}{1-2x}\frac{1}{1-x-y}\left(\frac{1}{1-z}-x\right)$$

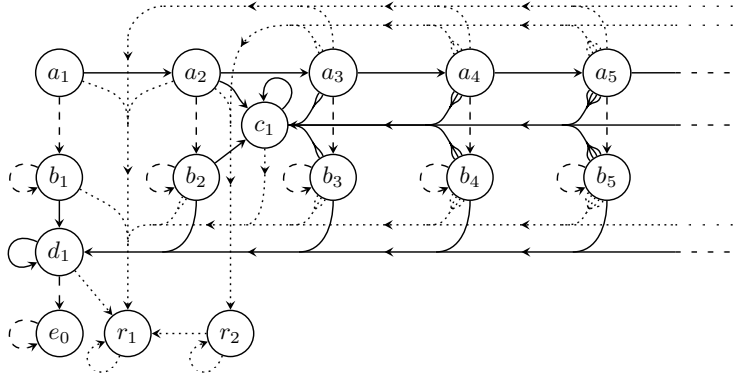

$$\text{Ext}_{213,321} = \frac{1}{1-x}\left(\frac{1}{1-y} + \frac{1}{1-z} - 1\right) + \frac{1+x}{(1-x)^2}\frac{y}{1-y}\frac{z}{1-z}$$
$$+ \frac{x}{(1-x)^2}\left(\frac{y}{(1-y)^2} + \frac{z}{(1-z)^2}\right) + \frac{x^2}{(1-x)^3}\frac{1}{1-y}\frac{1}{1-z}$$

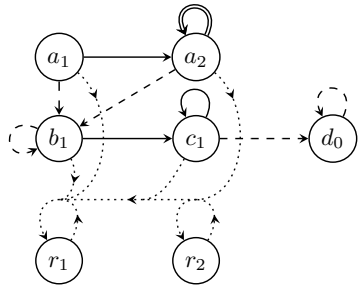

$$\text{Ext}_{231,312} = \frac{1}{1-2x}\left(x + \frac{1}{1-y}\frac{1}{1-z} - \frac{x(1-2y)}{(1-y)^2} - \frac{x(1-2z)}{(1-z)^2}\right)$$

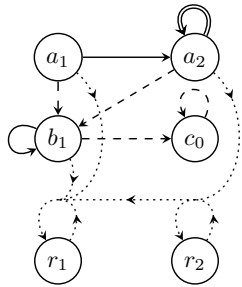

$$\text{Ext}_{312,321} = \frac{1}{1-z}\left(1 + \frac{y}{1-2x}\right) + \frac{1}{1-2x}\left(\frac{x}{(1-z)^2} + \frac{y^2}{1-y}\right)$$

## A.2. Three patterns of size three

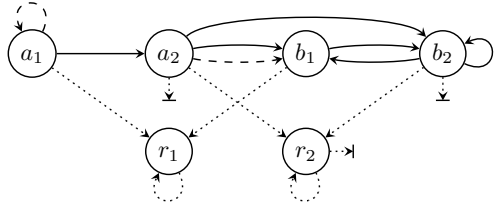

$$\mathrm{Ext}_{123,132,213} = \frac{1}{1-x-x^2}\frac{1+xy}{1-y}\frac{1+xz}{1-z}$$

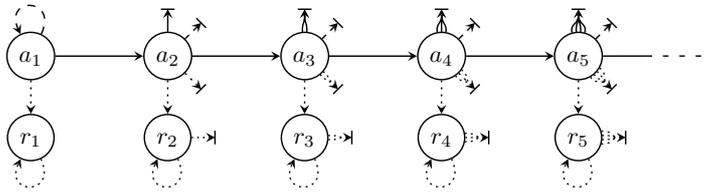

$$\mathrm{Ext}_{123,132,312} = \frac{1}{(1-x)^2}\frac{1}{1-y}\frac{1}{1-z} - \frac{x}{1-x}$$

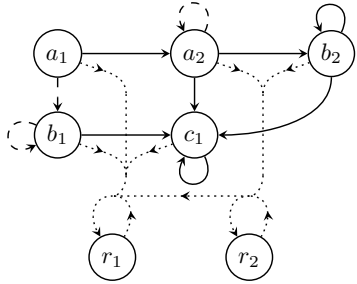

$$\mathrm{Ext}_{123,213,231} = \frac{1}{1-x}\frac{1}{1-y}\frac{1}{1-z}\left(\frac{1}{1-x} - x(1-y) + \frac{xz}{1-z}\right)$$

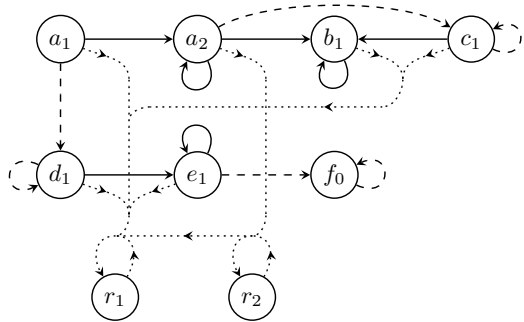

$$\mathrm{Ext}_{123,231,312} = \frac{1}{(1-x)^2}\frac{1}{1-y}\frac{1}{1-z} + \frac{x}{1-x}\left(\frac{y^2}{(1-y)^2} + \frac{z^2}{(1-z)^2} - 1\right)$$

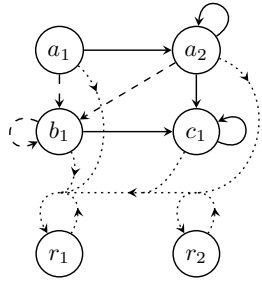

$$\text{Ext}_{132,213,231} = \frac{1}{(1-x)^2}\frac{1}{1-y}\frac{1}{1-z} + \frac{x}{1-x}\left(\frac{z}{(1-z)^2} - \frac{1}{1-z}\right)$$

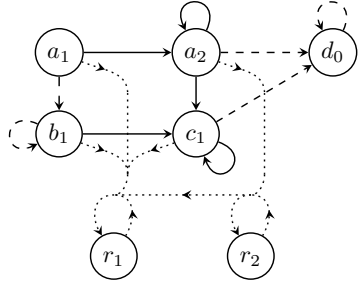

$$\text{Ext}_{132,231,312} = \frac{1}{1-x}\frac{1}{1-y}\frac{1}{1-z} + \frac{x^2}{(1-x)^2}\left(\frac{1}{1-y} + \frac{1}{1-z} - 1\right)$$
$$+ \frac{x}{1-x}\left(\frac{y}{(1-y)^2} + \frac{z}{(1-z)^2}\right)$$

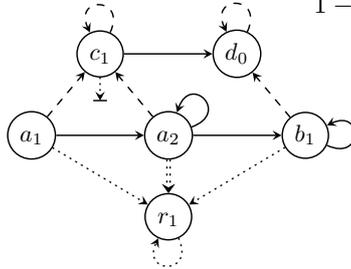

$$\text{Ext}_{132,231,321} = \frac{1}{1-x}\left(\frac{1}{1-x}\left(\frac{1}{1-y} + \frac{1}{1-z} - 1\right) + \frac{yz}{1-y} + \frac{xy^2}{(1-y)^2} - 1\right)$$

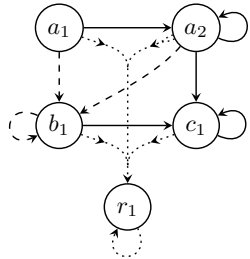

$$\text{Ext}_{213,231,312} = \frac{1}{(1-x)^2}\frac{1}{1-y}\frac{1}{1-z} - \frac{x}{1-x}$$

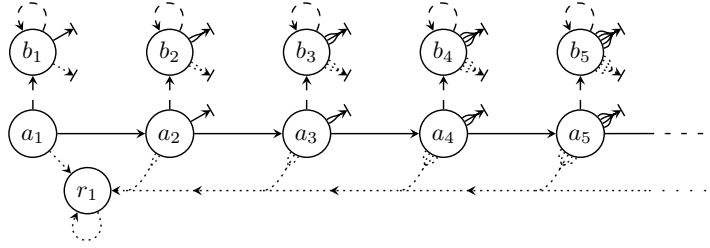

$$\text{Ext}_{213,231,321} = \frac{1+z}{(1-x)^2}\frac{1}{1-y} + \frac{1}{(1-x)^2}\frac{z^2}{1-z} - \frac{x}{1-x}$$

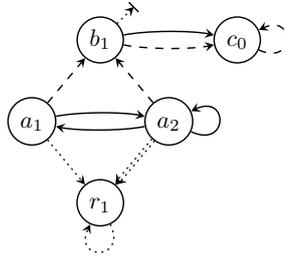

$$\text{Ext}_{231,312,321} = \frac{1}{1-x-x^2}\left((1+x)(\frac{1}{1-y} + \frac{1}{1-x} - 1) + yz - x\right)$$

## A.3. Four patterns of size three

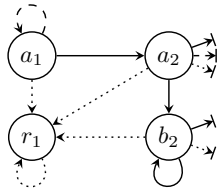

$$\text{Ext}_{123,132,213,312} = \frac{1}{1-x}\frac{1}{1-y}\frac{1}{1-z} + \frac{1}{1-y}\frac{x}{1-x}(1+z) - x$$

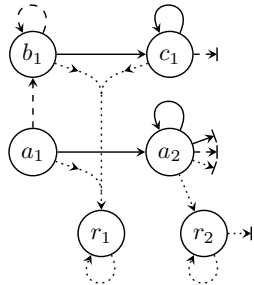

$$\text{Ext}_{123,132,231,312} = \frac{1}{1-x}\frac{1}{1-y}\frac{1}{1-z} + \frac{x}{1-x}\left(\frac{1}{1-y} + \frac{1}{1-z} - 1\right) - x$$

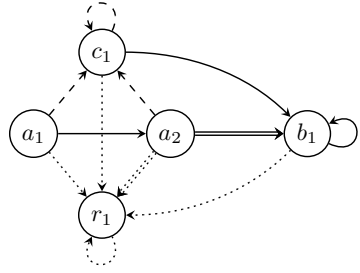

$$\text{Ext}_{123,213,231,312} = \frac{1+x}{1-x}\frac{1}{1-y}\frac{1}{1-z} - x$$

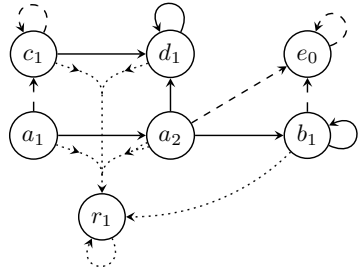

$$\text{Ext}_{132,213,231,312} = \frac{1}{1-x}\frac{1}{1-y}\frac{1}{1-z} + \frac{x}{1-x}\left(\frac{1}{1-y} + \frac{1}{1-z} - 1\right) - x$$

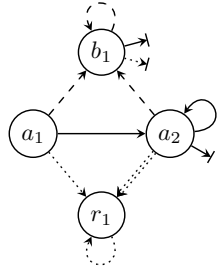

$$\text{Ext}_{132,213,231,321} = \frac{1+x}{1-x}\left(\frac{1}{1-y} + \frac{1}{1-z} - 1\right) + \frac{z}{1-x}\frac{y}{1-y} - x$$

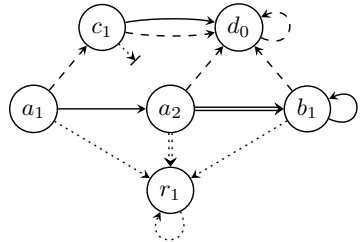

$$\text{Ext}_{132,231,312,321} = \frac{1+x}{1-x}\left(\frac{1}{1-y} + \frac{1}{1-z} - 1\right) + yz - x$$

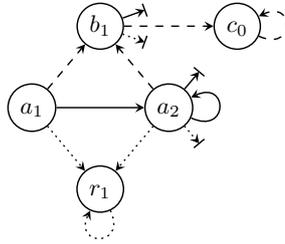

$$\text{Ext}_{213,231,312,321} = \frac{1}{1-x}\left(-1 + \frac{1}{1-y} + \frac{1}{1-z} + (x+y)(x+z)\right)$$

## A.4. Five patterns of size three

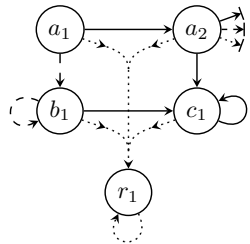

$$\text{Ext}_{123,132,213,231,312} = \frac{1}{1-x}\frac{1}{1-y}\frac{1}{1-z} + x(x+y+z)$$

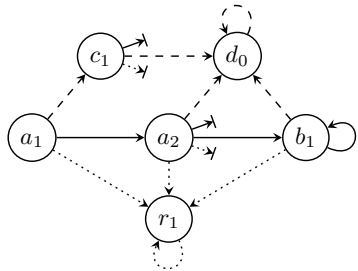

$$\text{Ext}_{132,213,231,312,321} = \frac{1}{1-x}\left(\frac{1}{1-y} + \frac{1}{1-z} - 1\right) + (x+y)(x+z)$$

Linköpings Universitet, 581 83 Linköping, Sweden
*E-mail address*: erouc@mai.liu.se

# Paper 3

**82**

# COMMUTATION RELATIONS FOR YOUNG TABLEAU INVOLUTIONS

ERIK OUCHTERLONY

ABSTRACT. In this paper we study a number of different involutions on Young tableaux. The relations between these involutions are classified and explained. The most remarkable result is that for Littlewood-Richardson tableaux, $a^3 = \mathrm{I}$, where $a$ is the composition of three different involutions: the fundamental symmetry map, the reversal involution and the rotation involution.

## 1. Preliminaries

### 1.1. Introduction

Over the years many bijections and involutions involving Young tableaux of different types have been found and investigated. In the paper by Pak and Vallejo [10], the authors studied a number of different bijection between Young tableaux in order to prove results about linear reductions, i.e., whether one bijection can be computed linearly in the time it takes to compute the other. The fact that almost all of the examined bijections were linearly equivalent (i.e. they are all pairwise linearly reducible to each other), indicate that the bijections have even tighter relationships than known before.

In this paper we study these relationships between the bijections in detail to get a better understanding of why the linear equivalences appear.

In Section 2 the companion and rotation involutions are defined and shown to commute. The section thereafter contains introductions to the well known tools, jeu de taquin, RSK and the tableau switching, which are used in Sections 4 and 5 to define the reversal and the fundamental symmetry map and show some commutation properties. From this the main result, Theorem 5.16, follows, which summarises all the commutation relations in the commutative diagram in Figure 21.

I thank Christian Krattenthaler for introducing me to the problem and for many interesting discussions. A thanks also goes to Bruce Sagan, Svante Linusson and Johan Wästlund for numerous comments and suggestions.

*Date*: December 17, 2005.

*Key words and phrases*. Young tableau, involutions, jeu de taquin, Robinson-Schensted correspondence, fundamental symmetry map.

The author was supported by the European Commission's IHRP Programme, grant HPRN-CT-2001-00272, "Algebraic Combinatorics in Europe".

### 1.2. **Basic Notation**

Let $\lambda \vdash n$ be a partition of $n$. Since there is a very natural bijection between partitions and Young diagrams, we use the same notation for these. The set difference between two Young diagrams, $\mu$ and $\lambda$, such that $\mu \subseteq \lambda$ is called a *skew diagram* and denoted $\lambda/\mu$. Filling the (skew) Young diagram, $\lambda/\mu$, with positive integers such that each row is weakly increasing and each column is strictly increasing gives us a *Young tableau* of *shape* $\lambda/\mu$. If $\mu$ is empty, we say that the Young tableau has *partition shape*. The set of Young tableaux of shape $\lambda/\mu$ is denoted $\mathrm{YT}(\lambda/\mu)$ and we say that $T \in \mathrm{YT}(\lambda/\mu)$ has *weight* $\boldsymbol{c} = (c_1, c_2, \ldots)$ if $c_i$ is the number of occurrences of the letter $i$ in $T$ and we write $\mathrm{YT}(\lambda/\mu, \boldsymbol{c})$ for all Young tableaux of shape $\lambda/\mu$ which have weight $\boldsymbol{c}$. If the weight is weakly decreasing, we may also use a partition instead of a vector to denote the weight. To indicate when we regard the partitions as vector, we will use bold face.

### 1.3. **Littlewood-Richardson tableaux**

From a tableau $T \in \mathrm{YT}(\lambda/\mu)$ we can construct a word, $w(T)$, by reading its entries row by row, starting from the top and reading each row from right to left. In the same way as for the tableaux, the weight of a word is a vector with the number of entries of each value, so we have $\mathrm{weight}(w(T)) = \mathrm{weight}(T)$.

We say that $w$ is a Yamanouchi word if for every prefix $v$ of $w$, $\mathrm{weight}(v)$ is a partition, i.e., weakly decreasing.

**Definition 1.1.** $T \in \mathrm{YT}(\lambda/\mu)$ is a Littlewood-Richardson tableau if $w(T)$ is a Yamanouchi word.

$$T = \begin{array}{ccc} & & \boxed{1} \\ & \boxed{1} & \boxed{2} \\ \boxed{1} & \boxed{2} & \boxed{3} \\ \boxed{2} & \boxed{4} \end{array} \qquad w(T) = 12132142$$

FIGURE 1. Example of a Littlewood-Richardson tableau.

### 1.4. **Schur functions**

One of the main motivations for Littlewood-Richardson and Young tableaux comes from the (skew) Schur functions.

$$s_{\lambda/\mu} = \sum_{T \in \mathrm{YT}(\lambda/\mu)} x^{\mathrm{weight}(T)}$$

The Schur functions corresponding to partition shapes constitute a basis for the symmetric functions, so we can express their product with respect to this basis:

$$s_\mu s_\nu = \sum_\lambda c^\lambda_{\mu\nu} s_\lambda \qquad (1)$$

The same coefficients appear when expanding the skew Schur functions:

$$s_{\lambda/\mu} = \sum_\nu c^\lambda_{\mu\nu} s_\nu \qquad (2)$$

The numbers $c^\lambda_{\mu\nu}$ are known as the Littlewood-Richardson coefficients due to the following proposition, called the Littlewood-Richardson rule, which gives a combinatorial interpretation in terms of Young tableaux, and was first stated, but not proved, by Littlewood and Richardson [9]. The first proof is due to Schützenberger [12], and today there are plenty of different proofs, see for example [14], [11], [8] and [17]. For more remarks about the history of the Littlewood-Richardson rule, see [16].

**Proposition 1.2.** *The Littlewood-Richardson coefficients $c^\lambda_{\mu\nu}$ are equal to the number of Littlewood-Richardson tableaux in* $\mathrm{YT}(\lambda/\mu, \boldsymbol{\nu})$.

From equation (1), it is evident that $c^\lambda_{\mu\nu} = c^\lambda_{\nu\mu}$. This symmetry is, however, not obvious in the combinatorial representation by Young tableaux. This implies that there should exists a bijection between Littlewood-Richardson tableaux in $\mathrm{YT}(\lambda/\mu, \boldsymbol{\nu})$ and those in $\mathrm{YT}(\lambda/\nu, \boldsymbol{\mu})$. Such a bijection is called a *fundamental symmetry map*.

The fact the equations (1) and (2) give the same coefficients implies that there are even more symmetries to examine.

## 2. **The companion involution and the rotation involution**

### 2.1. **The companion tableaux**

Let $T \in \mathrm{YT}$. The *recording matrix* of $T$, $M = \mathrm{rec}(T)$, is the matrix $M = (m_{ij})$, where $(m_{ij})$ is the number of $j$:s in row $i$. It is clear that the shape together with the recording matrix uniquely determine a tableau.

**Lemma 2.1.** *Let $T_1, T_2 \in \mathrm{YT}(\lambda/\mu)$. Then*

$$\mathrm{rec}(T_1) = \mathrm{rec}(T_2) \quad \Longleftrightarrow \quad T_1 = T_2$$

Following the notation of [16], we say that a *companion tableau* to $T$ is a tableau with $M^t$ (the transpose of $M$) as its recording matrix.

$$T_1 = \begin{array}{c}\end{array} \qquad \mathrm{rec}(T_1) = \begin{pmatrix} 1 & 1 & 0 & 0 \\ 1 & 1 & 1 & 0 \\ 0 & 0 & 2 & 0 \\ 0 & 1 & 0 & 2 \end{pmatrix} \qquad T_2 = \begin{array}{c}\end{array}$$

FIGURE 2. The Young tableaux $T_1$ and $T_2$ are companions, since $\mathrm{rec}(T_1) = \mathrm{rec}(T_2)^t$.

## 2.2. $\kappa$-dominance

Let $T \in \mathrm{YT}(\lambda/\mu)$. If $T$ has a companion in $\mathrm{YT}(\nu/\kappa)$, we say that $T$ is $\kappa$-dominant (or $\nu/\kappa$-dominant) and use the notation

$$\mathrm{YT}(\lambda/\mu, \nu/\kappa) = \{T \in \mathrm{YT}(\lambda/\mu, \boldsymbol{\nu} - \boldsymbol{\kappa}) : T \text{ is } \kappa\text{-dominant}\}.$$

A companion tableau is by no means unique since we only specify the recording matrix, not the shape. However, in $\mathrm{YT}(\lambda/\mu, \nu/\kappa)$, we can require the companion to have a specific shape, namely $\nu/\kappa$. This gives us an involution, by simply switching the tableau and its companion.

**Definition 2.2.** The *companion involution*, $\tau : \mathrm{YT}(\lambda/\mu, \nu/\kappa) \rightarrow \mathrm{YT}(\nu/\kappa, \lambda/\mu)$, maps a tableau to its unique companion of shape $\nu/\kappa$.

Note that $\tau(T)$ is in fact $\mu$-dominant, since by definition, $T$ is a companion to $\tau(T)$. We can also see that a tableau in $\mathrm{YT}(\lambda/\mu, \nu/\kappa)$ is uniquely determined by its recording matrix, since both the shape are fixed.

One reason this new notation works so well is the following characterisation of the Littlewood-Richardson tableaux, which has been used as the definition in e.g. [16].

**Proposition 2.3.** *$T$ is a Littlewood-Richardson tableaux* $\iff$ *$T$ is 0-dominant.*

A 0-dominant tableau is a tableau which has a partition shaped companion, so the companion involution gives us a bijection between the Littlewood-Richardson tableaux in $\mathrm{YT}(\lambda/\mu, \nu)$ and the partition shaped tableaux in $\mathrm{YT}(\nu, \lambda/\mu)$. In the current literature, the set of Littlewood-Richardson tableaux have mostly been denoted by LR. However, since they are so closely related, we would like to introduce a notation which reflects this connection.

**Definition 2.4.**

$$\mathrm{YT}^{\ulcorner} = \{T \in \mathrm{YT} : T \text{ has partition shape}\}$$

$$\mathrm{YT}_{\ulcorner} = \{T \in \mathrm{YT} : T \text{ is Littlewood-Richardson}\}$$

The superscript $\ulcorner$ is intended to be an indication that the tableau is partition shaped, and similarly, if it is a subscript, it indicates that the companion is partition shape, i.e., it is a Littlewood-Richardson tableau.

## 2.3. The rotation involution

Let $\square$ be a *rectangular shaped* Young diagram with $l$ rows and $k$ columns and let $\lambda/\mu \subseteq \square$, with $\lambda = (\lambda_1, \ldots, \lambda_l)$ and $\mu = (\mu_1, \ldots, \mu_l)$, where some of the trailing values may be zero. Then we define

$$(\lambda/\mu)^{\bullet} = \mu^c/\lambda^c,$$

where $\lambda^c = (k - \lambda_l, k - \lambda_{l-1}, \ldots, k - \lambda_1)$ and $\mu^c = (k - \mu_l, k - \mu_{l-1}, \ldots, k - \mu_1)$ are the complementary shapes with respect to $\square$, so that the skew diagram is rotated $180°$ with respect to $\square$. The shape $\square$ is called the *bounding box* for $\lambda/\mu$.

Let $T \in \mathrm{YT}(\lambda/\mu, \nu/\kappa)$ have the recording matrix $\mathrm{rec}(T) = M = (m_{i,j}) \in \mathbb{N}^{\ell(\lambda) \times \ell(\nu)}$, where $\ell(\lambda)$ is *length* of the partition $\lambda$, i.e., the number of parts. We define $T^\bullet \in \mathrm{YT}((\lambda/\mu)^\bullet, (\nu/\kappa)^\bullet)$ to be the unique tableau (given bounding boxes for $\lambda/\mu$ and $\nu/\kappa$) with recording matrix $M^\bullet = (m_{\ell(\lambda)+1-i, \ell(\mu)+1-j})$, i.e. a 180° rotation of $M$.

$$T = \begin{array}{c}\boxed{\begin{smallmatrix}1&2\end{smallmatrix}}\\ \boxed{\begin{smallmatrix}1&2&3\end{smallmatrix}}\\ \boxed{\begin{smallmatrix}3&3\end{smallmatrix}}\\ \boxed{\begin{smallmatrix}2&4&4\end{smallmatrix}}\end{array} \qquad \mathrm{rec}(T)^\bullet = \left(\begin{array}{cccc} 2 & 0 & 1 & 0 \\ 0 & 2 & 0 & 0 \\ 0 & 1 & 1 & 1 \\ 0 & 0 & 1 & 1 \end{array}\right) \qquad T^\bullet = \begin{array}{c}\boxed{\begin{smallmatrix}1&1&3\end{smallmatrix}}\\ \boxed{\begin{smallmatrix}2&2\end{smallmatrix}}\\ \boxed{\begin{smallmatrix}2&3&4\end{smallmatrix}}\\ \boxed{\begin{smallmatrix}3&4\end{smallmatrix}}\end{array}$$

FIGURE 3. Example of the rotation.

Note that the rotation is not involutive in general as shown in Figure 4.

$$T = \begin{array}{c}\vdots\\ \boxed{2}\\ \boxed{2\ 3}\end{array} \qquad T^\bullet = \begin{array}{c}\boxed{1\ 2}\\ \boxed{2}\\ \vdots\end{array} \qquad S = (T^\bullet)^\bullet = \begin{array}{c}\boxed{1}\\ \boxed{1\ 2}\end{array} \neq T$$

FIGURE 4. Non-involutive rotation.

In order to correct this problem, there are two possible solutions. One could regard $S$ and $T$ as equal by introducing an equivalence relation. However, this would give us problems later with the fundamental symmetry map, which we define in section 5, since $\rho(S) \neq \rho(T)$.

Instead, we keep track of which rectangular shape the tableau was rotated with respect to, so we associate to each tableau two bounding boxes, one for the tableau shape and one for the companion shape.

**Definition 2.5.** Let $(\lambda/\mu)_\square$ denote the shape $\lambda/\mu \subseteq \square$, bounded by the rectangular shaped Young diagram $\square$.

**Lemma 2.6.** *The rotation map $T \mapsto T^\bullet$ is an involution* $\mathrm{YT}((\lambda/\mu)_{\square_1}, (\nu/\kappa)_{\square_2}) \to \mathrm{YT}((\lambda/\mu)^\bullet_{\square_1}, (\nu/\kappa)^\bullet_{\square_2})$

*Proof.* Let $\square_1 = (i^l)$ and $\square_2 = (j^k)$. Then $\mathrm{rec}(T), \mathrm{rec}(T^\bullet) \in \mathbb{N}^{l \times k}$ and $\mathrm{rec}(T^\bullet) = \mathrm{rec}(T)^\bullet$, so $\mathrm{rec}((T^\bullet)^\bullet) = \mathrm{rec}(T)$. Since $(T^\bullet)^\bullet$ has shape $((\lambda/\mu)^\bullet_{\square_1})^\bullet_{\square_1} = \lambda/\mu$, we have $(T^\bullet)^\bullet = T$. $\square$

### 2.4. The companion and rotation involutions commute

**Lemma 2.7.** *Let $T \in \mathrm{YT}((\lambda/\mu)_{\square_1}, (\nu/\kappa)_{\square_2})$, where $\square_1 = (i^l)$ and $\square_2 = (j^k)$. Then*

$$\tau(T^\bullet) = \tau(T)^\bullet$$

*Proof.* Let $M = \mathrm{rec}(T)$, $M \in \mathbb{N}^{l \times k}$. Then both sides have the same shape, $(\nu/\kappa)^\bullet_{\square_2}$, companion shape, $(\lambda/\mu)^\bullet_{\square_1}$, and recording matrix, $(M^t)^\bullet = (M^\bullet)^t$. $\square$

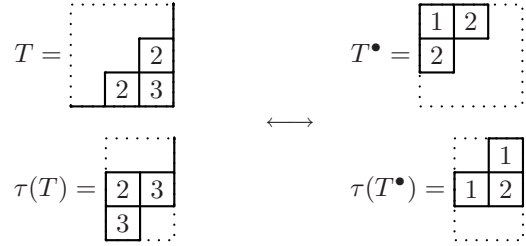

FIGURE 5. Involutive rotation.

Whenever there is a choice to make which bounding boxes to select for a given tableau, we can, for example, take the smallest possible. After the choice is made, the bounding boxes are fixed and involutions, like the rotation, will be made with respect to them. In the rest of the paper we will implicitly assume that any shape has an associated bounding box, so that the rotation is involutive.

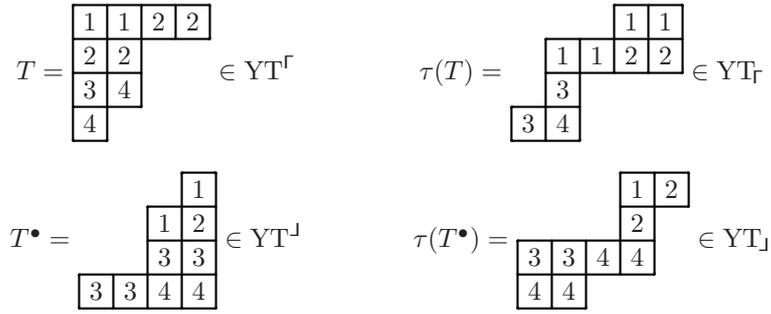

FIGURE 6. Examples of tableaux of (anti-)partition shape and of (anti-)Littlewood-Richardson tableaux.

2.5. **Anti-Littlewood-Richardson tableaux**

Using bounding boxes, it is natural to define the rotated versions of partition shape tableaux and Littlewood-Richardson tableaux.

**Definition 2.8.** The sets of *anti-Littlewood-Richardson tableaux*, $\mathrm{YT}_\lrcorner$, and *anti-partition* shaped tableaux, $\mathrm{YT}^\lrcorner$, are defined as

$$\mathrm{YT}_\lrcorner = \{T \in \mathrm{YT} : T^\bullet \in \mathrm{YT}_\ulcorner\},$$
$$\mathrm{YT}^\lrcorner = \{T \in \mathrm{YT} : T^\bullet \in \mathrm{YT}^\ulcorner\}.$$

See Figure 6 for examples.

2.6. **Canonical tableaux**

Let $\mathrm{YT}_\ulcorner^\ulcorner(\lambda) = \mathrm{YT}^\ulcorner(\lambda) \cap \mathrm{YT}_\ulcorner(\lambda)$. Then the recording matrix has to be diagonal, so that $\lambda$ only has ones in the first row, only twos in the second and so on. Therefore $\mathrm{YT}_\Gamma^\Gamma(\lambda)$ contains exactly one element which is called the *canonical tableau* of shape $\lambda$ and is denoted $\mathrm{can}(\lambda)$. In the same way, $\mathrm{can}(\lambda^\bullet), \mathrm{can}(\lambda)^\bullet$ and $\mathrm{can}(\lambda^\bullet)^\bullet$ are the unique elements in $\mathrm{YT}_\Gamma^\Gamma(\lambda^\bullet)$, $\mathrm{YT}_\lrcorner^\lrcorner(\lambda^\bullet)$ and $\mathrm{YT}_\lrcorner^\Gamma(\lambda)$, respectively. The last two are called *anti-canonical* tableaux, since they are both anti-Littlewood-Richardson tableaux. See Figure 2.6 for examples. Note that if $\lambda$ has rectangular shape, then $\mathrm{can}(\lambda) = \mathrm{can}(\lambda^\bullet) = \mathrm{can}(\lambda)^\bullet = \mathrm{can}(\lambda^\bullet)^\bullet$.

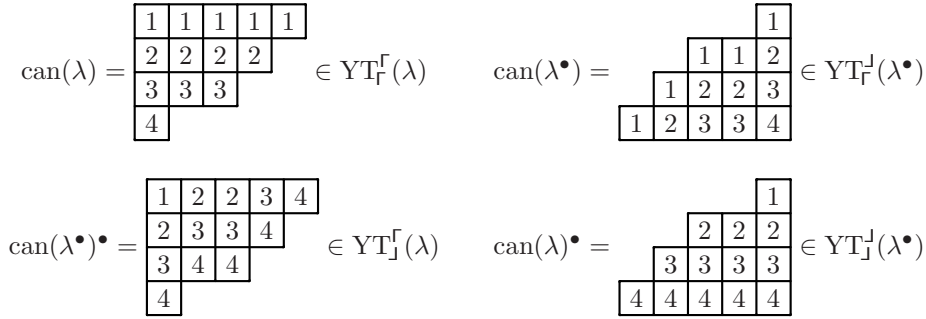

Figure 7. Examples of (anti-)canonical tableaux of different types for $\lambda = (5, 4, 3, 1)$.

## 3. Jeu de taquin and tableau switching

### 3.1. Jeu de taquin

A fundamental operation on Young tableaux is Schützenberger's [12] jeu de taquin. The name is French from the famous 15-puzzle, which well illustrates what is going on. To do a jeu de taquin slide on a tableaux, $T \in \mathrm{YT}(\lambda/\mu)$, we choose a box which is a south east corner of $\mu$ as the *blank* of the 15-puzzle. The blank is then used to do slides, such that the blank moves downwards or to the right until it ends up on the southeast side of the tableau. Each slide is uniquely determined by the tableau conditions since $\boxed{a\ b}$ is allowed iff $\boxed{\begin{smallmatrix} b \\ a \end{smallmatrix}}$ is not allowed, which means that the smallest element available is always moved.

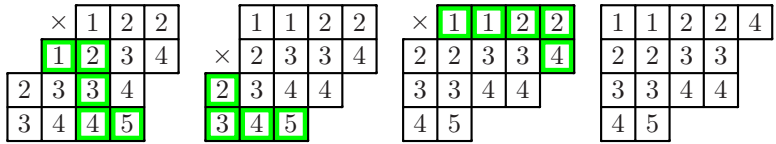

Figure 8. Examples of jeu de taquin slides, starting with $T$ and ending with $\mathrm{jdt}(T)$.

By repeatedly performing the slides on a tableau $T$, we finally end up with a partition shaped tableau, called the *jeu de taquin map* of $T$ and denoted by $\mathrm{jdt}(T)$.

It is an important property of jeu de taquin that jdt($T$) does not depend on the order we choose the corners. This is often referred to as the confluence of jeu de taquin.

**Proposition 3.1.** *The* jdt *map is well defined, .i.e., it does not depend on the order in which the jeu de taquin slides are made.*

For a proof, see, for example [16].

### 3.2. **Bender-Knuth transformations**

A close relative to the jeu de taquin slide are the Bender-Knuth transformations first defined in [1].

**Definition 3.2.** The Bender-Knuth map, $s_i : \text{YT}(\lambda/\mu, \boldsymbol{c}) \to \text{YT}(\lambda/\mu, \sigma_i(\boldsymbol{c}))$, where $\sigma_i$ switches the element $i$ and $i+1$ of the vector $\boldsymbol{c}$, is defined as follows:

(1) All entries filled with a number different from $i$ and $i+1$ are left unchanged.
(2) All columns which contain both an $i$ and an $i+1$ are also unchanged.
(3) For each row, if there are $k$ $i$:s and $l$ $i+1$:s in columns not satisfying (2), change the numbers in these positions into $l$ $i$:s and $k$ $i+1$:s.

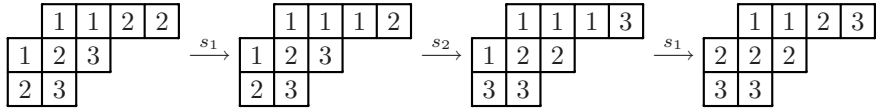

FIGURE 9. Example of Bender-Knuth transformations.

A sequence of Bender-Knuth transformations, $s_{m-1} s_{m-2} \dots s_1$, corresponds to doing jeu de taquin slides with several blank squares at the same time in the following way: We let all squares filled with the number 1 be the blanks. Then do the slides, starting with the right most blank. When finished, decrease all numbers with one and fill the blanks with the number $m$. The following combinations of these kind of sequences will later be used to define the tableau switching map and the Schützenberger involution:

$$z_m = (s_1 s_2 \cdots s_{m-1})(s_1 s_2 \cdots s_{m-2}) \cdots (s_1 s_2)(s_1)$$
$$t_{k,l} = (s_l s_{l+1} \cdots s_{l+k-1}) \cdots (s_2 s_3 \cdots s_{k+1})(s_1 s_2 \cdots s_k)$$
$$t_{k,l}^{(d)} = (s_{l+d} s_{l+d+1} \cdots s_{l+d+k-1}) \cdots (s_{d+2} s_{d+3} \cdots s_{d+k+1})(s_{d+1} s_{d+2} \cdots s_{d+k})$$

To see what these operation do, let $\boldsymbol{c} = (c_1, c_2, \dots, c_m)$ be a vector and $\boldsymbol{c}^r = (c_m, c_{m-1}, \dots, c_1)$ be the *reverse* vector of $\boldsymbol{c}$. If $T$ is a tableau whose weight vector is $\boldsymbol{c}$, then $z_m(T)$ is a tableau with the reverse weight $\boldsymbol{c}^r$, and for $t_{k,m-k}$, the upper and lower part of the weight vector become exchanged, so the new weight vector is $(c_{k+1}, c_{k+2}, \dots, c_m, c_1, c_2, \dots, c_k)$. The third operation is similar, but the $d$ first elements of the weight vector are left unchanged.

**Lemma 3.3.** *The Bender-Knuth transformations satisfy the following relations:*

(i) $s_i^2 = 1$,
(ii) $s_i s_j = s_j s_i$, *if* $|i - j| \geqslant 2$,
(iii) $z_m^2 = 1$,
(iv) $t_{k,l} \, t_{l,k} = 1$,
(v) $z_{l+k} = z_k \, t_{l,k} \, z_l$,
(vi) $t_{l+k,m} = t_{l,m} \, t_{k,m}^{(l)}$.

*Proof.* (*i*) and (*ii*) are obvious from the definition and (*iii*), (*iv*) and (*v*) are proved in [10]. For (*vi*), we can use (*ii*) to get

$$t_{l+k,m} = (s_m \cdots s_{m+l+k-1}) \cdots (s_2 \cdots s_{l+k+1})(s_1 \cdots s_{l+k})$$
$$= \Big((s_m \cdots s_{m+l+k-1}) \cdots (s_2 \cdots s_{l+k+1})(s_1 \cdots s_l)\Big)\Big((s_{l+1} \cdots s_{l+k})\Big)$$
$$= \Big((s_m \cdots s_{m+l+k-1}) \cdots (s_2 \cdots s_{l+1})(s_1 \cdots s_l)\Big)\Big((s_{l+2} \cdots s_{l+k+1})(s_{l+1} \cdots s_{l+k})\Big)$$
$$\vdots$$
$$= \Big((s_m \cdots s_{m+l-1}) \cdots (s_1 \cdots s_l)\Big)\Big((s_{m+l} \cdots s_{m+l+k}) \cdots (s_{l+1} \cdots s_{l+k})\Big)$$
$$= t_{l,m} \, t_{k,m}^{(l)}.$$

$\square$

### 3.3. **Tableau switching**

The *tableau switching map*, as defined by Benkart, Sottile and Stroomer [2], is an involution, $X : \mathrm{YT}(\nu/\mu, \boldsymbol{c}) \times \mathrm{YT}(\lambda/\nu, \boldsymbol{d}) \to \mathrm{YT}(\pi/\mu, \boldsymbol{d}) \times \mathrm{YT}(\lambda/\pi, \boldsymbol{c})$, defined by the procedure below, where the tableau switching is acting on the tableau pair $(S, T)$.

(1) Let $k = \ell(\boldsymbol{c})$ and $l = \ell(\boldsymbol{d})$.
(2) Increase the values in the tableau $T$ by $k$ to make sure they are all larger than all values in the tableau $S$.
(3) Take the union of $S$ and the modified $T$ to form a single tableau of shape $\lambda/\mu$.
(4) Make the switch by applying the Bender-Knuth transformations $t_{k,l}$.
(5) Split the tableau so that the tableau $T'$ contains all squares with values $1 \ldots l$ and $S'$ contains the rest.
(6) Reduce the values in $S'$ by $l$.

Now $X(S, T) := (T', S')$, where $T' \in \mathrm{YT}(\pi/\mu, \boldsymbol{d})$ and $S' \in \mathrm{YT}(\lambda/\pi, \boldsymbol{c})$.

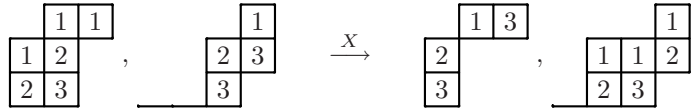

FIGURE 10. Example of the tableau switching map.

From the confluence of jeu de taquin it follows that tableau switching is involutive, so that $X(X(S, T)) = (S, T)$. It is also be possible to define the tableau switching map using outward jeu de taquin slides. Regard the squares in the inner tableau as blanks and slide them outwards starting with the square with highest value and, if there are several, choose the rightmost one first. One can see from the definition of the Bender-Knuth transformations that this give the same map. From this alternative description we get an interpretation of the first tableau of the image of the tableau switching map.

**Lemma 3.4.** *Let $T \in \mathrm{YT}(\lambda/\nu)$, $S \in \mathrm{YT}^\ulcorner(\nu)$ and $X(S, T) = (T', S')$. Then*

$$T' = \mathrm{jdt}(T)$$

Tableau switching can be illustrated with pictures like the one in Figure 3.3.

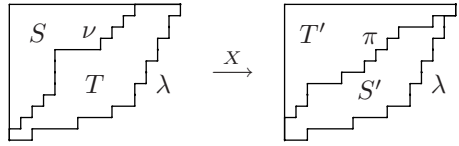

FIGURE 11. The tableau switching map applied to $S \in \mathrm{YT}(\nu)$ and $T \in \mathrm{YT}(\lambda/\nu)$.

## 4. **The Schützenberger involution and the reversal involution**

The Schützenberger involution was defined in [12] and is very important in the theory of the RSK correspondence.

**Definition 4.1.** Let $T \in \mathrm{YT}^\ulcorner(\lambda, \boldsymbol{c})$ be a partition shaped tableau with weight $\boldsymbol{c} = (c_1, \ldots, c_m)$, then the *Schützenberger involution*, $\xi : \mathrm{YT}^\ulcorner(\lambda, \boldsymbol{c}) \to \mathrm{YT}^\ulcorner(\lambda, \boldsymbol{c}^r)$, is the map

$$\xi(T) = z_m(T) = (s_1 s_2 \cdots s_m)(s_1 s_2 \cdots s_{m-1}) \cdots (s_1 s_2)(s_1)(T).$$

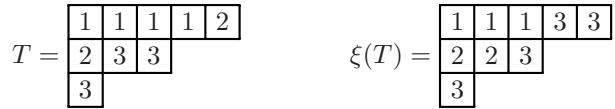

FIGURE 12. Example of the Schützenberger involution.

The Schützenberger involution can also be defined in terms of the jeu de taquin slides by the following procedure. Let the rightmost box with the smallest value, $v$, be the blank and slide it to the right side of the tableau with jeu de taquin slides. Then change its value to $m + 1 - v$ and mark it as being outside the tableau, so that it will not participate in the rest of the procedure. Next repeat this with a new box until all the boxes have been used as the blank precisely once. The resulting tableau is $\xi(T)$.

### 4.1. **The RSK correspondence**

The Robinson-Schensted correspondence is a bijection between permutations and a pair of standard tableaux. It was generalised by Knuth [6] to a bijection between nonnegative integer matrices and a pair of tableaux in $\mathrm{YT}^{\ulcorner}$:

$$\mathrm{RSK} : \mathrm{Mat} \longrightarrow \bigcup_\lambda \mathrm{YT}^{\ulcorner}(\lambda) \times \mathrm{YT}^{\ulcorner}(\lambda).$$

Let $M = (m_{ij})$ be a matrix with nonnegative integer entries. We define the RSK correspondence by an algorithm which transforms the matrix $M$ into the pair of partition shaped tableaux, $P$ and $Q$. At the start of the algorithm both $P$ and $Q$ are empty.

(1) Let $i$ be the first non-zero row of $M$ and let $j$ be the first column in row $i$ which has a positive value.
(2) Reduce $m_{ij}$ by one.
(3) Add the value $j$ to the tableau $P$ through the a bumping process which starts on the first row. If $j$ is greater or equal to all the entries of the row, add a box at the end, filled with the value. Otherwise, find the leftmost box that is has a value $k > j$ and put the value $j$ into that box. Next, repeat the bumping on the row below, but with $k$ instead of $j$.
(4) Add a new box to the tableau $Q$ at the same place as was just added to $P$ and fill the box with the value $i$.
(5) Repeat until $M$ is the zero matrix.

Much more theory about the RSK correspondence can be found, e.g., in the books by Fulton [4] and Stanley [13]. The RSK correspondence is very important when studying Young tableaux, and we will use it many times. The following classical result shows the relationship between RSK and the Schützenberger involution.

**Theorem 4.2** (Knuth)**.** *Assume* $\mathrm{RSK}(M) = (P, Q)$, *then*

$$\mathrm{RSK}(M^t) = (Q, P)$$
$$\mathrm{RSK}(M^\bullet) = (\xi(P), \xi(Q)).$$

The RSK correspondence, the jeu de taquin map and the Schützenberger involution are all very tightly related, as was shown already by Schützenberger [12]. Using the notation $M^r$ for the reverse of the matrix $M$, i.e., the matrix turned upside down, we get

**Lemma 4.3.** *Let* $T \in \mathrm{YT}(\lambda/\mu, \nu/\kappa)$ *with* $\mathrm{rec}(T) = M$ *and* $\mathrm{RSK}(M^r) = (P, Q)$. *Then*

(i) $\mathrm{jdt}(T) = P$
(ii) $T \in \mathrm{YT}^{\ulcorner} \implies T = P$
(iii) $\mathrm{jdt}(T^\bullet) = \xi(P)$
(iv) $T \in \mathrm{YT}^{\ulcorner} \implies \mathrm{jdt}(T^\bullet) = \xi(T)$
(v) $Q = \mathrm{jdt}(\tau(T^\bullet)) = \xi(\mathrm{jdt}(\tau(T)))$.

*Proof.* The first statement was proved by Schützenberger [12] and the second follows from the first, since $\mathrm{jdt}(T) = T$ if $T$ already has partition shape. From Theorem 4.2, we get $\mathrm{RSK}((M^\bullet)^r) = \mathrm{RSK}((M^r)^\bullet) = (\xi(P), \xi(Q))$, which in the light of (i), shows (iii), from which (iv) is an easy consequence. For the last part, we have, again by Theorem 4.2, that $\mathrm{RSK}(((M^t)^\bullet)^r) = \mathrm{RSK}((M^r)^t) = (Q, P)$. Hence $Q = \mathrm{jdt}(\tau(T^\bullet)) = \xi(\mathrm{jdt}(\tau(T)))$, where the last equality follows from (i) and (iii). $\square$

It may seem strange that we only define the Schützenberger involution for partition shaped tableaux, when the definition would in fact work fine for skew shapes as well. The problem is that the more general version does not behave very well, for example it does not preserve the structure of the companion, so that $T \in \mathrm{YT}(\lambda/\mu, \nu/\kappa)$ does not imply that $\xi(T) \in \mathrm{YT}(\lambda/\mu, (\nu/\kappa)^\bullet)$. Furthermore, the previous lemma does not generalise to skew shapes.

### 4.2. Reversal

To get a generalisation of the Schützenberger involution which behaves better, we use tableau switching to bring it into partition shape and back. This definition is due to Benkart, Sottile and Stroomer [2].

**Definition 4.4.** The *reversal* involution, $\chi^S : \mathrm{YT}(\lambda/\mu, \boldsymbol{c}) \to \mathrm{YT}(\lambda/\mu, \boldsymbol{c}^r)$, is a conjugation of the Schützenberger involution by tableau switching. In other words, if $X(S, T) = (T', S')$ and $X(\xi(T'), S') = (S'', T'')$, then $\chi^S(T) = T''$. See Figure 13 for an illustration.

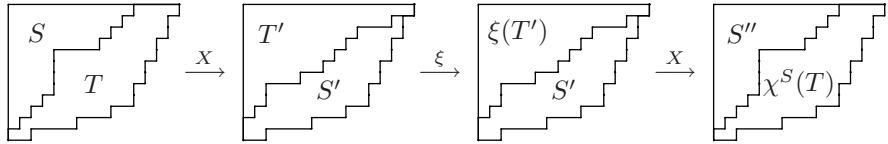

FIGURE 13. Definition of the reversal.

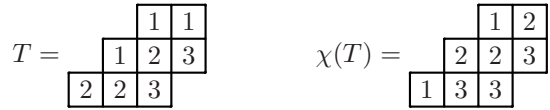

FIGURE 14. Example of the reversal.

**Lemma 4.5.** $\chi^S : \mathrm{YT}(\lambda/\mu, \nu/\kappa) \longrightarrow \mathrm{YT}(\lambda/\mu, (\nu/\kappa)^\bullet)$ *is an involution.*

It is clear that $\chi^S$ is an involution since both Schützenberger's involution and tableau switching are. So what we need to show is that if $T \in \mathrm{YT}(\lambda/\mu, \nu/\kappa)$ then $\chi^S(T)$ has shape $\lambda/\mu$ and a companion in $(\nu/\kappa)^\bullet$. The second part is easiest to do with the RSK algorithm, so we postpone the proof until Corollary 5.12 in section 5.3. In for example [10], the reversal was defined using the canonical tableau as $S$, so it might seem like we have made a more general definition here. This is not the case, however, as we shall see later, since it turns out that the reversal does not depend on $S$. This independence implies that the shape stay the same and is another indication that the reversal is the proper way to generalise the Schützenberger involution. In order to prove this we need some more tools.

### 4.3. Jeu de taquin equivalence and dual equivalence

When studying the jdt map, a natural question to ask is which tableaux will be mapped to the same partition shaped tableau.

**Definition 4.6.** $T_1$ and $T_2$ are *jeu de taquin equivalent* if $\mathrm{jdt}(T_1) = \mathrm{jdt}(T_2)$.

A very important property is that all jeu de taquin equivalent tableaux have companions of the same shape. For a proof see [16].

**Proposition 4.7.** *Let* $S, T \in \mathrm{YT}$ *and* $\mathrm{jdt}(S) = \mathrm{jdt}(T)$, *then*

$$T \text{ is } \nu/\kappa\text{-}dominant \iff S \text{ is } \nu/\kappa\text{-}dominant.$$

**Corollary 4.8.** *Let* $T \in \mathrm{YT}$, *then*

   (i) $T \in \mathrm{YT}_{\ulcorner} \iff \mathrm{jdt}(T)$ *is canonical,*
   (ii) $T \in \mathrm{YT}_{\lrcorner} \iff \mathrm{jdt}(T)$ *is anti-canonical.*

*Proof.* For the first statement, let $S = \mathrm{jdt}(T)$ and use the previous proposition. The second equivalence follows from the first by substituting $T$ for $T^{\bullet}$ and using Lemma 4.3. $\qquad\square$

Jeu de taquin equivalence measures if two tableaux become the same after the jeu de taquin slides. The dual equivalence, on the other hand, checks if two tableaux act the same way during a tableau slide sequence. Or, in other words, from the point of view of the tableau $T$, it does not matter if it was switched with $S$ or any tableau dual equivalent to $S$. The result would be the same. Dual equivalence was first defined by Haiman [5], but the presentation given here follows closely the one in van Leeuwen [16], and proofs for the propositions below can be found there.

**Definition 4.9.** The tableaux $S_1$ and $S_2$ of the same shape, are *dual equivalent* if for every $T$ with $X(S_1, T) = (T_1', S_1')$ and $X(S_2, T) = (T_2', S_2')$ one has $T_1' = T_2'$.

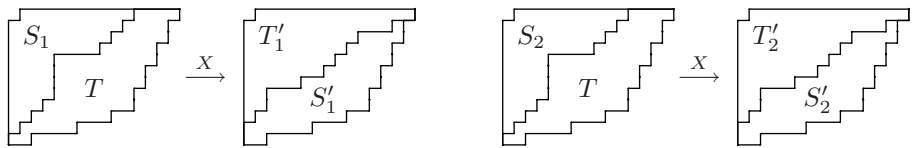

FIGURE 15. Illustration of the dual equivalence.

At a first glance it might appear that this definition is not very useful, since it would seem improbable that this property could hold for all different $T$. But, surprisingly, it is not at all hard to find examples of dual equivalent tableaux. The following propositions give an explanation.

First, a tableau switching in the opposite direction also behaves in a similar manner.

**Proposition 4.10.** *If $X(S_1, T) = (T', S_1')$ and $X(S_2, T) = (T', S_2')$, then $S_1$ and $S_2$ are dual equivalent iff $S_1'$ and $S_2'$ are dual equivalent.*

This means that given any pair of dual equivalent tableaux, we can create more pairs by switching both of them with the same tableau. Next, a large source of dual equivalences comes from the following proposition on partition shaped tableaux.

**Proposition 4.11.** *All tableaux in $\mathrm{YT}^{\Gamma}(\lambda)$ are dual equivalent.*

This proposition is in fact equivalent to Proposition 3.1, since the choice of order of the jeu de taquin slides when calculating $\mathrm{jdt}(T)$ corresponds to the choice of $S$ in the definition of dual equivalence.

The jeu de taquin and the dual equivalences are complementary to each other in the sense that if two tableaux are both jeu de taquin and dual equivalent they must be equal. We will see later that this corresponds to the Robinson-Schensted-Knuth map, where jeu de taquin equivalences means that the $P$ tableaux are the same and dual equivalence that the $Q$ tableaux are the same.

**Proposition 4.12.** *If $T_1$ and $T_2$ are both jeu de taquin equivalent and dual equivalent, then $T_1 = T_2$.*

**Proposition 4.13.** *If $T_1$ and $T_2$ are jeu de taquin equivalent, then $\tau(T_1)$ and $\tau(T_2)$ are dual equivalent.*

4.4. **Uniqueness of the reversal**

**Lemma 4.14.** *$\chi^S(T)$ does not depend on $S$.*

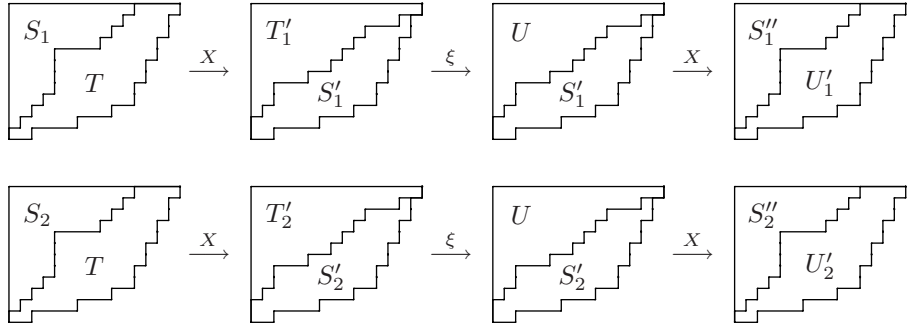

Figure 16. Reversal of $T$ using two different $S$ tableaux.

*Proof.* Let $T \in \mathrm{YT}(\lambda/\mu)$ and $S_1, S_2 \in \mathrm{YT}^{\Gamma}(\mu)$ be arbitrary. We need to show that $\chi^{S_1}(T) = \chi^{S_2}(T)$. Now put $(T_1', S_1') = X(S_1, T)$ and $(T_2', S_2') = X(S_2, T)$ as in Figure 16. $S_1$ and $S_2$ are dual equivalent by Proposition 4.11, so $T_1' = T_2'$ and by Proposition 4.10, $S_1'$ and $S_2'$ are dual equivalent. Let $U = \xi(T_1') = \xi(T_2')$. Now since $S_1'$ and $S_2'$ are dual equivalent, we have $\chi^{S_1}(T) = U_1' = U_2' = \chi^{S_1}(T)$, where $(S_1'', U_1') = X(U, S_1')$ and $(S_2'', U_1') = X(U, S_2')$. □

Note also that since $U$ and $T_i'$ are dual equivalent, we have $S_1'' = S_1$ and $S_2'' = S_2$. By this lemma, we can define the reversal independently of $S$, by using e.g. the canonical tableau as $S$, $\chi(T) = \chi^{\mathrm{can}(\mu)}(T)$.

## 5. **The fundamental symmetry map**

Next we define the fundamental symmetry map, by using Corollary 4.8 and Lemma 3.4, which tells us that if $X(S,T) = (T', S')$, where $S \in \mathrm{YT}^{\Gamma}$, then

$$T \in \mathrm{YT}_{\Gamma} \quad \Longleftrightarrow \quad T' \text{ is canonical.}$$

Using this fact, we can define an involution $\rho : \mathrm{YT}_{\Gamma}(\lambda/\mu, \nu) \to \mathrm{YT}_{\Gamma}(\lambda/\nu, \mu)$, called the *fundamental symmetry map*:

$$X(\mathrm{can}(\mu), T) = (\mathrm{can}(\nu), \rho(T)).$$

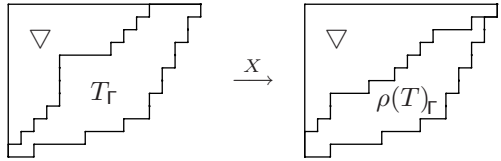

FIGURE 17. Definition of the fundamental symmetry map. A triangle with the tip pointing downwards stands for a canonical tableau of appropriate shape and the indices on $T$ and $\rho(T)$ indicate that are Littlewood-Richardson.

By the second part of Corollary 4.8, we have that if $X(S,T) = (T', S')$, then $T \in \mathrm{YT}_{\lrcorner}$ iff $T'$ is anti-canonical. Hence we can extend the definition of $\rho$ to allow tableaux of type $\mathrm{YT}_{\lrcorner}$ as well, $\rho : \mathrm{YT}_{\lrcorner}(\lambda/\mu, \nu^{\bullet}) \to \mathrm{YT}_{\lrcorner}(\lambda/\nu, \mu^{\bullet})$, such that

$$X(\xi(\mathrm{can}(\mu)), T) = (\xi(\mathrm{can}(\nu)), \rho(T)),$$

where $\xi(\mathrm{can}(\mu)) = \mathrm{can}(\mu^{\bullet})^{\bullet}$ is the anti-canonical tableau of shape $\mu$.

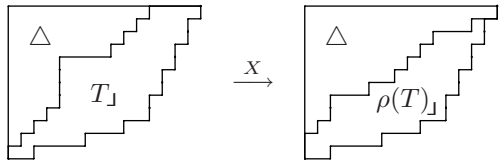

FIGURE 18. Definition of the fundamental symmetry map for anti-Littlewood-Richardson tableaux. An upside down triangles stands for a anti-canonical tableau of appropriate shape and the indices on $T$ and $\rho(T)$ indicate that are anti-Littlewood-Richardson.

Note that there is an ambiguity when $T \in \mathrm{YT}_{\Gamma} \cap \mathrm{YT}_{\lrcorner}$, so in order to make it a well defined involution we let the definition set be $\mathrm{YT}_{\Gamma} \uplus \mathrm{YT}_{\lrcorner}$, which ensures that we can differentiate between the two cases.

### 5.1. **A dual fundamental symmetry map**

It is possible to define a 'dual' version of the fundamental symmetry map, $\rho'$, by exchanging the use of the canonical/anti-canonical tableaux when applying the tableau switch.

$$X(\xi(\mathrm{can}(\mu)), T) = (\mathrm{can}(\nu), \rho'(T)), \quad \text{if } T \in \mathrm{YT}_\ulcorner,$$
$$X(\mathrm{can}(\mu), T) = (\xi(\mathrm{can}(\nu)), \rho'(T)), \quad \text{if } T \in \mathrm{YT}_\lrcorner.$$

**Lemma 5.1.** *Let* $T \in \mathrm{YT}_\ulcorner \uplus \mathrm{YT}_\lrcorner$. *Then* $\rho' \circ \rho(T) = \rho \circ \rho'(T) = \chi(T)$ *and hence* $\chi$ *commutes with both* $\rho$ *and* $\rho'$.

*Proof.* The composition $\rho' \circ \rho$ is exactly the reversal applied on the canonical tableau and when changing the order to $\rho \circ \rho'$ one gets the reversal applied to the anti-canonical tableau. Hence by the uniqueness of the reversal, they are the same. □

### 5.2. **The $\Omega$-involution**

The last involution we define is the $\Omega$-involution, which is defined by three consecutive tableau switchings. We show in Lemma 5.15 that this involution is the same as the composition of the rotation involution and the reversal.

**Definition 5.2.** Let $\Omega : \mathrm{YT}(\lambda/\mu, \nu/\kappa) \to \mathrm{YT}((\lambda/\mu)^\bullet, \nu/\kappa)$ be defined by $\Omega(T) = T''$, where $X(\triangledown, T) = (T', S)$, $X(S, \triangledown) = (U, \triangledown)$ and $X(T', U) = (\triangledown, T'')$. See Figure 19 for an illustration.

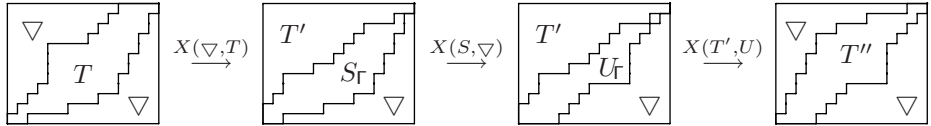

Figure 19. Illustration of the $\Omega$ involution.

Since tableau switching in involutive it is clear that $\Omega$ is an involution. In terms of the Bender-Knuth transforms, $\Omega$ can be written

$$\Omega = t_{l,m}\, t_{k,m}^{(l)}\, t_{k,l},$$

due to the fact that the middle tableau switching, which goes in the lower right direction, can be written as a conjugation of tableau switching by the rotation involution, $t_{k,m}^{(l)}(T) = (t_{m,k}(T^\bullet))^\bullet$. If $T$ is a Littlewood-Richardson tableau, then all three tableau switchings are fundamental symmetry maps, so we get

**Lemma 5.3.** *If* $T \in \mathrm{YT}_\ulcorner$, *then*

$$\Omega(T) = \rho(\rho(\rho(T)^\bullet)^\bullet).$$

5.3. **An RSK correspondence for tableaux**

The jeu de taquin equivalence and dual equivalence have a very direct interpretation by the RSK algorithm. To make this clearer, we define an alternative version of the RSK algorithm, which is defined on tableaux instead of matrices.

**Definition 5.4.** We define $\mathrm{rsk} : \mathrm{YT} \to \bigcup_\lambda \mathrm{YT}^\ulcorner(\lambda) \times \mathrm{YT}^\ulcorner(\lambda)$ as

$$\mathrm{rsk}(T) = \big(\mathrm{jdt}(T), \ \mathrm{jdt}(\tau(T))\big).$$

Shifting rows in a tableau horizontally must give a jeu de taquin equivalent tableau, since a one step shift of a row that does not destroy the semistandardness is precisely a jeu de taquin slide along that row. Induction over the number of shifts gives:

**Lemma 5.5.** Let $T_1, T_2 \in \mathrm{YT}$. Then

$$\mathrm{rec}(T_1) = \mathrm{rec}(T_2) \quad \Longrightarrow \quad \mathrm{jdt}(T_1) = \mathrm{jdt}(T_2).$$

This means that rsk is well defined and also that it is closely related to the classical algorithm. As a direct consequence of Lemma 4.3, we get

**Lemma 5.6.** Let $T \in \mathrm{YT}$, $M = \mathrm{rec}(T)^r$ and $\mathrm{rsk}(T) = (P, Q)$. Then

$$\mathrm{RSK}(M) = (P, \xi(Q)).$$

**Lemma 5.7.** Let $T_1, T_2 \in \mathrm{YT}(\lambda/\mu, \nu/\kappa)$, $\mathrm{rsk}(T_1) = (P_1, Q_1)$ and $\mathrm{rsk}(T_2) = (P_2, Q_2)$. Then

$$
\begin{aligned}
P_1 = P_2 &\quad \Longleftrightarrow \quad T_1 \text{ and } T_2 \text{ are jeu de taquin equivalent,} \\
Q_1 = Q_2 &\quad \Longleftrightarrow \quad T_1 \text{ and } T_2 \text{ are dual equivalent,} \\
T_1 = T_2 &\quad \Longleftrightarrow \quad P_1 = P_2 \text{ and } Q_1 = Q_2.
\end{aligned}
$$

*Proof.* The first statement follows directly from the definition. The second was proved by van Leeuwen in [16], and in [15] for the picture version. The last equivalence is just a combination of the two previous, together with Proposition 4.12. $\square$

**Corollary 5.8.** Let $T_1, T_2 \in \mathrm{YT}(\lambda/\mu, \nu/\kappa)$. Then

$$T_1 = T_2 \quad \Longleftrightarrow \quad \mathrm{rsk}(T_1) = \mathrm{rsk}(T_2).$$

**Lemma 5.9.** Suppose $\mathrm{rsk}(T) = (P, Q)$. Then $T \in \mathrm{YT}(\lambda/\mu, \nu/\kappa)$ if and only if $P \in \mathrm{YT}^\ulcorner(\sigma, \nu/\kappa)$ and $Q \in \mathrm{YT}^\ulcorner(\sigma, \lambda/\mu)$, for some shape $\sigma$.

*Proof.* This follows directly from Proposition 4.7, since

$$
\begin{aligned}
P = \mathrm{jdt}(T) \text{ is } \nu/\kappa\text{-dominant} &\quad \Longleftrightarrow \quad T \text{ is } \nu/\kappa\text{-dominant} \\
Q = \mathrm{jdt}(\tau(T)) \text{ is } \lambda/\mu\text{-dominant} &\quad \Longleftrightarrow \quad \tau(T) \text{ is } \lambda/\mu\text{-dominant.}
\end{aligned}
$$

$\square$

**Lemma 5.10.**

$$\mathrm{rsk}(T) = (P, Q) \quad \Longleftrightarrow \quad X(\bigtriangledown, T) = (P, \rho(\tau(Q))).$$

*Proof.* Let $X(\bigtriangledown, T) = (P', S)$. Firstly we have from the definitions that $P = \mathrm{jdt}(T) = P'$. Also, $S$ only depends on the dual equivalence class of $T$, since if $T_1$ and $T_2$ are dual equivalent, with $\mathrm{rsk}(T_1) = (P_1, Q)$ and $\mathrm{rsk}(T_2) = (P_2, Q)$, we have by the definition of dual equivalence, $X(\bigtriangledown, T_1) = (P_1, S_1)$ and $X(\bigtriangledown, T_2) = (P_2, S_2)$, where $S_1 = S_2$. Thus, since $T$ and $\tau(Q) \in \mathrm{YT}_\Gamma$ are dual equivalent by Proposition 4.13,

$$(\bigtriangledown, S(Q)) = X(\bigtriangledown, \tau(Q)) = (\bigtriangledown, \rho(\tau(Q))).$$

$\square$

This lemma is useful for computing the inverse of rsk. Given $P, Q \in \mathrm{YT}^\Gamma(\lambda)$,

$$(\bigtriangledown, \mathrm{rsk}^{-1}(P, Q)) = X(P, \rho(\tau(Q))).$$

Also the inverse of the classical RSK correspondence can be computed using this, in the light of Lemma 5.6,

$$\mathrm{RSK}^{-1}(P, Q) = \mathrm{rec}(\mathrm{rsk}^{-1}(P, \xi(Q)))^r.$$

**Lemma 5.11.** *Suppose that* $\mathrm{rsk}(T) = (P, Q)$, *then*

$$\mathrm{rsk}(T^\bullet) = (\xi(P), \xi(Q)),$$
$$\mathrm{rsk}(\tau(T)) = (Q, P),$$
$$\mathrm{rsk}(\chi(T)) = (\xi(P), Q).$$

*Proof.* The first two relations follows directly from Lemma 5.6 and Theorem 4.2. For the last identity, let $\mathrm{rsk}(\chi(T)) = (\tilde{P}, \tilde{Q})$. Then, by Lemma 5.10, $X(\bigtriangledown, T) = (P, \rho(\tau(Q)))$ and $X(\bigtriangledown, \chi(T)) = (\tilde{P}, \rho(\tau(\tilde{Q})))$, so that, by the definition of $\chi$, $\tilde{P} = \xi(P)$ and $\tilde{Q} = Q$. $\square$

**Corollary 5.12.** *If* $T \in \mathrm{YT}(\lambda/\mu, \nu/\kappa)$, *then* $\chi(T)$ *has a companion in* $(\nu/\kappa)^\bullet$.

*Proof.* From the previous lemma we have that $\mathrm{jdt}(\chi(T)) = \mathrm{jdt}(T^\bullet)$. Therefore $\chi(T)$ is $(\nu/\kappa)^\bullet$-dominant by Proposition 4.7 since $T^\bullet$ is. $\square$

### 5.4. **Commutation relations for the reversal**

Using the rsk correspondence, we can find the commutation relations for the reversal.

**Corollary 5.13.** *Let* $T \in \mathrm{YT}(\lambda/\mu, \nu/\kappa)$. *Then*

$$\chi(T^\bullet) = \chi(T)^\bullet.$$

*Proof.* If $\mathrm{rsk}(T) = (P, Q)$, both sides of the equation are mapped to $(P, \xi(Q))$ by rsk, hence, by Lemmas 5.7 and 5.11, they are equal. $\square$

**Lemma 5.14.** *Let* $T \in \mathrm{YT}(\lambda/\mu, \nu/\kappa)$. *Then*

$$\chi(\tau(T)) = \tau(\chi(T^\bullet)), \quad \textit{i.e.,} \quad \chi \circ \tau \circ \chi \circ \tau(T) = T^\bullet.$$

*Proof.* Again, both sides have the same rsk-tableaux, namely $(\xi(Q), P)$. $\square$

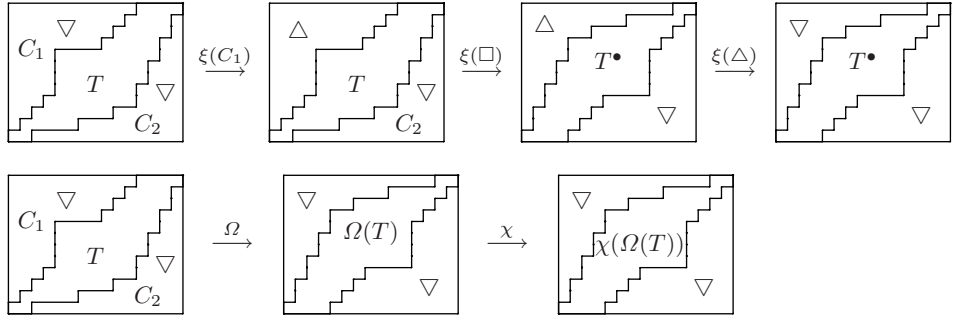

Figure 20. Illustration of the $\Omega$-lemma.

### 5.5. The $\Omega$-lemma

**Lemma 5.15.** *Let* $T \in \mathrm{YT}(\lambda/\mu)$. *Then* $\Omega(T) = \chi(T^\bullet)$ *(which in turn is equal to* $\chi(T)^{\bullet}$).

*Proof.* We prove that the two rows in Figure 20 actually result in the same tableaux by showing that both are applying the same sequence of Bender-Knuth transformations. By Lemma 4.3(iv) we have that if $T$ has rectangular shape, then $T^\bullet = \mathrm{jdt}(T^\bullet) = \xi(T)$, which explains the middle map in the first row.

Let $C_1$ and $C_2$ be canonical tableaux as in the figure, $k = \ell(C_1)$, $l = \ell(\tau(T))$ and $m = \ell(C_2)$. Evaluating the first row in terms of Bender-Knuth transformations and simplifying, using Lemma 3.3, gives

$$
\begin{aligned}
z_m \, z_{k+l+m} \, z_k &= z_m \left( z_m \, t_{k+l,m} \, z_{k+l} \right) z_k && \text{(Lemma 3.3(v))} \\
&= t_{l,m} \, t_{k,m}^{(l)} (z_l \, t_{k,l} \, z_k) \, z_k && \text{(Lemma 3.3(iii),(v) and (vi))} \\
&= t_{l,m} \, z_l \, t_{k,m}^{(l)} \, t_{k,l} && \text{(Lemma 3.3(ii) and (iii))} \\
&= (t_{l,m} \, z_l \, t_{m,l}) \, (t_{l,m} \, t_{k,m}^{(l)} \, t_{k,l}), && \text{(Lemma 3.3(iv))}
\end{aligned}
$$

where the last expression is the Bender-Knuth transformations of the second row of Figure 20. So we get $T^\bullet = \chi(\Omega(T))$, i.e., $\Omega(T) = \chi(T^\bullet)$. $\qquad\square$

### 5.6. The commutative diagram

Finally, we can now create a commutative diagram, which summarises most of the previous lemmas for (anti-)Littlewood-Robinson tableaux and partition shaped tableaux. The most striking part is the threefold symmetry, which reflects the fact that it requires three partitions to specify the shape and the weight of a Littlewood-Robinson tableau. Similar type of threefold symmetry has also been studied by Berenstein and Zelevinsky [3], using BZ-triangles, and by Knutson and Tao [7], using hives and honeycombs.

**Theorem 5.16.** *The diagram in Figure 21 is commutative.*

*Proof.* From Lemmas 2.7, 5.13 and 5.14 we know that the three rectangular parts are commuting diagram. Also, by Lemmas 5.15, 5.1 and 5.3, we have, if $r$ denotes the map $T \mapsto T^{\bullet}$,

$$\chi \circ r \circ \rho \circ \chi \circ r \circ \rho \circ \chi \circ r \circ \rho$$
$$= r \circ \chi \circ \rho \circ r \circ \rho \circ r \circ \rho \quad (\chi \text{ commutes with both } r \text{ and } \rho)$$
$$= \Omega \circ \Omega \quad\quad\quad\quad (\Omega = r \circ \chi = \rho \circ r \circ \rho \circ r \circ \rho)$$
$$= I.$$

$\square$

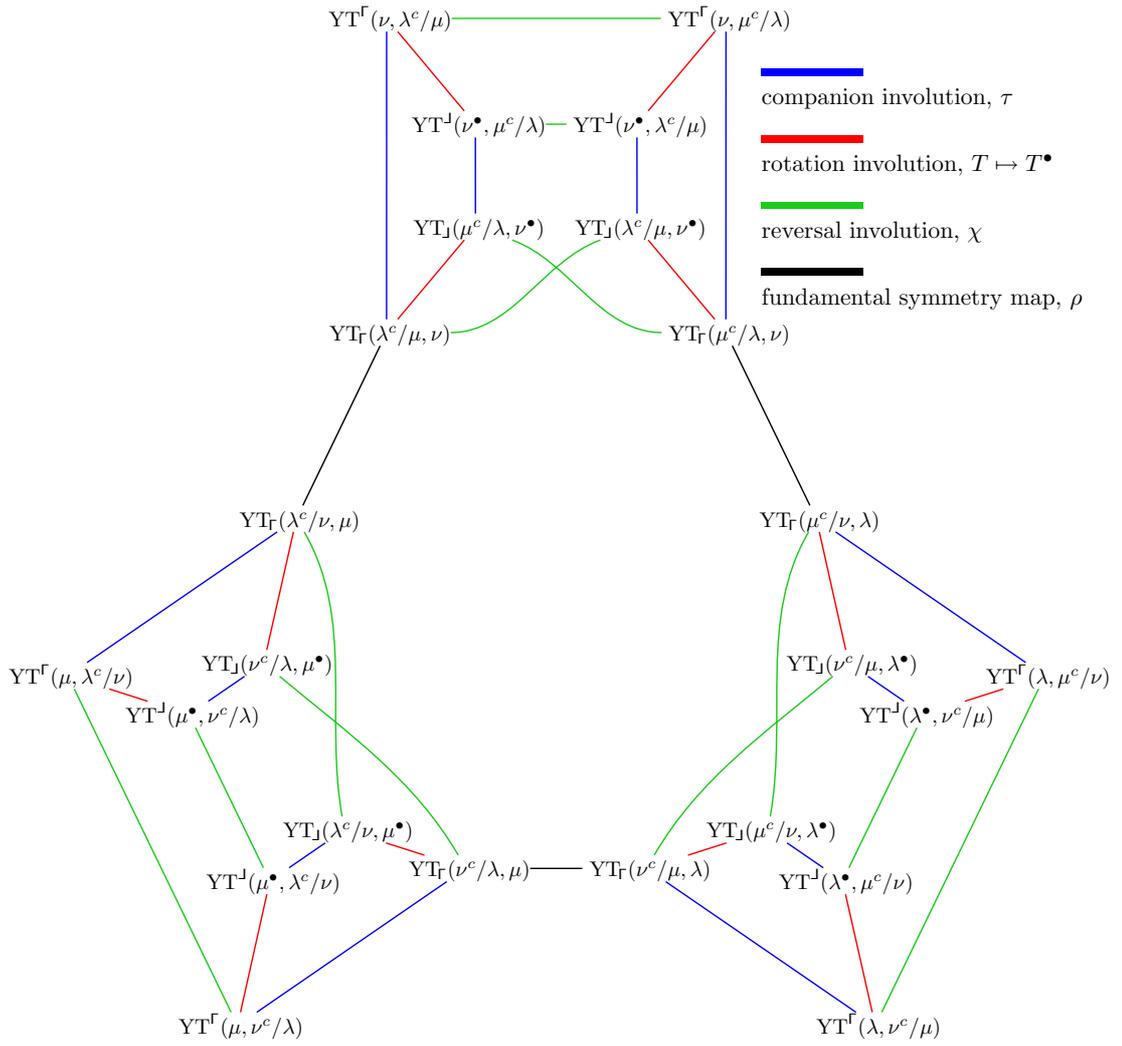

Figure 21. Commutative diagram.

Linköpings Universitet, 581 83 Linköping, Sweden
*E-mail address*: `erouc@mai.liu.se`